\documentclass[12pt, twoside]{article}
\usepackage{amsmath,amsthm,amssymb}
\usepackage{times}
\usepackage{enumerate}

\theoremstyle{definition}
\newtheorem{thm}{Theorem}[section]

\newtheorem*{xrem}{Remark}

\numberwithin{equation}{section}

\newcommand{\subjclass}[1]{\bigskip\noindent\emph{2010 Mathematics Subject Classification:}\enspace#1}
\newcommand{\keywords}[1]{\noindent\emph{Keywords:}\enspace#1}
\newtheorem{problem}[thm]{Problem}
\newtheorem{algorithm}[thm]{Algorithm}

\newtheorem{conjecture}[thm]{Conjecture}

\usepackage{color}
\definecolor{gray}{gray}{0.5}
\definecolor{newblue}{RGB}{94,89,144}
\definecolor{newblue2}{cmyk}{1,0.6,0,0.06}
\usepackage[pdfborder={0 0 0},
  colorlinks=true,
  citecolor=newblue,
  linkcolor=newblue2,
  urlcolor=newblue]{hyperref}

\usepackage[numbers,longnamesfirst]{natbib}
\bibpunct{(}{)}{;}{a}{}{,}

\usepackage{subcaption} % for multiple floats
\usepackage{verbatim}

\newcommand{\eps}{\varepsilon}
\newcommand{\R}{\mathbb{R}}
\newcommand{\abs}[1]{\left| #1 \right|}
\newcommand{\norm}[1]{\left\| #1 \right\|}
\newcommand{\dd}{\, \mathrm{d}}

\newcommand{\T}{\mathcal{T}}
\newcommand{\N}{\mathcal{N}}
\renewcommand{\vec}[1]{\ensuremath{\boldsymbol{#1}}}

\newcommand{\E}{\mathcal{E}}

\newcommand{\Epf}{\mathcal{E}_{\mathrm{SEG}}}

\usepackage{graphicx}
\usepackage{color}
\usepackage{hyperref}
\usepackage{array}
\begin{document}

%%%%% To ease editing, add:

\baselineskip=17pt

%%%%%%%%%%%%%%%%

\title{A computational approach to an optimal partition problem on surfaces}

\author{Charles M. Elliott\\
Mathematics Institute, Zeeman Building, University of Warwick. CV4 7AL. UK.\\
C.M.Elliott@warwick.ac.uk\\
Thomas Ranner\\
School of Computing, EC Stoner Building, University of Leeds. LS2 9JT. UK.\\
T.Ranner@leeds.ac.uk}

\date{}

\maketitle

%%%%%%%%

\begin{abstract}
  We explore an optimal partition problem on surfaces using a computational approach. The problem is to minimise the sum of the first Dirichlet Laplace--Beltrami operator eigenvalues over a given number of partitions of a surface. We consider a method based on eigenfunction segregation and perform calculations using modern high performance computing techniques. We first test the accuracy of the method in the case of three partitions on the sphere then explore the problem for higher numbers of partitions and on other surfaces.

% 49Q10 geometric minimisation not minimal surfaces 
% 49R50 Variational methods for eigenvalues of operators
% 35R01 pdes on manifolds
% 65M60 finite elements


\end{abstract}

\section{Introduction}

In this paper, we use the surface finite element method to tackle an eigenvalue optimal partition problem for  $n$-dimensional 
hypersurfaces in $\mathbb R^{n+1}$. Our computations are restricted to $n=2$.
We denote by $\Gamma$ a closed, smooth, connected $n$-dimensional hypersurface embedded in $\R^{n+1}$. 
For a given positive integer $m$, we say that $\{ \Gamma_i \}_{i=1}^m$ is an $m$-partition 
of $\Gamma$ if $\Gamma_i \subset \Gamma$ for 
$i = 1,\ldots,m$, $\Gamma_i \cap \Gamma_j = \emptyset$ for $i,j = 1,\ldots,m$ 
with $i \neq j$ and $\bigcup_{i = 1,\ldots,m} \overline\Gamma_i = \Gamma$.

\begin{problem}
  \label{pb:partition}
  Given a positive integer $m$ and a smooth surface $\Gamma$, divide $\Gamma$ into 
  an $m$-partition $\{ \Gamma_i \}_{i=1}^m$ to minimise the energy:
  \begin{equation}
    \label{eq:1}
    \E( \{ \Gamma_i \}_{i=1}^m ) = \sum_{i=1}^m \lambda_1 ( \Gamma_i ),
  \end{equation}
  where $\lambda_1( \Gamma_i )$ is the first eigenvalue of the Dirichlet Laplace-Beltrami operator over $\Gamma_i$.
\end{problem}

 This is a generalisation of a similar problem considered in various formulations  over a Cartesian 
 domain $\Omega$ with appropriate boundary conditions. 
The flat problem was studied in the context of shape optimisation in the 1990's by \citet{ButMaz93,Sve93,BucZol95,BucButHen98}. 
A key challenge is how to define an appropriate space of admissible partitions and how to equip this space 
with a topology so that one can define an absolute minimiser. By restricting to quasi-open sets, \citet{BucButHen98} 
show existence of a optimal partition as a consequence of a more general result. Quasi-open sets are sets which are 
close to open sets in the sense that given a quasi-open set  there is an open set 
such that their symmetric difference has arbitrarily small capacity \citep{CafLin07}. 
Formally speaking, these are a class of general sets which can be 
used to define a weak form of elliptic equations. For example, all open sets are quasi-open. The set
 $\mathcal A(\Omega)$ of 
quasi-open sets in a domain $\Omega$ can be equipped with a notion of weak convergence  by defining that a sequence of quasi-open 
sets $\{ A_n \}$ weakly converges to $A \in \mathcal A(\Omega)$ if $\eta_{A_n} \to \eta_A$ weakly in $H^1(\Omega)$ and $A=\{\eta_{A}>0\}$  
where 
$\eta_\omega \in H^1(\Omega)$ is the extension to $\Omega$ by zero of 
the unique weak solution of
\begin{equation*}
  - \Delta \eta_\omega = 1 \quad \mbox{ in } \omega 
  \qquad \mbox{ and } \qquad 
  \eta_\omega = 0 \mbox{ on } \partial \omega.
\end{equation*}
Using these notions  it is possible to establish that the spectral functional is lower semi-continuous with 
respect to weak convergence in $\mathcal A(\Omega)$ and existence of an $m$-partition into quasi-open sets follows from the 
direct method of the calculus of variations \citep{CafLin07}.

An alternative method is based on  using the eigenfunctions to partition 
the domain using an approach formulated by \citet{CafLin07}.
The energy \eqref{eq:1} is transformed into a functional form as a constrained Dirichlet energy:
\begin{problem}
  Given a positive integer $m$ and a smooth surface 
  $\Gamma$, find $\vec{u} = (u_1, \ldots, u_m) \in H^1( \Gamma, \Xi )$ with $\norm{ u_i }_{L^2(\Gamma)} = 1$ for $i = 1, \ldots, m$, to minimise
  \begin{equation}
    \label{eq:7}
    \Epf^0( \vec{u} ) = \sum_{i=1}^m \int_{\Gamma} \abs{ \nabla_\Gamma u_i }^2 \dd \sigma,
  \end{equation}
  where $\Xi \subset \R^m$ is the singular set
  \begin{equation*}
    \Xi = \left\{ \vec{y} = ( y_1, \ldots, y_m ) \in \R^m : \sum_{i=1}^m \sum_{i \neq j} y_i^2 y_j^2 = 0 \mbox{ and } y_{i}\ge 0, i=1,2,...m \right\}.
  \end{equation*}
\end{problem}
It was shown by  \citet{CafLin07} that, when $\Gamma$ is a Cartesian domain in $\R^n$, \eqref{eq:7} is equivalent to \eqref{eq:1} when we restrict to $m$-partitions of $\Gamma$ 
in which $\Gamma_{i}$ are quasi-open sets.
The proof can be adapted to the surface case also.
Let $\{ \Gamma_i \}_{i=1}^m$ be a minimiser of \eqref{eq:1} consisting of quasi-open sets, 
then if $u_i$ is the first eigenfunction of the Dirichlet Laplace--Beltrami operator over $\Gamma_i$, 
for $i=1,\ldots,m$, the vector quantity $\vec{u} = ( u_1, \ldots, u_m )$ is a minimiser of \eqref{eq:7}. 
Conversely, let the function $\vec{u} = ( u_1, \ldots, u_m ) \in H^1( \Gamma, \Xi )$ be a minimiser of \eqref{eq:7}, 
then  setting $\Gamma_i = \{ u_i > 0 \}$, for $i = 1,\ldots,m$, the collection of quasi-open sets  $\{ \Gamma_i \}_{i=1}^m$ is an $m$-partition of $\Gamma$ which is a minimiser of \eqref{eq:1} and
\begin{equation*}
  \lambda_1( \Gamma_ i ) = \int_\Gamma \abs{ \nabla_\Gamma u_i }^2 \dd \sigma \quad \mbox{ for } i = 1,\ldots,m.
\end{equation*}
The authors \citet{CafLin07} use this formulation to show existence of minimisers and regularity of the interface between partitions.

Other works by \citet{ConTerVer02,ConTerVer03} and \citet{CafLin07,CafLin08} have focused on regularity and more 
qualitative aspects of the problem for a Cartesian domain. Conti, Terracini and Verzini derive optimality conditions, such as 
the gradient of eigenfunctions should match at partition boundaries, and also that the partition consists of open sets. Caffarelli 
and Lin obtain regularity results, such as $C^{1,\alpha}$-smoothness of the partition boundaries away from a set of codimension two, 
and also an estimate of the behaviour in the limit of large $m$. In particular, they prove that the optimal energy is bounded above and below by a constant times the $m$-th eigenvalue on $\Gamma$ and conjecture that for large $m$ the optimal partition 
will be asymptotically close to a hexagonal tiling in the case of a planar domain. The problem can be seen as a strong competition 
limit of segregating species either in Bose-Einstein condensate \citep{ChaLinLin04}, population 
dynamics \citep{ConTerVer05,ConTerVer05a} or materials science \citep{Che02} in curved geometries.

Numerical studies of this type of problem have so far been limited to the planar case. We mention in particular the study of 
\cite{ChaLinLin04} and  some special algorithms in the case of small $m$ given by \cite{BozAra13} and \cite{Boz09}. Also \citet{BouBucOud10} considered the problem for large values of $m$ using a fictitious domain approach. This 
problem has also been considered on graphs \citep{CoiLaf06,OstWhiOud14-pp} with applications in big data segmentation. Finally, we 
mention the study which will be the basis of our work in the paper: an eigenfunction segregation 
approach \citep{DuLin09}. We will describe the algorithm in more detail in the following.

\begin{figure}[tb]
  \centering
  \includegraphics[width=0.3\textwidth]{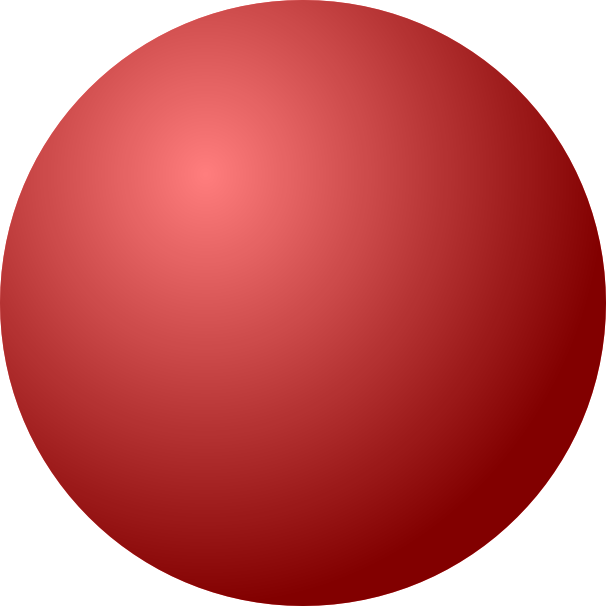}
  \quad
  \includegraphics[width=0.3\textwidth]{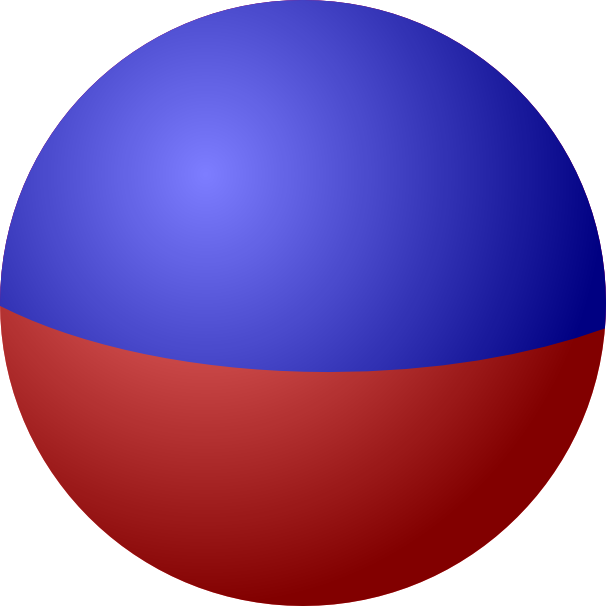}
  \quad
  \includegraphics[width=0.3\textwidth]{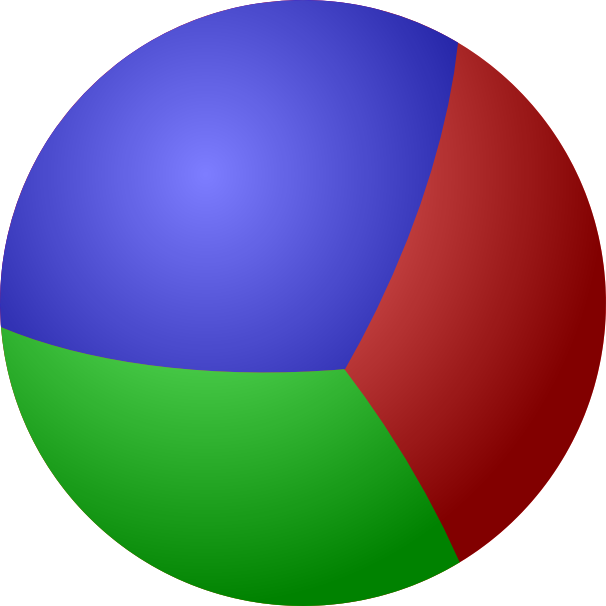}
  
  \caption{Plots of solution of Problem~\ref{pb:partition} when $\Gamma$ is a sphere, $m = 1$ (left), $m=2$ (center) and $m=3$ (right, Bishop's conjecture) \citep{HelOstTer10}.}
  \label{fig:known-sphere}
\end{figure}

The curved hypersurface problem has been studied analytically in the case that $\Gamma$ is a sphere.
For $m=1$, the result is clear and for $m=2$ the solution is two hemispheres leading to total energy $2$.
The case $m=3$ on the sphere leads to the Bishop conjecture \citep{Bis92}.
\begin{conjecture}
  \label{conj:bishop}
  The minimal $3$-partition for Problem~\ref{pb:partition}, with $\Gamma = $ sphere, corresponds to the Y-partition whose boundary is given, up to a fixed rotation, by the intersection of $\Gamma$ with the three half planes defined in polar coordinates by $\phi = 0, \frac{2\pi}{3}, \frac{-2\pi}{3}$ (see Figure~\ref{fig:known-sphere} and Section~\ref{Ypartition}).
\end{conjecture}

A similar problem to Problem~\ref{pb:partition} has been considered by exchanging the sum in \eqref{eq:1} to an $\ell^p$-norm for $p \in [1,\infty]$.
\begin{problem}
  \label{pb:partition-p}
  Given a positive integer $m$ and a smooth surface $\Gamma$, divide $\Gamma$ into an $m$-partition $\{ \Gamma_i \}_{i=1}^m$ to minimise the energy
  \begin{equation}
    \label{eq:partition-p}
    \E_p( \{ \Gamma_i \}_{i=1}^m ) =
    \begin{cases}
      \left( \frac{1}{m} \sum_{i=1}^m \lambda_1( \Gamma_i )^{p} \right)^{\frac{1}{p}}
      & \quad p \in [1,\infty) \\
      \max_{i=1,\ldots,m} \lambda_1( \Gamma_i )
      & \quad p = \infty.
    \end{cases}
  \end{equation}
\end{problem}

The differences between this more general problem and the case $p=1$ have been studied by \citet{HelHof09} in the case of Cartesian domains.
In particular they show a monotonicity formula for optimal partitions: Denoting by $\mathcal{P}_p$ the optimal partition for the energy $\E_p$, for $p \in [1,\infty]$, then we have
\begin{equation*}
  \E_p( \mathcal{P}_p ) \le \E_q( \mathcal{P}_q )
  \qquad
  \mbox{ if } p \le q.
\end{equation*}
It is well known that the optimal partition for the case $p = \infty$ is equi-spectral (equal $\lambda_1$ for each set in the partition).
This implies that if a partition is optimal for $p=1$ and is equi-spectral then it is optimal for $p=\infty$ \citep[Proposition~2.1]{HelHof09}.

The case $p = \infty$ has been studied on the sphere in the recent work 
of \cite{HelOstTer10}. They show the optimal partition is given by two hemispheres for the case $m=2$ and the
Y-partition for $m=3$; see Figure~\ref{fig:known-sphere} and Section \ref{Ypartition}.
The authors also conjecture that for $m=4$ the optimal partition is a spherical projection of a regular tetrahedron.
Furthermore, they show that for each $m$ there is 
an optimal partition which satisfies an equal angle condition which says that the boundary arcs that meet at a critical point do so with equal angles.
Computations for the $p=\infty$ case on a flat torus can be found in \cite{Len14-pp}.

We derive computational approaches using the surface finite element method \citep{Dzi88,DziEll07a} to find solutions to these 
problems. A review of computational techniques for partial differential equations on surfaces is given by \citet{DziEll13a}. Our methods 
will be one of the algorithms given by \citet{DuLin09} applied with the surface finite element method in order to explore Problem~\ref{pb:partition}.

We believe some of the techniques used in this paper, such as operator splitting and parallel computing, could be applied in a wide 
range of multiphase problems; for example \cite{GraKorSac14}. In these problems, one typically has a large system of reaction 
diffusion systems to solve with small parameter $\eps$ indicating an interfacial width. The small parameter $\eps$ acts with 
nonlinear terms to separate different phases. Our methods are designed to be transferable to this type of problem also. In contrast 
to many multiphase problems, the dynamic problem considered in this paper is based on non-local interface motion.

\subsection{Approximation approach}

One could try to directly compute the gradient flow of the energy $\Epf^0$ in \eqref{eq:7}; see \citet{May98} for 
analytic considerations of this approach. However, this would lead to equations which would be hard to discretise. 
We instead relax the constraint that $\vec{u}$ takes values in $\Xi$ by adding a penalty term to the energy functional 
following \citet{CafLin08}. In this way, we consider the extended energy functional:
\begin{equation*}
  \Epf^\eps( \vec{u}^\eps ) = \sum_{i=1}^m\frac{1}{2} \int_{\Gamma} \abs{ \nabla_\Gamma u_i^\eps }^2 \dd \sigma + \int_{\Gamma} F_\eps( \vec{u}^\eps ) \dd \sigma,
  \qquad
  F_\eps( \vec{u}^\eps ) = \frac{1}{\eps^2} \sum_{i=1}^m \sum_{\stackrel{j=1}{j \neq i}}^m (u_i^\eps)^2 (u_j^\eps)^2.
\end{equation*}

\begin{problem}
  \label{pb:phase-cts}
  Given a positive integer $m$, a smooth surface $\Gamma$ and $\eps > 0$, find $\vec{u}^\eps = (u_1^\eps, \ldots, u_m^\eps) \in H^1( \Gamma, \R^m )$ with $\norm{ u_i^\eps }_{L^2(\Gamma)} = 1$ for $i = 1, \ldots, m$, to minimise
  \begin{equation}
    \label{eq:phase-cts}
    \Epf^\eps( \vec{u}^\eps ) = \sum_{i=1}^m \frac{1}{2} \int_{\Gamma} \abs{ \nabla_\Gamma u_i^\eps }^2 \dd \sigma + \int_{\Gamma} F_\eps( \vec{u}^\eps ) \dd \sigma.
  \end{equation}
\end{problem}

We will now compute the gradient flow of this relaxed problem. We seek a time dependent 
function $\vec{u}^\eps \colon \Gamma \times \R_+ \to \R^m$ and $\lambda^\eps \colon \R_+ \to \R^m$ satisfying
\begin{subequations}
  \label{eq:phase-gf}
  \begin{align}
    \label{eq:9}
    \partial_t u_i^\eps & = \Delta_\Gamma  u_i^\eps + \lambda_i u_i^\eps - \frac{2}{\eps^2} \left( \sum_{j \neq i} (u_j^\eps)^2 \right) u_i^\eps && \mbox{ on } \Gamma \times \R_+, \mbox{ for } i = 1,\ldots,m, \\
    \vec{u}^\eps( \cdot, 0 ) & = \vec{u}^0 && \mbox{ on } \Gamma,
  \end{align}
\end{subequations}
subject to the constraint
\begin{equation}
  \label{eq:10}
  \int_\Gamma \abs{ u_i^\eps }^2 \dd \sigma = 1 \quad \mbox{ for } i = 1, \ldots, m.
\end{equation}
Here, we suppose that the initial condition partitions $\Gamma$ and has unit norm:
\begin{equation*}
  \vec{u}^0 \in H^1( \Gamma, \Xi ), \qquad \int_\Gamma \abs{ u^0_i }^2 \dd \sigma = 1 \quad \mbox{ for } i = 1, \ldots, m.
\end{equation*}
We remark that  $u^0_i \ge 0$ implies $u^\eps_i \ge 0$ for $i= 1,\ldots,m$.

This gradient flow problem was studied by \citet{CafLin09} for Cartesian geometries. The proofs can be easily transferred onto surfaces. We recall their results stated on surfaces:
\begin{equation*}
  \lambda_i^\eps(t) = \int_\Gamma \abs{ \nabla_\Gamma u_i^\eps }^2 + \frac{2}{\eps^2} \left( \sum_{j \neq i} (u_j^\eps)^2 \right) (u_i^\eps)^2 \dd \sigma,
\end{equation*}
and
\begin{equation*}
  \Epf^\eps( \vec{u}^\eps ) \le \sum_{i=1}^m \lambda_i^\eps( t ) = \Epf^\eps( \vec{u}^\eps ) + 2 \int_\Gamma F_\eps( \vec{u}^\eps ) \dd \sigma.
\end{equation*}
Furthermore, they show that $\Epf^\eps( \vec{u}^\eps )$ is a monotone decreasing function of time for $\vec{u}^\eps$ the solution of \eqref{eq:phase-gf}. This implies the existence of a unique global strong solution $\vec{u}^\eps \in L^\infty( \R_+, H^1( \Gamma, \R^m ) )$ for each $\eps > 0$. Finally, they give estimates of interest when considering the sharp interface limit: Denoting by $\bar{\vec{u}}^\eps$ the minimiser of the $\eps$-problem, for any $0 < t_1 < t_2$, we have
\begin{equation*}
  \int_{t_1}^{t_2} \int_\Gamma F_\eps( \bar{\vec{u}}^\eps ) \dd \sigma \dd t \to 0 \quad \mbox{ as } \eps \to 0,
\end{equation*}
and that the limit of minimising functions as $\eps \to 0$, $\bar{\vec{u}}^\eps$ converges strongly in $H^1( \Gamma \times \R_+)$ to a suitable weak solution of the constrained gradient flow of \eqref{eq:7}.
Further asymptotic analysis of the limit $\eps \to 0$ has been considered by \citet{DuZha11} and \citet{BerLinWei13}.

A key advantage of this approach is that we are trying to approximate smooth functions $\vec{u}^\eps$ in place of the domains $\Gamma_i$. The limiting function $\vec{u}^*= ( u_1^*, \ldots, u_m^* )$, the limit of $\vec{u}^\eps$ as $\eps \to 0$, partitions $\Gamma$ so we can define $\Gamma_i = \{ u^*_i > 0 \}$ and  $u_{j}^*=0$ in $\Gamma_{j}$, $j\ne i$. We note also that setting $v^{*}_{i}:=u^{*}_{i}- \sum_{j \neq i}u^{*}_{j}$ we have $\Gamma_i = \{ v^*_i > 0 \}$.
A possible  disadvantage of this method is that it is not clear how to relate $\vec{u}^\eps$ to a partition $\{ \Gamma_i \}$ when $\eps$ is fixed. Possibilities for defining $\Gamma_{i}^{\eps}$ include $\Gamma_{i}^{\eps}=\{ u_i ^{\eps}> c(\eps) \}$ or $\Gamma_{i}^{\eps}=\{v_{i}^{\eps}>0\}$ where $v^{\eps}_{i}:=u^{\eps}_{i}- \sum_{j \neq i}u^{\eps}_{j}$.

\subsection{Outline}

In the remainder of this paper, we will give a suitable discretisation of this approach using the surface finite element method. We will propose an algorithm to solve the discretised optimisation problem and give practical details of how we implement this method. Our experience is  that the eigenfunction segregation method 
performs very well. 
Our results section consists of three parts. First, we will test our algorithm in the case of three partitions on the sphere for which we know the absolute minimiser. We will then compute partitions of the sphere  for  larger values of $m$ and make some observations about the structure. Finally we consider other surfaces to see the different effects of curvature and different genus surfaces. The computations lead to some natural conjectures. 

%%% Local Variables:
%%% mode: latex
%%% TeX-master: "../partition"
%%% End:

\section{Computational method}

\subsection{Discretisation}

We start the discretisation by taking a polyhedral approximation $\Gamma_h$ of $\Gamma$. We assume that $\Gamma_h$ consists of a shape regular triangulation $\T_h$ where $h$ is the maximal diameter of a simplex (triangle for $n=2$) in $\T_h$. We will denote by $\N_h$ the vertices of $\Gamma_h$ and call $\Gamma_h$ a triangulated surface. We suppose that $\Gamma_{h}$ interpolates $\Gamma$ in the sense that the vertices of triangles of $\Gamma_{h}$ lie on $\Gamma$.

Over this triangulation, we define  two continuous finite element spaces, a space of scalar valued functions $S_h$ and a space of vector valued functions $\vec{S}_h$. These are given by
\begin{align*}
  S_h & = \{ \chi_h \in C( \Gamma_h ) : \chi_h |_T \mbox{ is affine linear, for all } T \in \T_h \} \\
  \vec{S}_h & = \{ \vec{\eta}^h = ( \eta^h_1, \ldots, \eta^h_m ) \in C( \Gamma_h; \R^m ) : \eta^h_i \in S_h \mbox{ for } i = 1,\ldots,m \}.
\end{align*}

We can directly formulate the discrete version of Problem~\ref{pb:phase-cts}.

\begin{problem}
  Given a positive integer $m$, a triangulated surface $\Gamma_h$ and $\eps > 0$, find $\vec{u}^{\eps,h} = ( u_1^{\eps,h}, \ldots, u_m^{\eps,h} ) \in \vec{S}_h$ to minimise
  \begin{equation}
    \Epf^{\eps,h}( \vec{u}^{\eps,h} ) = \frac{1}{2} \sum_{i=1}^m \int_{\Gamma_h} \abs{ \nabla_{\Gamma_h} u^{\eps,h}_i }^2 \dd \sigma_h + \int_{\Gamma_h} F_\eps( \vec{u}^{\eps,h} ) \dd \sigma_h.
  \end{equation}
\end{problem}

Our optimisation strategy will be to directly solve a discretisation of the gradient flow equations. Discretising in space first, we seek a time dependent finite element function $\vec{u}^{\eps,h} \in C^1( \R_+ ;\vec{S}_h )$ and $\lambda^{\eps,h} \colon \R_+ \to \R^m$ satisfying $||u_{i}^{\eps,h}||^{2}_{\Gamma_{h}}=1~~i=1,2,....m,$
\begin{equation}
  \label{eq:uh-gf}
  \begin{aligned}
    & \int_{\Gamma_h} \partial_t u^{\eps,h}_i \chi_h + \nabla_{\Gamma_h} u^{\eps,h}_i \cdot \nabla_{\Gamma_h} \chi_h \dd \sigma_h \\
    & \qquad = \int_{\Gamma_h} \lambda_i^{\eps,h} u_i^{\eps,h} \chi_h - \frac{2}{\eps^2} \left( \sum_{j \neq i} ( u_j^{\eps,h} )^2 \right) u_i^{\eps,h} \chi_h \dd \sigma_h
    && \mbox{ for all } \chi_h \in S_h \\
    & \vec{u^{\eps,h}}( \cdot, 0 ) = \vec{u}^{h,0}.
  \end{aligned}
\end{equation}
Here, $\vec{u}^{h,0} = ( u^{h,0}_1, \ldots, u^{h,0}_m )$ is initial data in $\vec{S}_h$ such that $\sum_{j \neq i} ( u^{h,0}_i )^2 ( u^{h,0}_j )^2 = 0$ for $i=1,\ldots,m$.

We discretise in time using a operator splitting strategy similar to a scheme proposed by \cite{DuLin09}. At each time step, we first solve one step of the heat equation, then solve an ordinary differential equation for the nonlinear terms, and use a projection to deal with the Lagrange multiplier.

\subsection{Computational method}

The operator splitting method is as follows.

\begin{algorithm}
  Given $\eps > 0$, a positive integer $m$, a time step $\tau > 0$ and an initial condition $\vec{u}^{h,0} = ( (u^{h,0}_1), \ldots, (u^{h,0}_m) ) \in \vec{S}_h$ with $\sum_{j \neq i} u^{h,0}_i(z)^2 u^{h,0}_j(z)^2 = 0$ for all $z \in \N_h$ and $i=1,\ldots,m$, for $k=0,1,2,\ldots$,
  \begin{enumerate}
  \item Solve one time step of the heat equation for $i=1,\ldots,m$ using implicit Euler. We wish to find $\vec{\tilde{u}}^{\eps,h} = ( \tilde{u}^{\eps,h}_1, \ldots, \tilde{u}^{\eps,h}_m ) \in \vec{S}_h$
    \begin{equation*}
      \int_{\Gamma_h} \tfrac{1}{\tau} \big( \tilde{u}^{\eps,h}_i - ( u^{\eps,h}_i )_k \big) \chi_h
      + \nabla_{\Gamma_h} \tilde{u}^{\eps,h}_i \cdot \nabla_{\Gamma_h} \chi_h \dd \sigma_h = 0 \quad \mbox{ for all } \chi_h \in S_h, i = 1,\ldots,m.
    \end{equation*}
  \item Solve the nonlinear terms exactly as ordinary differential equation at each node. For all nodes $z \in \N_h$ and $i = 1,\ldots,m$, find $\hat{u}^{\eps,h}_i(z) \colon [ t^k, t^{k+1} ] \to \R$ such that
    \begin{equation*}
      \frac{\mathrm{d}}{\mathrm{d}t} \left( \hat{u}^{\eps,h}_i( z )( t ) \right)
      = -\left( \frac{2}{\eps^2} \sum_{j\neq i} (\tilde{u}^{\eps,h}_j(z))^2 \right) \hat{u}^{\eps,h}_i( z )( t ),
      \quad
      \hat{u}^{\eps,h}_i( z )( t^k ) = \tilde{u}^{\eps,h}_i( z ).
    \end{equation*}
  \item Find the new solution $(\vec{u}^{\eps,h})_{k+1}$ by normalising the final time solution $(\hat{u}^{\eps,h}_1(\cdot)(t^{k+1}), \ldots, \hat{u}^{\eps,h}_m(\cdot)(t^{k+1}))$:
    \begin{equation*}
      ( u^{\eps,h}_i( z ) )_{k+1} = \frac{ \hat{u}^{\eps,h}_i( z )( t^{k+1} ) }{ \norm{ \hat{u}^{\eps,h}_i( \cdot )( t^{k+1} ) }_{L^2(\Gamma_h)} } \quad \mbox{ for all } z \in \N_h, i = 1, \ldots, m.
    \end{equation*}
  \end{enumerate}
\end{algorithm}

Similarly to \citet{BaoDu04}, one can show  an energy decreasing property for this scheme.
The method is the same as the scheme of \citet{DuLin09} except we exchange a Gauss-Seidel iteration in step 2 for a Jacobi iteration. The ordinary differential equation from step 2 can be solved exactly to give:
\begin{equation*}
  \hat{u}_i^{\eps,h}( z )( t^{k+1} ) = \tilde{u}^{\eps,h}_i( z ) \exp\left( - \frac{\tau}{\eps^2} \sum_{j\neq i} (\tilde{u}^{\eps,h}_j(z)(t))^2 \right).
\end{equation*}
Using this solution, we write a more practical version of step 2 as
\begin{enumerate}
\setcounter{enumi}{1}
\item For each node $z \in \N_h$,
  \begin{enumerate}
  \item For $i = 1,\ldots,m$, compute $\tilde{u}^{\eps,h}_i(z)^2$;
  \item Find $S = \sum_{i=1}^m \tilde{u}^{\eps,h}_i(z)^2$;
  \item For $i = 1,\ldots,m$, compute $\hat{u}_i^{\eps,h}( z )( t^{k+1} )$ by
    \begin{equation*}
      \hat{u}_i^{\eps,h}( z )( t^{k+1} ) = \tilde{u}_i^{\eps,h}( z ) \exp\left( - \frac{2\tau}{\eps^2} ( S - \tilde{u}_i^{\eps,h}( z )^2 ) \right)  
    \end{equation*}
  \end{enumerate}
\end{enumerate}

We stop the computation when the change in energy is less than $10^{-6}$. In order to  reduce the computational cost this is only calculated every $M_{\tau}$ iterations where $0.1=M_{\tau}\tau$.  

Since, in general, we do not know the configuration of the optimal domains, we initialise the computations with a random initial condition. We loop over the grid nodes $z \in \N_h$ and uniformly at random choose one value $i \in \{ 1, \ldots, m \}$ and set $(u^h_0)( z )_{i} = 1$ and $(u^h_0)( z )_{j} = 0$ for $j \neq i$ then normalise each component, $(u^h_0)_i$, for $i = 1,\ldots,m$, in $L^2(\Gamma)$. As a result the first linear solve for the heat equation step will take more iterations, however the difference is not significant in this case.

\begin{xrem}
  In practice, we find this operator splitting method to be stable and efficient. If we discretised \eqref{eq:uh-gf} in time directly using the Lagrange multiplier, we would have the choice to take  the Lagrange multiplier  implicitly or explicitly. An implicit discretisation  would leave a fully coupled system of equations to solve, which would not be so easily implemented using parallel  high performance computing techniques. An explicit discretisation  would imply a time step restriction based on the size of the maximum $H^1$-semi norm of each component. We wish to start with a random initial condition in order to avoid local minima, however this has a very large $H^1$-semi norm which would give an unfeasible time step restriction. All three methods are considered for the flat problem in the time discrete-space continuous case by \citet{DuLin09}.
\end{xrem}

\subsection{Parallel computations}

The algorithm has been formulated so that we can use high performance computing to implement the optimisation. The key idea is to store the solution over $m$ parallel processors and perform most of the computations on a single processor. Communication between processors is kept to a minimum.

We distribute the solution $\vec{u}^{\eps,h}$ over $m$ processors so that processor $i$ stores $u_i^{\eps,h}$. At each time step, each processor performs one linear solve (step 1), one loop over all nodes communicating with all other nodes to perform the sum in step 2(b) (step 2), then one more loop over all nodes to normalise the solution (step 3). In particular, computing sum in step 2(b) over all $j$ is more efficient then computing the sum over all other $j \neq i$.

A similar approach was also taken to parallelisation by \citet{BouBucOud10} who computed up to $512$ partitions. Our approach performs very well for $m \le 32$. At the moment we restricted to this number of partitions  because we wish to have a meaningful number of elements in each partition. It is possible that one may gain efficiency by using  an adaptive mesh refinement on the unstructured grids enabling  sufficiently accurate computations with a larger number of partitions. This is left for future work.

All test cases were implemented using the Distributed and Unified Numerics Environment (DUNE) \citep{basbladed08,basbladed08a}. Matrices are assembled using the DUNE-FEM \citep{dedklonol10} and solved using a conjugate gradient method preconditioned with algebraic multigrid Jacobi preconditioner from DUNE-ISTL \citep{BasBla07}. Parallelisation is performed using \textsf{MPI}. All visualisation is performed in ParaView \citep{paraview}. The code we have written for the simulations in this paper is available at
\begin{quote}
  \centering
  \url{http://users.dune-project.org/projects/dune-partition}
\end{quote}

%%% Local Variables:
%%% mode: latex
%%% TeX-master: "../partition"
%%% End:

\section{Results}

\subsection{Convergence tests for three partitions of the sphere}\label{Ypartition}

Bishop's conjecture (Conjecture~\ref{conj:bishop}) suggests that the Y-partition is optimal in the case $m=3$ on the sphere. This corresponds (up to rotations of the sphere) to $\Gamma_1 = \{ 0 < \varphi < 2 \pi / 3 \}$, $\Gamma_2 \{ -2 \pi / 3 < \varphi < 0 \}$ and $\Gamma_3 = \{ \abs{ \varphi } > 2 \pi / 3 \}$. We can compute that the first eigenfunctions are:
\begin{align*}
  u_1( \theta, \varphi ) & = \sin( \tfrac{3\phi}{2} ) ( \sin \theta )^{\frac{3}{2}} && \mbox{ on } \Gamma_1 \\
  u_2( \theta, \varphi ) & = -\sin( \tfrac{3\phi}{2} ) ( \sin \theta )^{\frac{3}{2}} && \mbox{ on } \Gamma_2 \\
  u_3( \theta, \varphi ) & = \sin( \tfrac{3 \abs{\phi} }{2} - \pi ) ( \sin \theta )^{\frac{3}{2}} && \mbox{ on } \Gamma_3.
\end{align*}
Each of these eigenfunctions has eigenvalue $15/4$. We will test our scheme by checking the rate of convergence to the Y-partition.

\begin{figure}[tbh]
  \centering
  \includegraphics{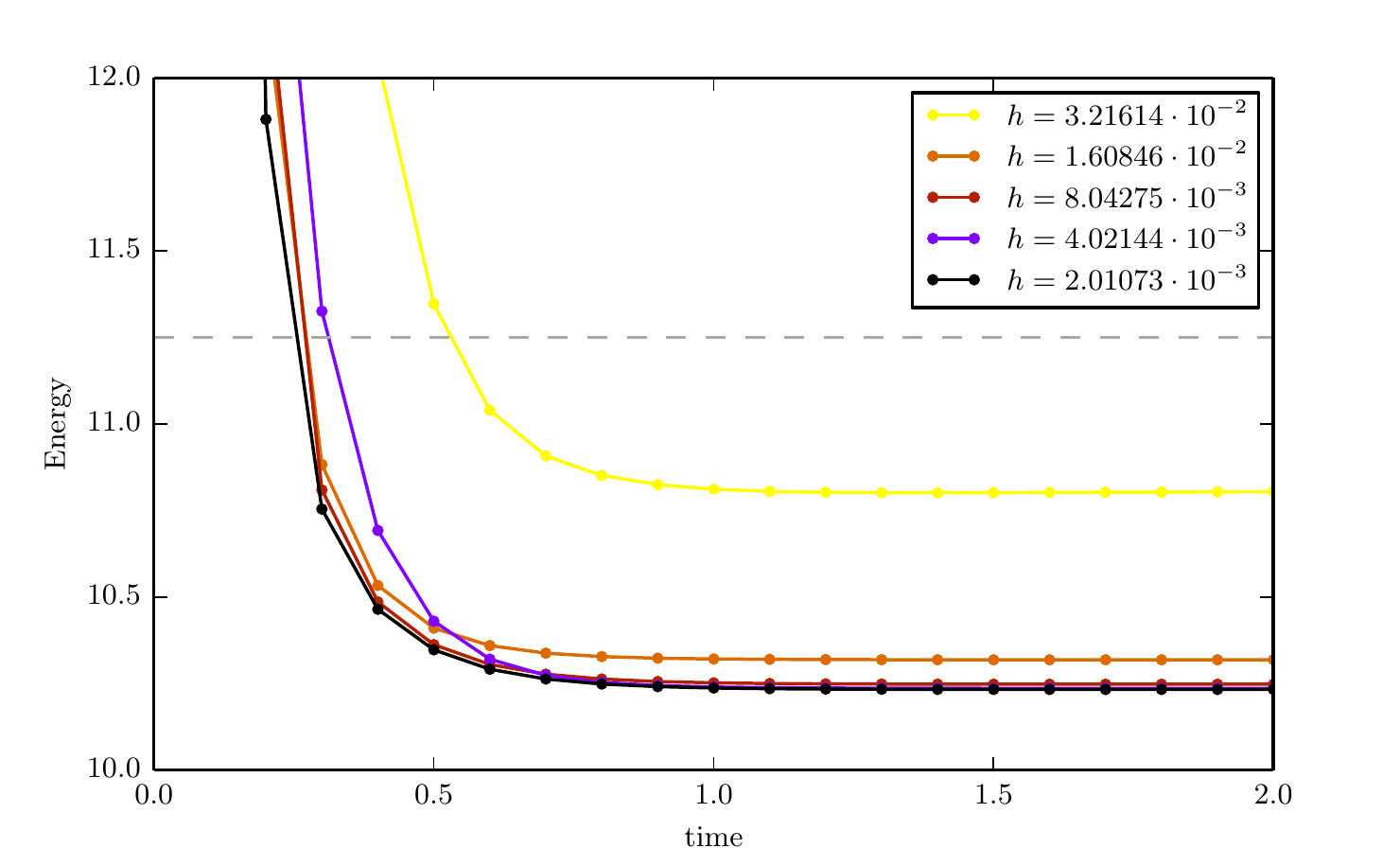}
  \caption{Convergence with respect to discretisation parameters for $\eps = 5 \cdot 10^{-3}$ to the Y-partition on the sphere. The dashed grey line is the Y-partition energy for $\eps = 0$.}
  \label{fig:h-conv}
\end{figure}

We first test convergence with respect to the discretisation parameters. We perform our algorithm at $\eps = 5 \cdot 10^{-3}$ and $\tau = 10^{-4}$ over five levels of mesh refinement, reducing from $h = 3.21614 \cdot 10^{-2}$ to $h = 2.01073 \cdot 10^{-3}$. We compute until $t = 2$. We have plotted the energy along the time evolution in Figure~\ref{fig:h-conv} and see good convergence. We have also included a dashed line at the Y-partition energy $45/4$ for $\eps = 0$. We see that for a given $\eps$ the error in energy can be large.

To test the convergence of the regularisation we compute the minimizer for a sequence for values for $\eps$. We start on a coarse mesh ($h=3.21614 \cdot 10^{-2}$) with $\tau = 8 \cdot 10^{-4}$, once we have reached a minimizer, we refine the mesh by bisecting elements once (two bisections reduces $h$ roughly by half) and reduce $\tau$ by a factor $1/\sqrt{2}$. Instead of computing a new random initial condition after each refinement, we use the previous minimiser as the new initial condition.

We define $S_\eps$ to be part of the energy associated with regularisation:
\begin{equation}
  \label{eq:Seps-u}
  S_\eps( \vec{u}^{\eps,h} ) := \int_{\Gamma_h} F_\eps( \vec{u}^{\eps,h} ) \dd \sigma_h = \frac{1}{\eps^2} \int_{\Gamma_h} \sum_{i=1}^m \sum_{j \neq i} (u^{\eps,h}_i)^2 (u_j^{\eps,h})^2 \dd \sigma_h.
\end{equation}
These values illustrate the convergence of the relaxation to the exact problem. We expect $S_\eps \to 0$ as we know that we recover a minimiser of the partition problem as $\eps \to 0$.

\begin{figure}%[tbh]
%  \showthe\textwidth > 475.16101pt.
  \centering
  \includegraphics{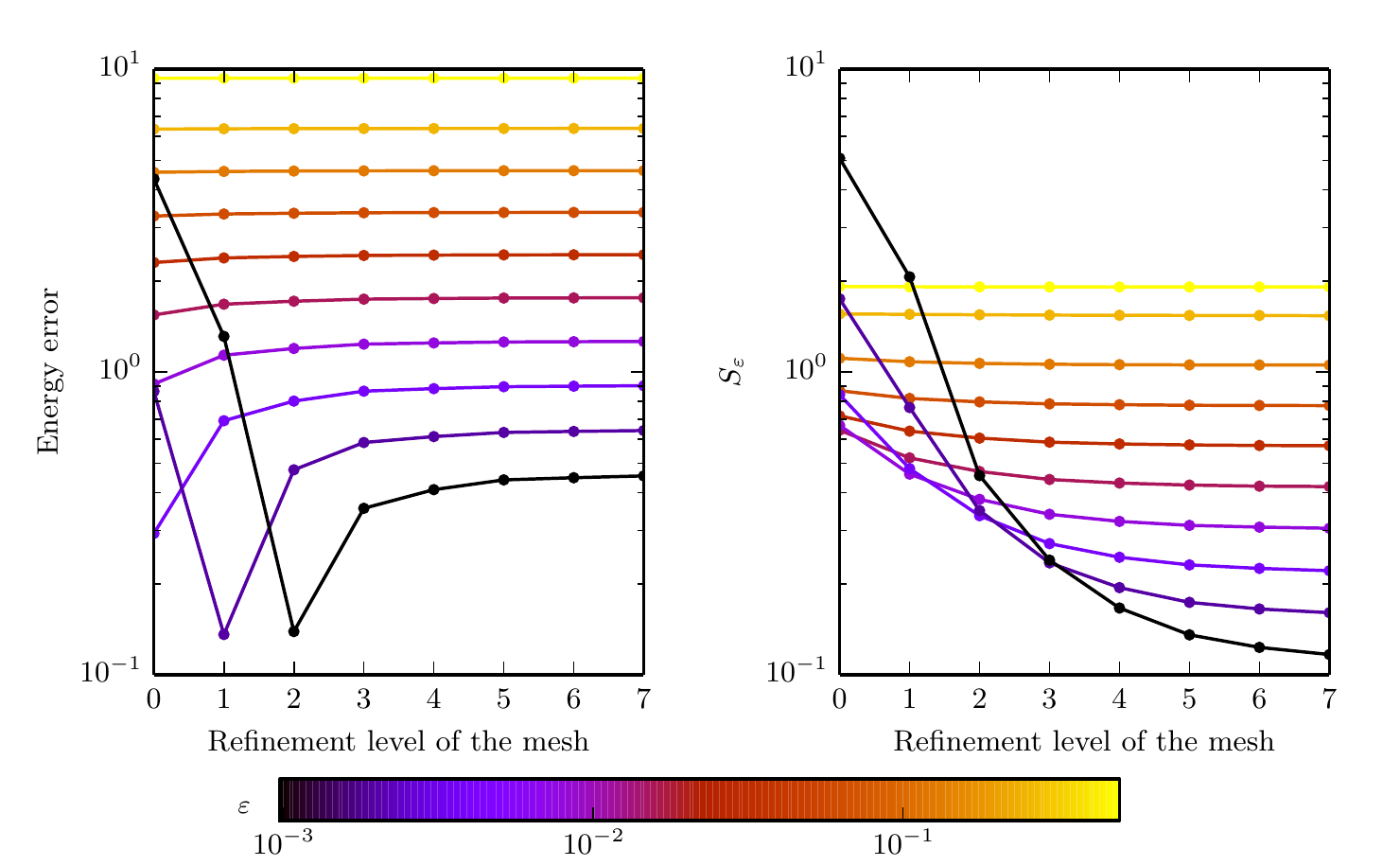}
  \caption{Convergence with respect to $\eps$ to the Y-partition on the sphere. The energy error is difference to $45/4$, the Y-partition energy for $\eps = 0$, and $S_\eps$ is given by \eqref{eq:Seps-u}.}
  \label{fig:u-eps}
\end{figure}

We have computed the full and regularisation energy at each minimiser. The results are shown in Table~\ref{tab:u-eps} and Figure~\ref{fig:u-eps}. The tables also show the experimental order of convergence (eoc) which is computed via the formula
\begin{equation*}
  (\mathrm{eoc})_i = \frac{\log( \mathrm{error}_i / \mathrm{error}_{i-1} )}{\log( 1/2 )}.
\end{equation*}
where $\mathrm{error}_i$ is the error in energy against the Y-partition at refinement level $i$.

\begin{table}%[tbh]
  \centering
  \begin{tabular}{cccccc}
    \hline
$\eps$ & Energy & Energy error & (eoc) & $S_\eps$ & (eoc) \\
\hline$5.00000 \cdot 10^{-1}$ & $1.9100$ & $9.3400$ & --- & $1.9098$ & --- \\
$2.50000 \cdot 10^{-1}$ & $4.8759$ & $6.3741$ & $0.5512$ & $1.5350$ & $0.3151$ \\
$1.25000 \cdot 10^{-1}$ & $6.6257$ & $4.6243$ & $0.4630$ & $1.0548$ & $0.5413$ \\
$6.25000 \cdot 10^{-2}$ & $7.8829$ & $3.3671$ & $0.4577$ & $0.7751$ & $0.4444$ \\
$3.12500 \cdot 10^{-2}$ & $8.8095$ & $2.4405$ & $0.4643$ & $0.5714$ & $0.4400$ \\
$1.56250 \cdot 10^{-2}$ & $9.4907$ & $1.7593$ & $0.4721$ & $0.4188$ & $0.4482$ \\
$7.81250 \cdot 10^{-3}$ & $9.9880$ & $1.2620$ & $0.4793$ & $0.3050$ & $0.4576$ \\
$3.90625 \cdot 10^{-3}$ & $10.3487$ & $0.9013$ & $0.4856$ & $0.2209$ & $0.4652$ \\
$1.95312 \cdot 10^{-3}$ & $10.6088$ & $0.6412$ & $0.4912$ & $0.1605$ & $0.4605$ \\
$9.76562 \cdot 10^{-4}$ & $10.7958$ & $0.4542$ & $0.4974$ & $0.1168$ & $0.4591$ \\
    \hline
  \end{tabular}
  \caption{Results of convergence test in $\eps$ for numerical tests for three partition case. Energy is $\Epf^\eps$ at the best computed partition, energy error is the difference to $45/4$, the Y-partition energy for $\eps = 0$, and $S_\eps$ is given by \eqref{eq:Seps-u}.}
  \label{tab:u-eps}
\end{table}

The eigenfunction segregation approach performs very well with respect to convergence in $\eps$. We observe order $\eps^{\frac{1}{2}}$ convergence both for the full energy and also for $S_\eps$. The errors are still quite large for reasonable sized values of $\eps$ so we must take very small values of $\eps$ to trust any predictions of energy values using this method.

%%%%%%%%%%%%%%%%%%%%%%%%%%%%%%%%
\subsection{Computed partitions of the sphere  for $m \ge 3$}
\label{sec:many-partition}
%%%%%%%%%%%%%%%%%%%%%%%%%%%%%

We proceed with the following refinement rules. We initialise the problem with a random initial condition for $\eps_0 = \frac12$, $\tau_0 = 8 \cdot 10^{-4}$ on a mesh $\Gamma_{h,0}$, then for $l = 0,1,2,\ldots$, we find a minimiser of the $\eps_l$-problem on $\Gamma_{h,l}$, then refine the mesh globally by bisecting all elements, and find $\eps_{l+1}$ and $\tau_{l+1}$ as
\begin{equation*}
  \eps_{l+1} = \sqrt{2} \eps_l \qquad \tau_{l+1} = \sqrt{2} \tau_l.
\end{equation*}
We use the optimal function for level $l-1$ as the initial condition on level $l$. The final parameters are given in Table~\ref{tab:dofs-sphere}.

\begin{table}
  \centering
  \begin{tabular}{cccc}
    \hline
    $m$ & $l$ & Degrees of freedom & $\eps$ \\
    \hline
    3 & 9 & $579\,830$ & $6.25 \cdot 10^{-4}$ \\
    4 & 9 & $786\,440$ & $6.25 \cdot 10^{-4}$ \\
    5 & 9 & $983\,050$ & $6.25 \cdot 10^{-4}$ \\
    6 & 9 & $1\,179\,660$ & $6.25 \cdot 10^{-4}$ \\
    7 & 9 & $1\,376\,270$ & $6.25 \cdot 10^{-4}$ \\
    8 & 7 & $196\,624$ & $1.25 \cdot 10^{-3}$ \\
    16 & 7 & $393\,248$ & $1.25 \cdot 10^{-3}$ \\
    32 & 7 & $786\,496$ & $1.25 \cdot 10^{-3}$ \\
    \hline
  \end{tabular}
  
  \caption{Final parameters for computations on the sphere.}
  \label{tab:dofs-sphere}
\end{table}

\begin{figure}[p]
  \centering

  \begin{subfigure}[b]{0.3\textwidth}
    \includegraphics[width=\textwidth]{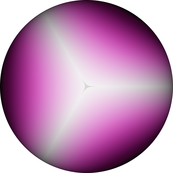}
    \caption{$m=3$, 3 lens (pink)}
  \end{subfigure}
  \hfill
  \begin{subfigure}[b]{0.3\textwidth}
    \includegraphics[width=\textwidth]{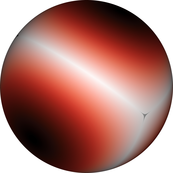}
    \caption{$m=4$, 4 triangles (red)}
  \end{subfigure}
  \hfill
  \begin{subfigure}[b]{0.3\textwidth}
    \includegraphics[width=\textwidth]{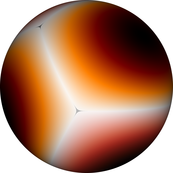}
    \caption{$m=5$, 2 triangles (red) and 3 quadrilaterals (orange)}
  \end{subfigure}

  \vfill

  \begin{subfigure}[b]{0.3\textwidth}
    \includegraphics[width=\textwidth]{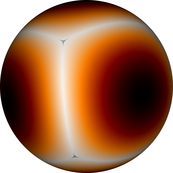}
    \caption{$m=6$, 6 quadrilaterals (orange)}
  \end{subfigure}
  \hfill
  \begin{subfigure}[b]{0.3\textwidth}
    \includegraphics[width=\textwidth]{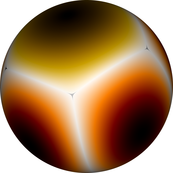}
    \caption{$m=7$, 5 quadrilaterals (orange) and 2 pentagons (yellow)}
  \end{subfigure}
  \hfill
  \begin{subfigure}[b]{0.3\textwidth}
    \includegraphics[width=\textwidth]{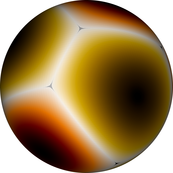}
    \caption{$m=8$, 4 quadrilaterals (orange) and 4 pentagons (yellow)}
  \end{subfigure}

  \vfill

  \hspace*{\fill}
  \begin{subfigure}[b]{0.3\textwidth}
    \includegraphics[width=\textwidth]{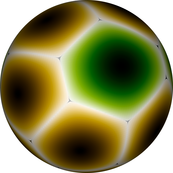}
    \caption{$m=16$, 8 A and 4 B pentagons (yellow) and 4 hexagons (green)}
  \end{subfigure}
  \hfill
  \begin{subfigure}[b]{0.3\textwidth}
    \includegraphics[width=\textwidth]{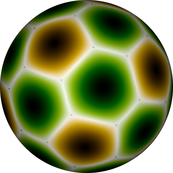}
    \caption{$m=32$, 12 pentagons (yellow) and 20 hexagons (green)}
  \end{subfigure}
  \hspace*{\fill}

  \caption{Plots of the minimising configurations $\{ \Gamma_i^{\eps,h} \}_{i=1}^m$ with void regions in grey. Colours only in the online version. Each partition is coloured according to the polygon type and shaded by the eigenfunction from white for $u_i = 0$ to black for $u_i$ at the maximum.}
  \label{fig:many-m}
\end{figure}

\begin{figure}
  \centering
  \includegraphics{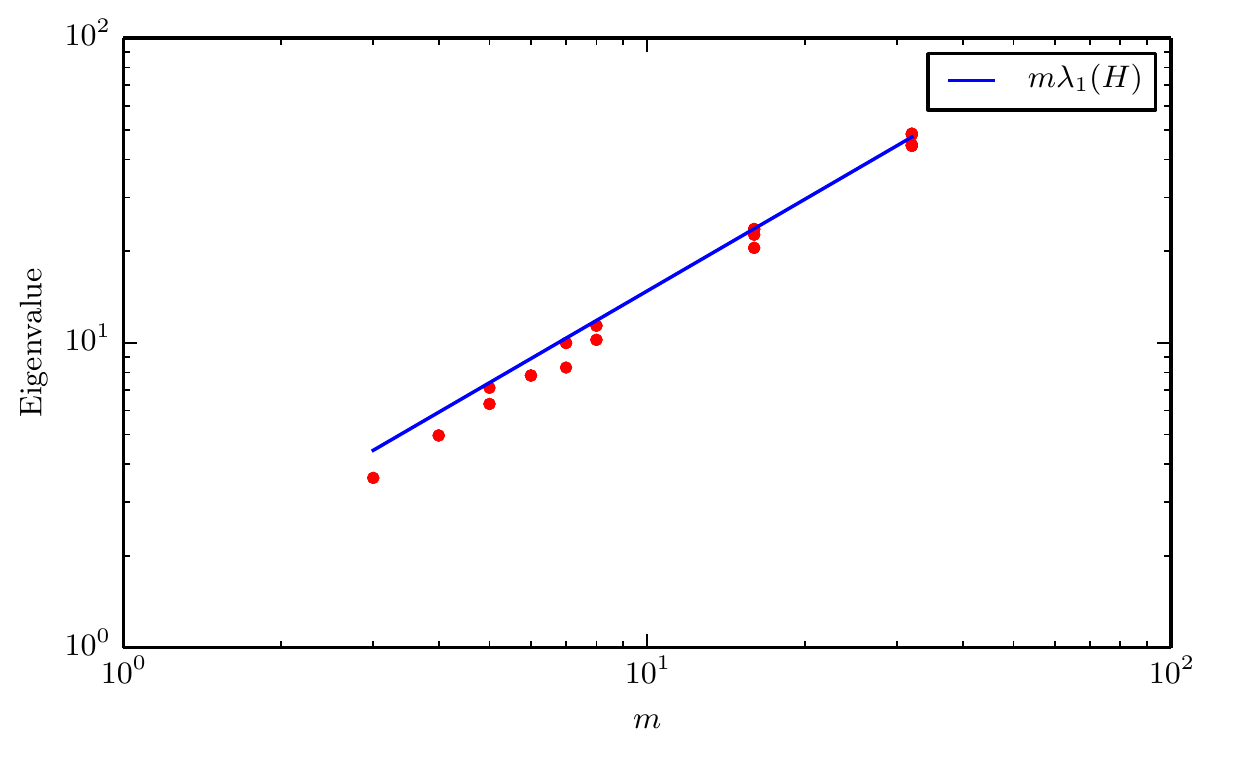}

  \caption{Plot of the eigenvalues at different values of $m$. The blue line is $m \lambda_1( H )$ where $H$ is the planar hexagon with area $4 \pi$ (equal to the surface area of the sphere).}
  \label{fig:many-m-plot}
\end{figure}

Plots of the solutions for several values of $m$ are given in Figure~\ref{fig:many-m}. Observe that the colour coding of these figures indicates the partitions using the computed values of the eigenfunctions. Eigenvalue estimates are computing by taking the mean $H^1$-semi norm of the components.  The computed eigenvalues are plotted in Figure~\ref{fig:many-m-plot}. Theorem~3 of the work by \citet{CafLin07} proves that the energy scales like $\lambda_m(\Gamma)$ up to a constant factor. Using Weyl's asymptotics, we see that in two space dimensions this means that the average eigenvalue is bounded above and below by $m$ times a constant. This is indicated by the blue line which is $m$ times the first eigenvalue corresponding to a hexagon $H$ of area $4 \pi$ (the surface area of the sphere) -- this is the conjectured average eigenvalue for large $m$ in the plane \citep{CafLin07}. Our results indicate a similar scaling property for the sphere.

\begin{figure}
  \centering
  \includegraphics[width=0.8\textwidth]{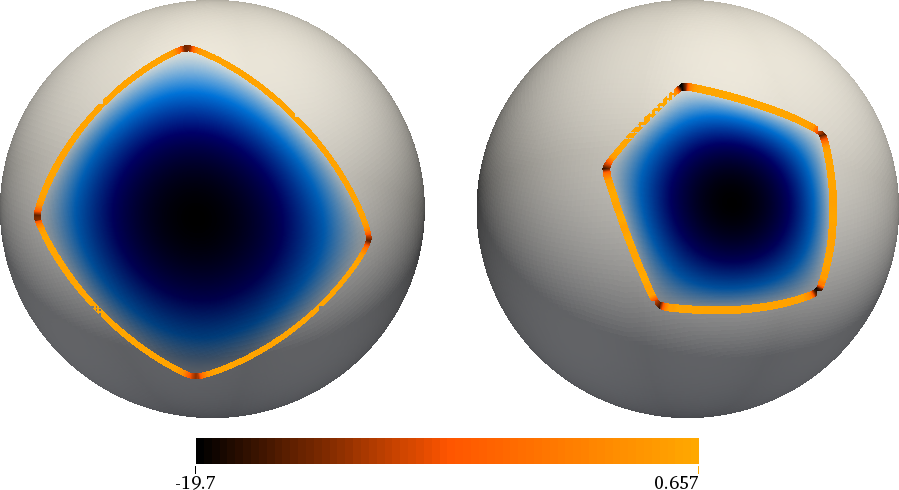}

  \caption{Plots of one partition and $\kappa_g$ for $m=8$ (left) and $m=16$ (right). The value of $u_i^{\eps,h}$ is shown on a black to white scale and $\kappa_g$ is plotted on the curve $\{ v_i^{\eps,h} = 0 \}$ on a black to orange scale.}
  \label{fig:kappag}
\end{figure}

Rather than just  using  the computed eigenfunction values, as mentioned earlier,  we may  define an approximate partition by
\begin{equation}
  \label{eq:vi}
  \Gamma_i^{\eps,h} := \left\{ x \in \Gamma : v_i^{\eps,h}(x) := u_i^{\eps,h}(x) - \sum_{j \neq i} u_j^{\eps,h}(x) >0\right\} \quad \mbox{ for } i = 1,\ldots,m.
\end{equation}
We motivate the use of this  definition by noting that each $u_i^{\eps,h}$ is positive and the supports of $\{ u_i^{\eps,h} \}$ overlap, hence this function is zero only surrounding one partition where $u_i^{\eps,h} = u_j^{\eps,h}$ for some $j \neq i$. Note that these sets will not cover $\Gamma$ and there will be a small void between regions. Furthermore we may use $ v_i^{\eps,h}$ in the following interesting way.
Suppose that $\gamma$ is a curve on $\Gamma$ defined by as the zero level set of a function $\phi$, $\gamma = \{ \phi =  0\}$, then the geodesic curvature of $\gamma$, which we denote by $\kappa_g$ is given by
\begin{equation}
  \label{eq:kappag}
  \kappa_g = \nabla_\Gamma \cdot \left( \frac{\nabla_\Gamma \phi}{\abs{ \nabla_\Gamma \phi }} \right).
\end{equation}
We can use ParaView's gradient reconstruction function to compute an approximation of $\kappa_g$ over the interface at the boundary of each partition $\Gamma_i$ using $\phi = v^{\eps,h}_i$. An example of this is shown in Figure~\ref{fig:kappag}. We see that this value is small away from junctions.

\begin{figure}
  \centering
  \includegraphics[width=0.45\textwidth]{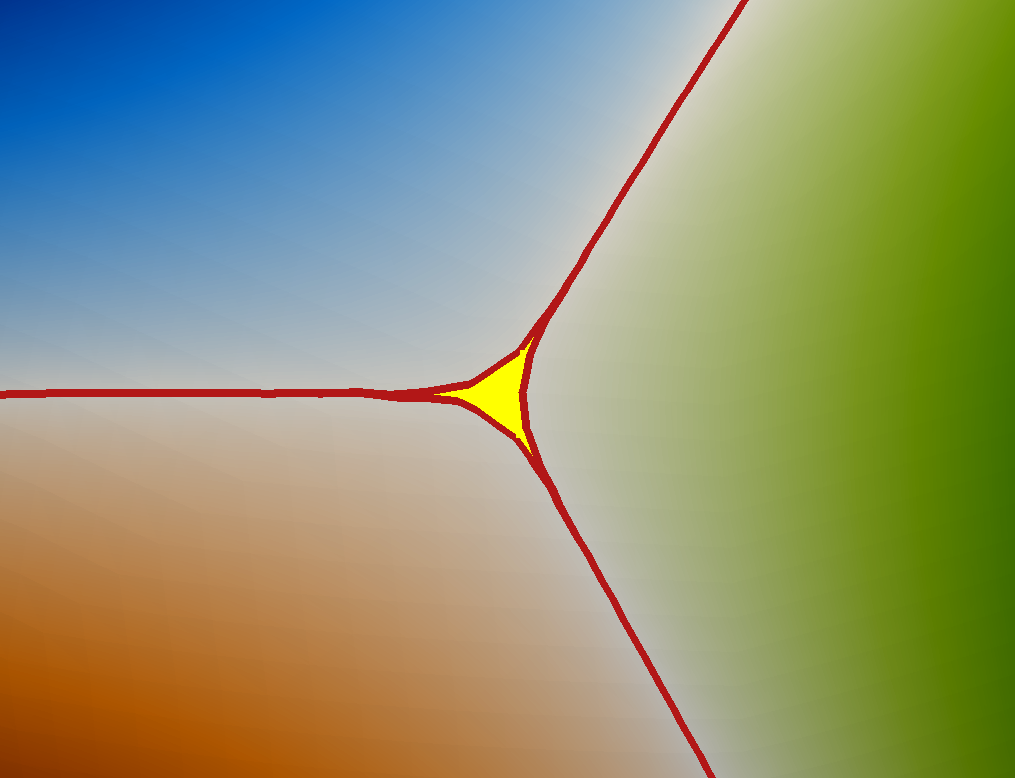}
  \quad
  \includegraphics[width=0.45\textwidth]{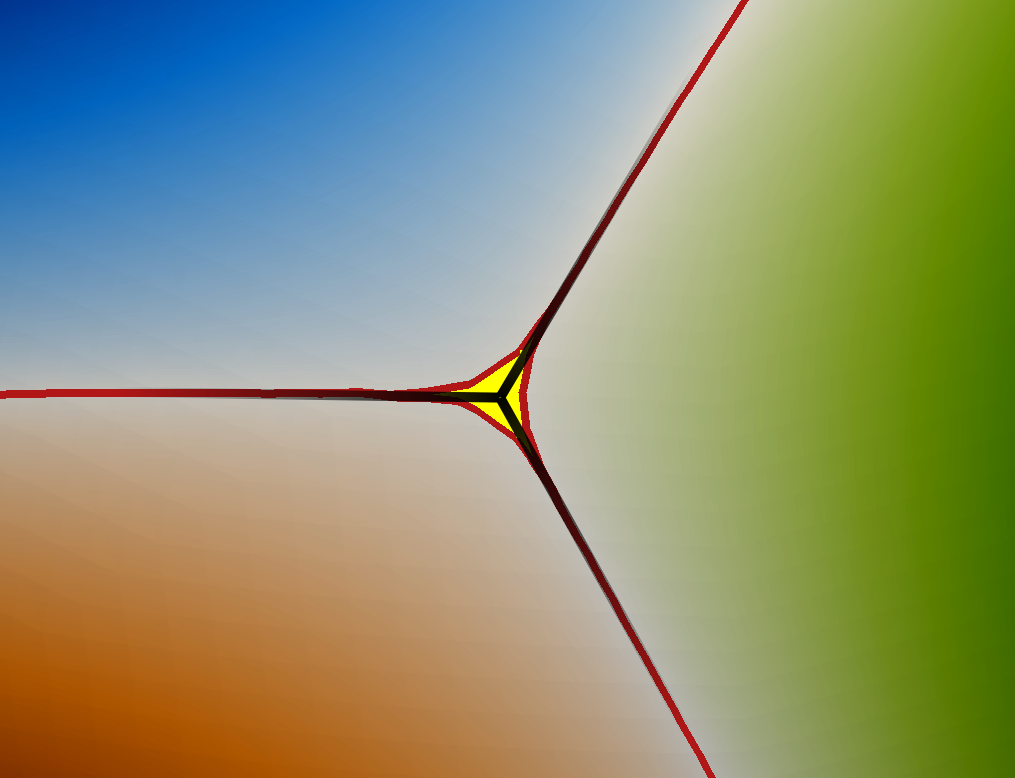}

  \caption{A zoom of a triple junction on the sphere. Three partitions $\{ v_i^{\eps,h} > 0 \}$ are coloured on blue, green and orange according to the eigenfunction $u_i^{\eps,h}$ with red boundaries at $\{ v_i^{\eps,h} = 0 \}$. The void region is shown in yellow. Additionally in the right plot we have added black lines which would correspond to an equal angle triple junction.}
  \label{fig:trip-point}
\end{figure}

We observe that at junctions three partitions coincide with equal angles. See, for example, Figure~\ref{fig:trip-point}. This is consistent with the results of \citet{HelOstTer10} who prove, for the case $p=\infty$, that all partitions have an equal angle property. From our results it is difficult to quantify this result since at any triple point there is a void region because of our regularisation. Also in Figure~\ref{fig:trip-point}, we have superimposed an equal angle triple junction which shows good agreement to results we have. We can consider a reduced problem of finding the first eigenvalue over partitions of the unit disk. We find with three equal partitions (similar to the Y-partition) the total energy is approximately $60.6 (= 3 \cdot 20.2)$ and for four partitions, one in each quadrant, the total energy is approximately $105.6 (=4 \cdot 26.4)$. Taking three partitions leads to a significant reduction in energy.

Table~\ref{tab:many-m} shows one representative of each polygon similarity class and more details of the best estimate of the energy and also the similarity classes of polygons. The energy calculation shows the values of each eigenvalue (mean and standard deviation for each similarity class of polygons) and also $S_\eps$ for each of the final configurations.
We note that for $m = 3,4,6$, our optimal configuration are equi-spectral and for the case $m=4$ we recover a spherically projected tetrahedron as conjectured by \citet{HelOstTer10}.
Thus we conjecture that there partitions are optimal for the case $p=\infty$ also.

\begin{table}[p]
  \centering
  \newcolumntype{C}[1]{>{\centering\let\newline\\\arraybackslash\hspace{0pt}}m{#1}}
  \begin{tabular}{|C{0.025\textwidth}|C{0.525\textwidth}|C{0.45\textwidth}|}
    \hline
    $m$ & Shape & Energy information \\
    \hline
    3 & %
    \begin{tabular}{c}
      3 lens \\ \includegraphics[width=0.75in]{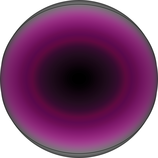}
    \end{tabular} &
    Lens eigenvalue: $3.605\: (2.59 \cdot 10^{-4})$

    $S_\eps$: $0.072$

    Total energy: $10.887$ \\
    \hline
    4 & %
    \begin{tabular}{c}
      4 triangles \\ \includegraphics[width=0.75in]{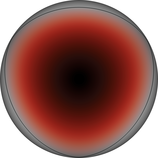}
    \end{tabular} &
    Triangle eigenvalue: $4.966\: (2.46 \cdot 10^{-4})$

    $S_\eps$: $0.121$

    Total energy: $19.987$  \\
    \hline
    5 &
    \begin{tabular}{c}
      2 triangles \\ \includegraphics[width=0.75in]{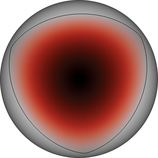}
    \end{tabular}
    and
    \begin{tabular}{c}
      3 quadrilaterals \\ \includegraphics[width=0.75in]{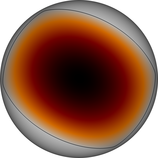}
    \end{tabular} &
    Triangle eigenvalue: $7.118\: (3.35 \cdot 10^{-4})$

    Quadrilateral eigenvalue: $6.302$

    $S_\eps$: $0.187$

    Total energy: $33.330$ \\
    \hline
    6 &
    \begin{tabular}{c}
      6 quadrilaterals \\ \includegraphics[width=0.75in]{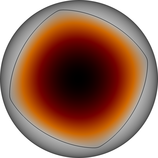}
    \end{tabular} &
    Quadrilateral eigenvalue: $7.812\: (7.22 \cdot 10^{-4})$

    $S_\eps$: $0.248$

    Total energy: $47.122$ \\
    \hline
    7 &
    \begin{tabular}{c}
      5 quadrilaterals \\ \includegraphics[width=0.75in]{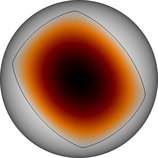}
    \end{tabular}
    and
    \begin{tabular}{c}
      2 pentagons \\ \includegraphics[width=0.75in]{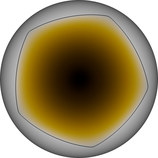}
    \end{tabular} &
    Quadrilateral eigenvalue: $9.988\: (1.63 \cdot 10^{-3})$

    Pentagon eigenvalue: $8.298\: (7.50 \cdot 10^{-5})$

    $S_\eps$: $0.322$
    
    Total energy: $66.859$ \\
    \hline
    8 &
    \begin{tabular}{c}
      4 quadrilaterals \\ \includegraphics[width=0.75in]{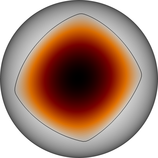} 
    \end{tabular}
    and
    \begin{tabular}{c}
      4 pentagons \\ \includegraphics[width=0.75in]{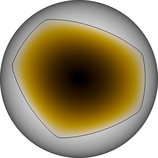}
    \end{tabular} &

    Quadrilateral eigenvalue: $11.380\: (5.31 \cdot 10^{-3})$
    
    Pentagon eigenvalue: $10.230\: (2.91 \cdot 10^{-3})$

    $S_\eps$: $0.650$

    Total energy: $87.102$ \\
    \hline
    16 &
    \begin{tabular}{c}
      8 A and 4 B pentagons\\ \includegraphics[width=0.75in]{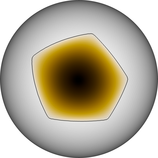} \includegraphics[width=0.75in]{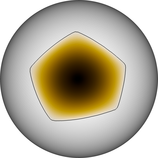}
    \end{tabular}
    and
    \begin{tabular}{c}
      4 hexagons \\ \includegraphics[width=0.75in]{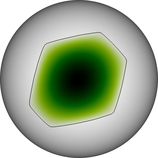} 
    \end{tabular}
    &
    Pentagon (A) eigenvalue: $22.647\: (1.05 \cdot 10^{-2})$

    Pentagon (B) eigenvalue: $23.610\: (2.43 \cdot 10^{-2})$
    
    Hexagon eigenvalue: $20.496\: (1.05 \cdot 10^{-2})$

    $S_\eps$: $1.264$

    Total energy: $362.718$ \\
    \hline
    32 &
    \begin{tabular}{c}
      12 pentagons \\ \includegraphics[width=0.75in]{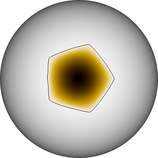}
    \end{tabular}
    and
    \begin{tabular}{c}
      20 hexagons \\ \includegraphics[width=0.75in]{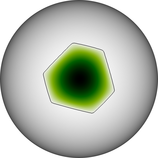}
    \end{tabular}
    &
    Pentagon eigenvalue: $48.436\: (1.46 \cdot 10^{-1})$

    Hexagon eigenvalue: $44.460\: (1.24 \cdot 10^{-1})$

    $S_\eps$: $2.496$

    Total energy: $1472.920$ \\
    \hline
  \end{tabular}

  \caption{More details of optimal partitions. In the small plots, we plot the corresponding $u_i^{\eps,h}$ with a black contour at $v_i^{\eps,h} = 0$.}
  \label{tab:many-m}
\end{table}

There are several striking features:
\begin{itemize}
\item All partitions consist of curvi-linear polygons;

\item The boundary of each partition consists of arcs with zero geodesic curvature (``straight lines'');
\item Each junction is a triple junction with an equal angle condition satisfied;
\item There are at most two types of polygon in the partition;
\item In the case of two different polygons, the polygon with more sides has lower eigenvalue;
\item As $m$ increases the number of edges in each polygon increases;
\item Each polygon has at most $6$ edges.
\end{itemize}

We define the dual polygon to a partition by considering the edges and vertices as a graph and taking the dual graph. In our case, since we always have triple junctions this defines a triangulation of the sphere. Let $V$ be the number of vertices, $E$ the number of edges and $F$ the number of faces in the dual polygon to a partition $\{ \Gamma_i \}_{i=1}^m$. We know that this will satisfy Euler's identity, $V - E + F = \chi$, where $\chi$ is the Euler characteristic ($2$ in the case of a sphere), and also that
\begin{equation*}
  2 E = \sum_{k=0}^\infty k n_k, \quad
  3 F = \sum_{k=0}^\infty k n_k, \quad
  V = \sum_{k=0}^\infty n_k,
\end{equation*}
where $n_k$ is the degree of a vertex in the dual polygon. The degree of a vertex is equal to the number of edges of the corresponding partition. Using these equations in Euler's identity gives
\begin{equation}
  \label{eq:gauss-bonnet}
  4 n_2 + 3 n_3 + 2 n_4 + n_5 = 6 \chi + \sum_{k=7}^\infty ( k - 6 ) n_k.
\end{equation}
This result is a special case of the Gauss-Bonnet theorem. We can think of this result as saying that polygons with less than six sides correspond to regions of positive Gauss curvature, hexagons correspond to zero Gauss curvature and polygons with more than six sides correspond to negative Gauss curvature.

This identity is consistent with the partitions in Table~\ref{tab:many-m}. Our computations suggest that the polygonal structure of the optimal partition consists of polygons with six or less sides. This agrees with the idea that the sphere has uniform positive Gauss curvature. We can deduce that if an $m$-partition of the sphere consists of only pentagons and hexagons, then there will be $12$ pentagons and $m-12$ hexagons. We expect this to be the optimal partition for large values of $m$.

%%%%%%%%%%%%%%%%%%%%%%%%%%%
\subsection{Computed partitions of other surfaces}
%%%%%%%%%%%%%%%%%%%%%%%%%%%

We consider two other surfaces to see if these conclusions persist on a large class of surfaces. The first example, surface ({\it D}), is taken from the work of \citet{Dzi88} where the surface is given by $\Gamma = \{ x \in \R^3 : \Phi(x) = 0\}$ for $\Phi$ given by
\begin{equation*}
  \Phi( x_1, x_2, x_3 ) := ( x_1 - x_3^2 )^2 + x_2^2 + x_3^2 - 1.
\end{equation*}
This has the same genus as a sphere but has large changes in curvature. The second example is given by a torus ({\it T}) with inner radius $0.6$ and outer radius $1$. This has different genus to the sphere. We proceed with the same refinement strategy as on the sphere. Details of the parameters are given in Table~\ref{tab:dofs-other}.

\begin{table}
  \centering
  \begin{tabular}{|c|ccc|ccc|}
    \hline
    $m$ & & {\raggedleft Surface (\emph{D})} & & & Torus & \\
     & $l$ & Degrees of freedom & $\eps$ & $l$ & Degrees of freedom & $\eps$ \\
    \hline
    3 & 12 & $311\,982$ & $3.125 \cdot 10^{-4}$ &
        12 & $393\,216$ & $3.125 \cdot 10^{-4}$ \\
    4 & 12 & $415\,976$ & $3.125 \cdot 10^{-4}$ &
        12 & $524\,288$ & $3.125 \cdot 10^{-4}$ \\
    5 & 10 & $256\,365$ & $6.25 \cdot 10^{-4}$ &
        10 & $326\,680$ & $6.25 \cdot 10^{-4}$ \\
    6 & 9 & $150\,900$ & $8.883 \cdot 10^{-4}$ &
        10 & $393\,216$ & $6.25 \cdot 10^{-4}$ \\
    7 & 9 & $176\,050$ & $8.883 \cdot 10^{-4}$ &
        10 & $458\,752$ & $6.25 \cdot 10^{-4}$ \\
    8 & 9 & $201\,200$ & $8.883 \cdot 10^{-4}$ &
        10 & $524\,288$ & $6.25 \cdot 10^{-4}$ \\
    16 & 9 & $402\,400$ & $8.883 \cdot 10^{-4}$ &
        10 & $1\,048\,576$ & $6.25 \cdot 10^{-4}$ \\
    32 & 9 & $804\,800$ & $8.883 \cdot 10^{-4}$ &
         8 & $1\,045\,696$ & $8.883 \cdot 10^{-4}$ \\
    \hline
  \end{tabular}
  
  \caption{Final parameters for computations on the surface (\emph{D}) and the torus.}
  \label{tab:dofs-other}
\end{table}

\begin{figure}
  \centering
  \includegraphics{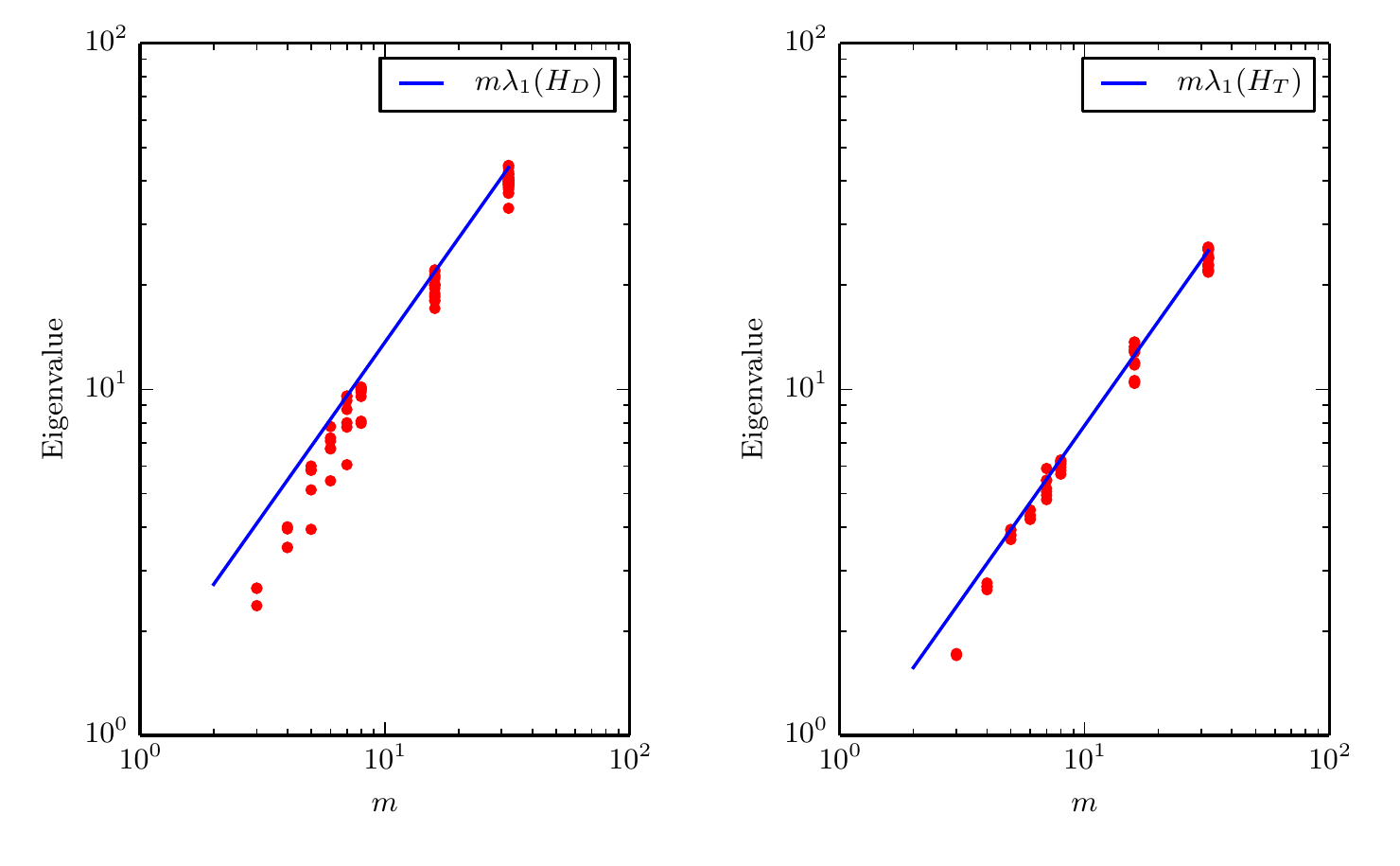}

  \caption{Plot of the eigenvalues at different values of $m$. Left for the surface (\emph{D}) and right for the torus (\emph{T}). The blue line indicates the scaled eigenvalue corresponding to a hexagon $H$ of equal area to each surface -- this is the conjectured average eigenvalue for large $m$ in the plane \citep{CafLin07}.}
  \label{fig:many-m-other}
\end{figure}

We plot for the eigenvalues corresponding to the optimal partition in Figure~\ref{fig:many-m-other}. We compute the eigenvalue as the $H^1(\Gamma)$ semi-norm of each component. We have also included the line at $m \lambda_1( H_D )$ and $m \lambda_1( H_T )$ in each plot, where $H_{D}$ and $H_{T}$ are the regular hexagons with area equal to the surface ({\it D}) and the torus ({\it T}). We do not have direct access to the eigenvalues on either of these surfaces so do not add that to this plot.

\begin{figure}[p]
  \begin{subfigure}[b]{0.3\textwidth}
    \includegraphics[width=\textwidth]{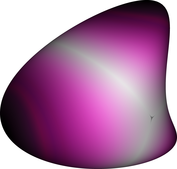}
    \caption{$m=3$, 3 lens (pink)}
  \end{subfigure}
  \hfill
  \begin{subfigure}[b]{0.3\textwidth}
    \includegraphics[width=\textwidth]{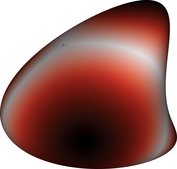}
    \caption{$m=4$, 4 triangles (red)}
  \end{subfigure}
  \hfill
  \begin{subfigure}[b]{0.3\textwidth}
    \includegraphics[width=\textwidth]{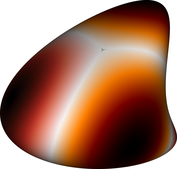}
    \caption{$m=5$, 2 triangles (red) and 3 quadrilaterals (orange)}
  \end{subfigure}
  
  \begin{subfigure}[b]{0.3\textwidth}
    \includegraphics[width=\textwidth]{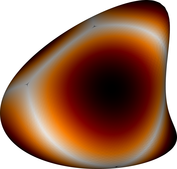}
    \caption{$m=6$, 6 quadrilaterals (orange)}
  \end{subfigure}
  \hfill
  \begin{subfigure}[b]{0.3\textwidth}
    \includegraphics[width=\textwidth]{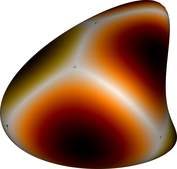}
    \caption{$m=7$, 5 quadrilaterals (orange) and 2 pentagons (yellow)}
  \end{subfigure}
  \hfill
  \begin{subfigure}[b]{0.3\textwidth}
    \includegraphics[width=\textwidth]{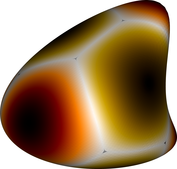}
    \caption{$m=8$, 4 quadrilaterals (orange) and 4 pentagons (yellow)}
  \end{subfigure}

  \hspace*{\fill}
  \begin{subfigure}[b]{0.3\textwidth}
    \includegraphics[width=\textwidth]{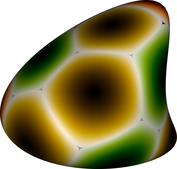}
    \caption{$m=16$, 3 quadrilateral (orange), 6 pentagons (yellow), 7 hexagons (green)}
  \end{subfigure}
  \hfill
  \begin{subfigure}[b]{0.3\textwidth}
    \includegraphics[width=\textwidth]{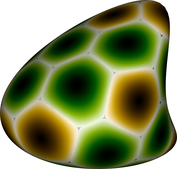}
    \caption{$m=32$, 12 pentagons (yellow) and 20 hexagons (green).}
  \end{subfigure}
  \hspace*{\fill}

  \caption{Plots of the minimising configurations on the surface (\emph{D}). Same colouring as Figure~\ref{fig:many-m}}
  \label{fig:many-m-dziuk}
\end{figure}

\begin{figure}[p]
  \begin{subfigure}[b]{0.3\textwidth}
    \includegraphics[width=\textwidth]{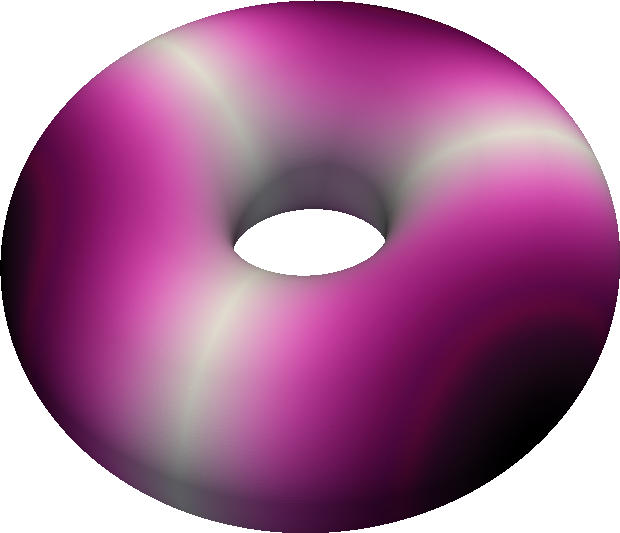}
    \caption{$m=3$, 3 cylinders (pink)}
  \end{subfigure}
  \hfill
  \begin{subfigure}[b]{0.3\textwidth}
    \includegraphics[width=\textwidth]{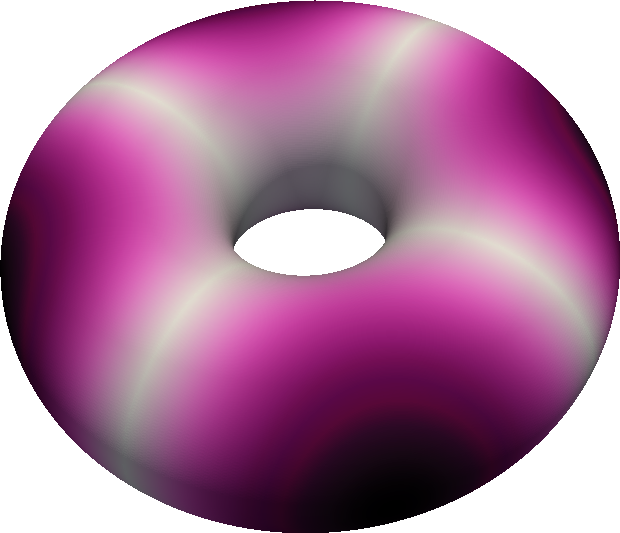}
    \caption{$m=4$, 4 cylinders (pink)}
  \end{subfigure}
  \hfill
  \begin{subfigure}[b]{0.3\textwidth}
    \includegraphics[width=\textwidth]{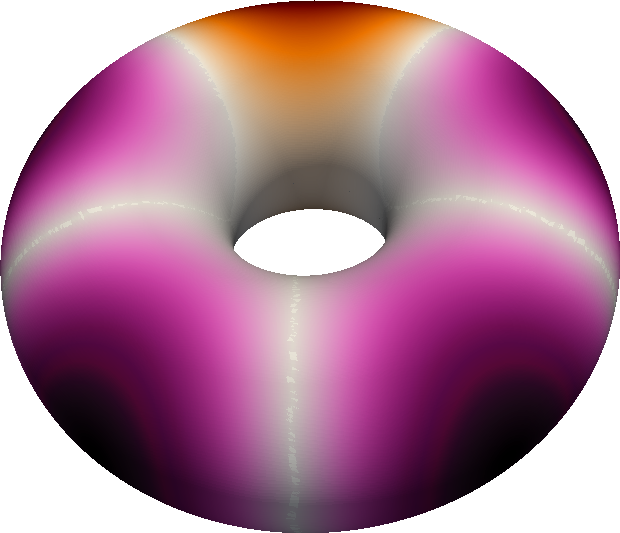}
    \caption{$m=5$, 4 two-sided shapes (pink) and 1 quadrilateral (orange)}
  \end{subfigure}

  \begin{subfigure}[b]{0.3\textwidth}
    \includegraphics[width=\textwidth]{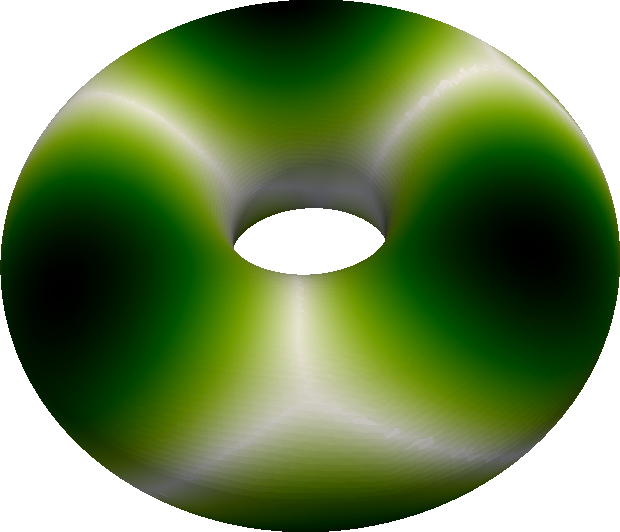}
    \caption{$m=6$, 6 hexagons (green)}
  \end{subfigure}
  \hfill
  \begin{subfigure}[b]{0.3\textwidth}
    \includegraphics[width=\textwidth]{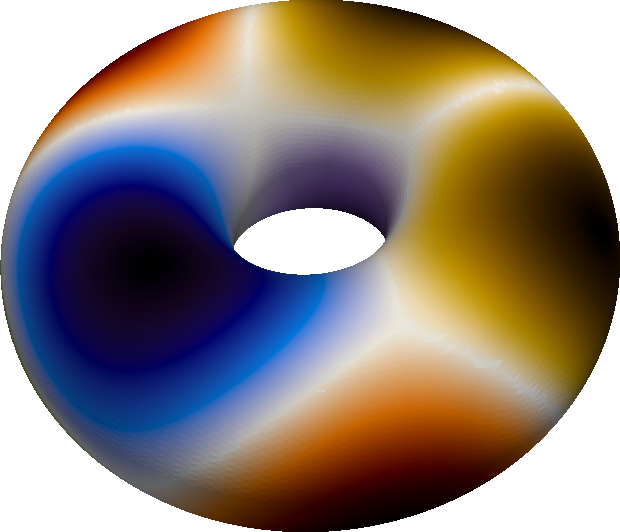}
    \caption{$m=7$, 2 quadrilaterals (orange), 2 pentagons (yellow), 1 hexagon (green), 1 octagon (blue), 1 decagon (purple)}
  \end{subfigure}
  \hfill
  \begin{subfigure}[b]{0.3\textwidth}
    \includegraphics[width=\textwidth]{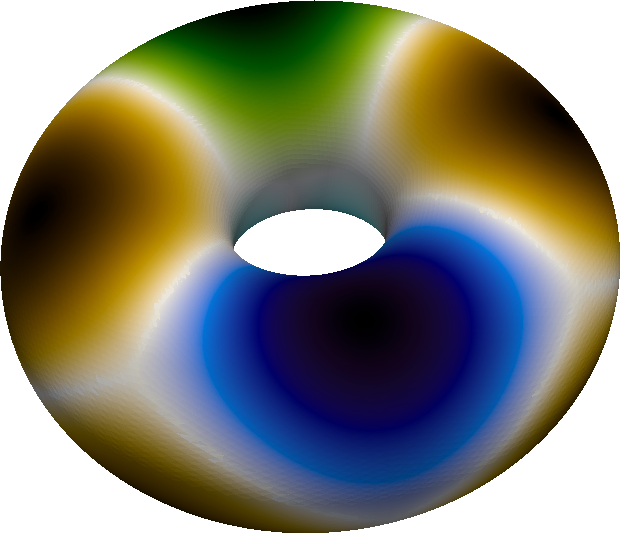}
    \caption{$m=8$, 4 pentagons (yellow), 1 hexagon (green), 2 heptagons (cyan), 1 octagon (blue)}
  \end{subfigure}

  \hspace*{\fill}
  \begin{subfigure}[b]{0.3\textwidth}
    \includegraphics[width=\textwidth]{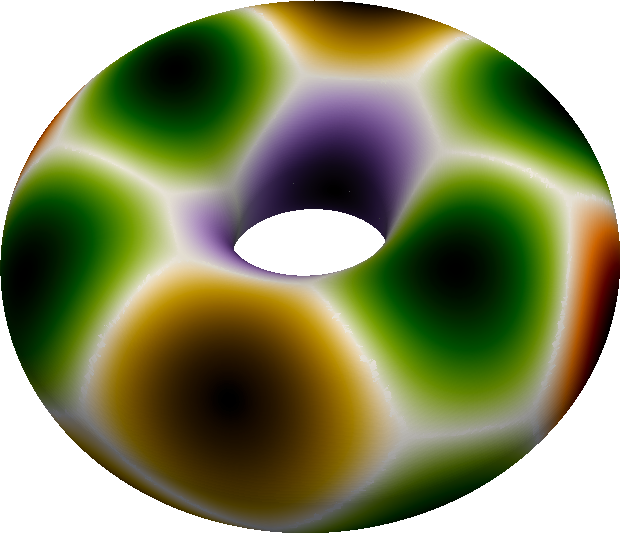}
    \caption{$m=16$, 2 quadrilaterals (orange), 4 pentagons (yellow), 8 hexagons (green) and 2 decagons (purple)}
  \end{subfigure}
  \hfill
  \begin{subfigure}[b]{0.3\textwidth}
    \includegraphics[width=\textwidth]{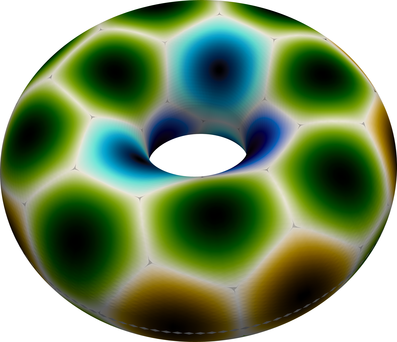}
    \caption{$m=32$, 8 pentagons (yellow), 18 hexagons (green), 4 heptagons (cyan) and 2 octagons (blue).}
  \end{subfigure}
  \hspace*{\fill}

  \caption{Plots of the minimising configurations on the torus. Same colouring as Figure~\ref{fig:many-m}}
  \label{fig:many-m-torus}
\end{figure}

For surface (\emph{D}), we plot the optimal configurations in Figure~\ref{fig:many-m-dziuk} with more details given, including eigenvalues and energy, in Table~\ref{tab:many-m-dziuk}. For the torus, we plot the optimal configurations in Figure~\ref{fig:many-m-torus} with more details given, include eigenvalues and energy, in Table~\ref{tab:many-m-torus}.

\begin{table}[p]
  \centering
  \newcolumntype{C}[1]{>{\centering\let\newline\\\arraybackslash\hspace{0pt}}m{#1}}
  \begin{tabular}{|C{0.025\textwidth}|C{0.9725\textwidth}|}
    \hline
    $m$ & Partition \\
    \hline
    3 & %
    \begin{tabular}{c}
      lens \\ \includegraphics[width=1in,height=0.9in,keepaspectratio]{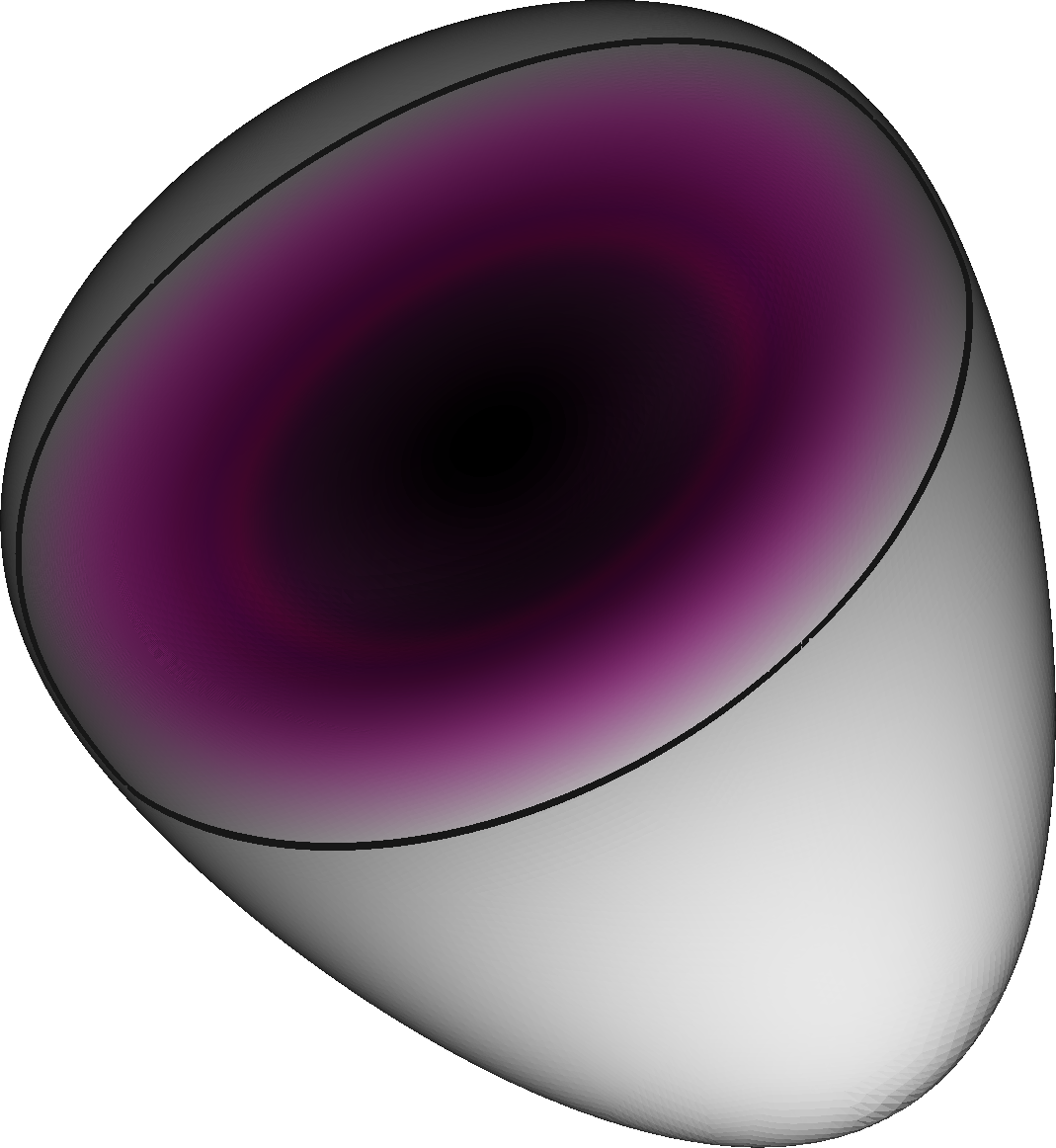} \\ 2.664
    \end{tabular}
    \begin{tabular}{c}
      crescent \\ \includegraphics[width=1in,height=0.9in,keepaspectratio]{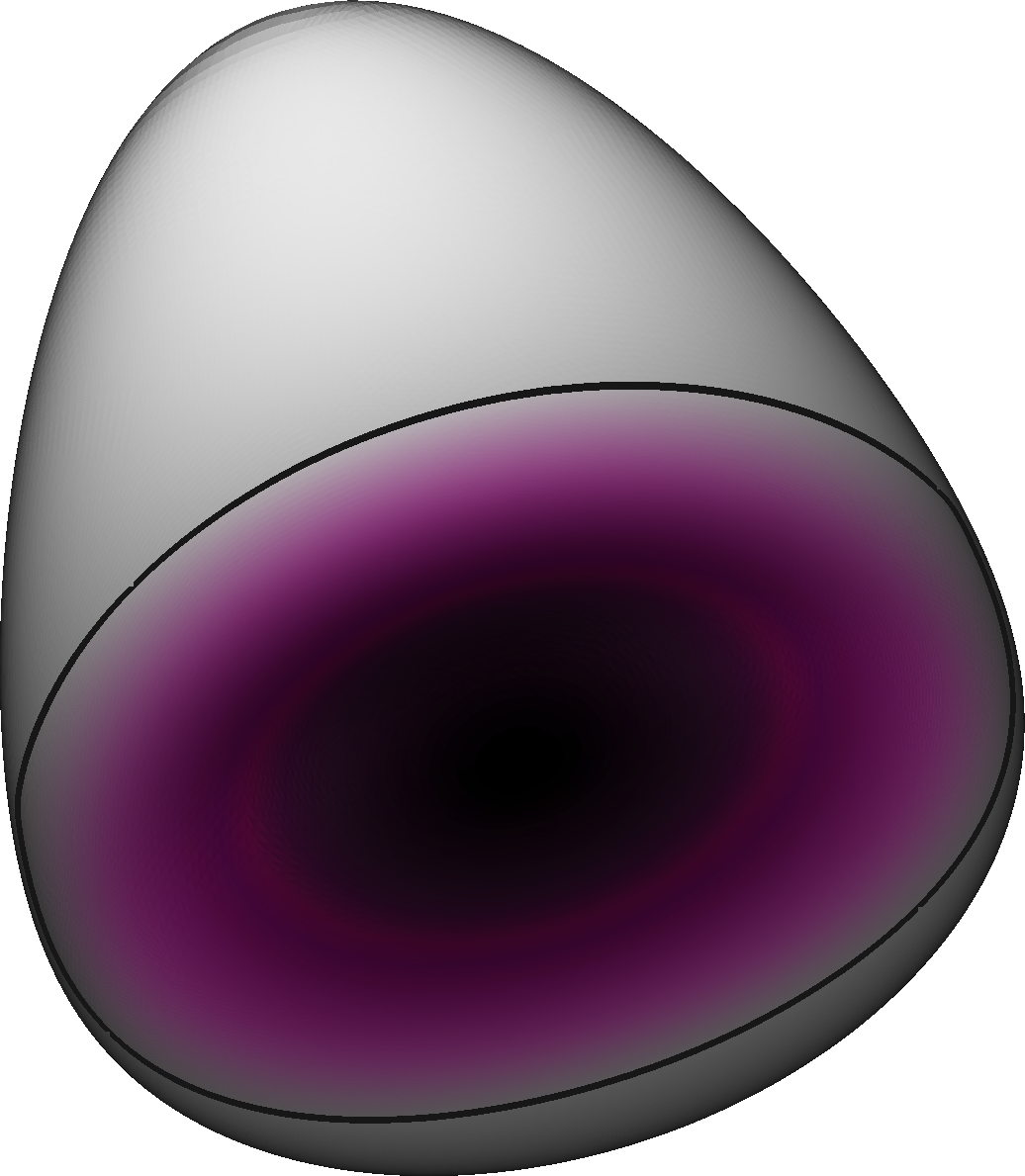} \\ 2.664
    \end{tabular}
    \begin{tabular}{c}
      crescent \\ \includegraphics[width=1in,height=0.9in,keepaspectratio]{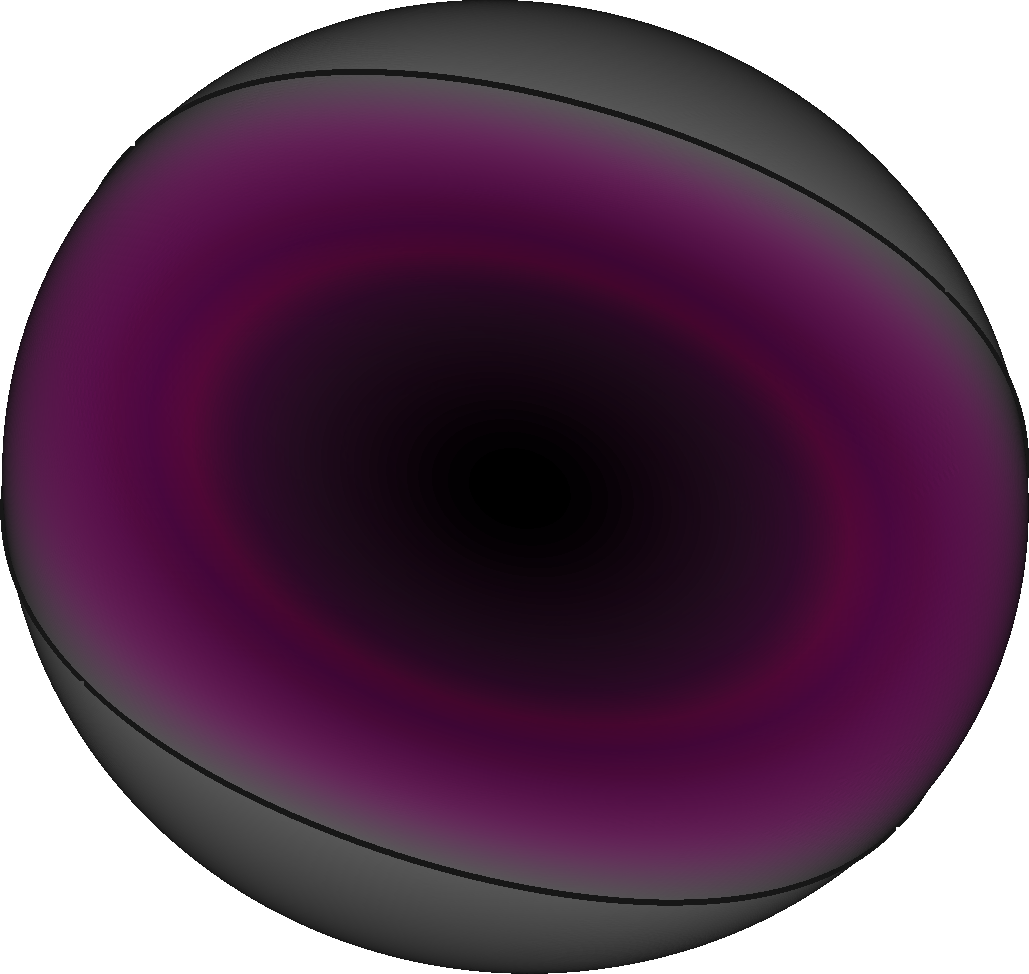} \\ 2.372
    \end{tabular}
    
    $S_\eps$: $0.040$

    Total energy: $7.741$
    \\
    \hline
    4 & %
    \begin{tabular}{c}
      triangle \\ \includegraphics[width=1in,height=0.9in,keepaspectratio]{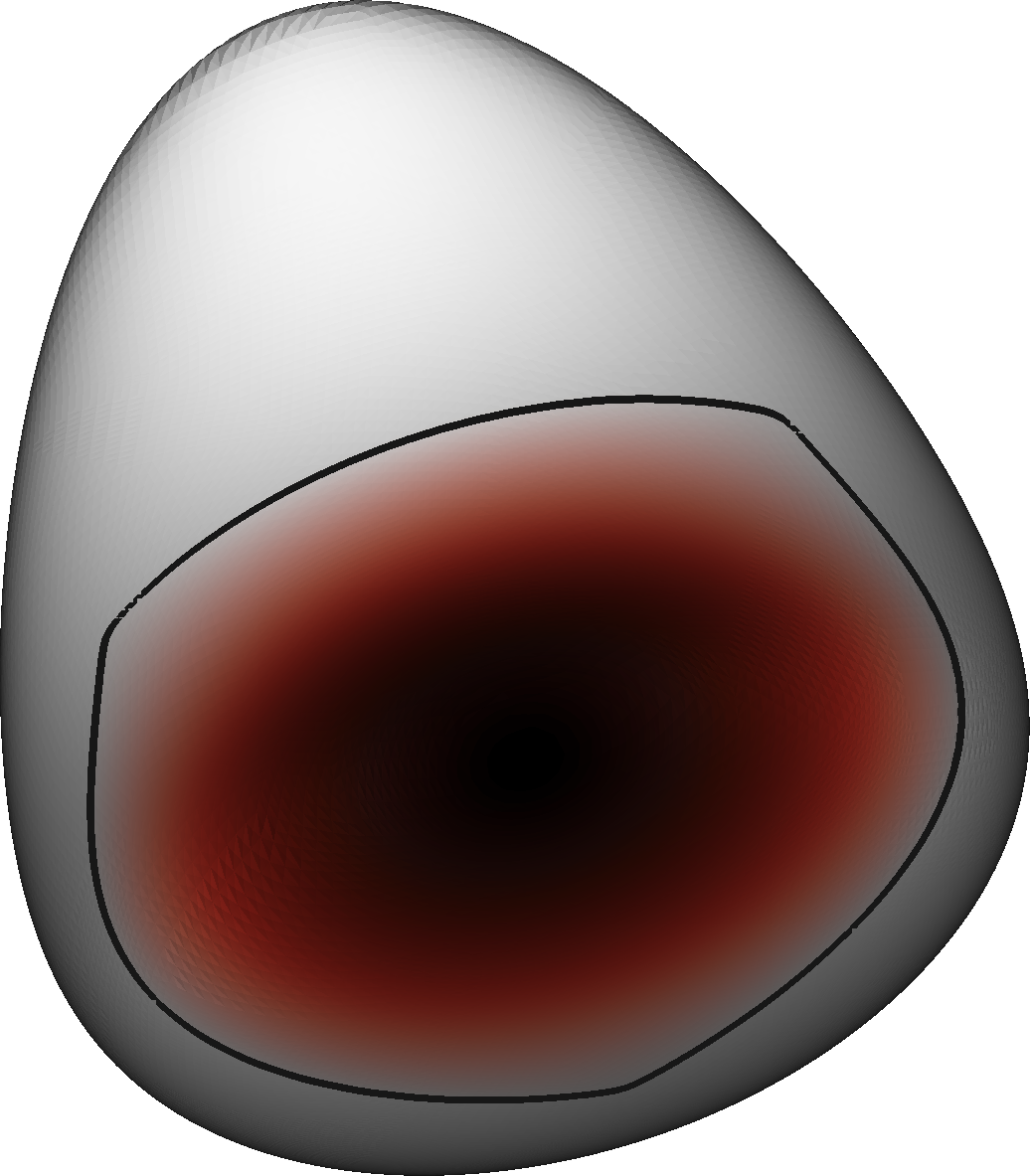} \\ 3.493
    \end{tabular}
    \begin{tabular}{c}
      triangle \\ \includegraphics[width=1in,height=0.9in,keepaspectratio]{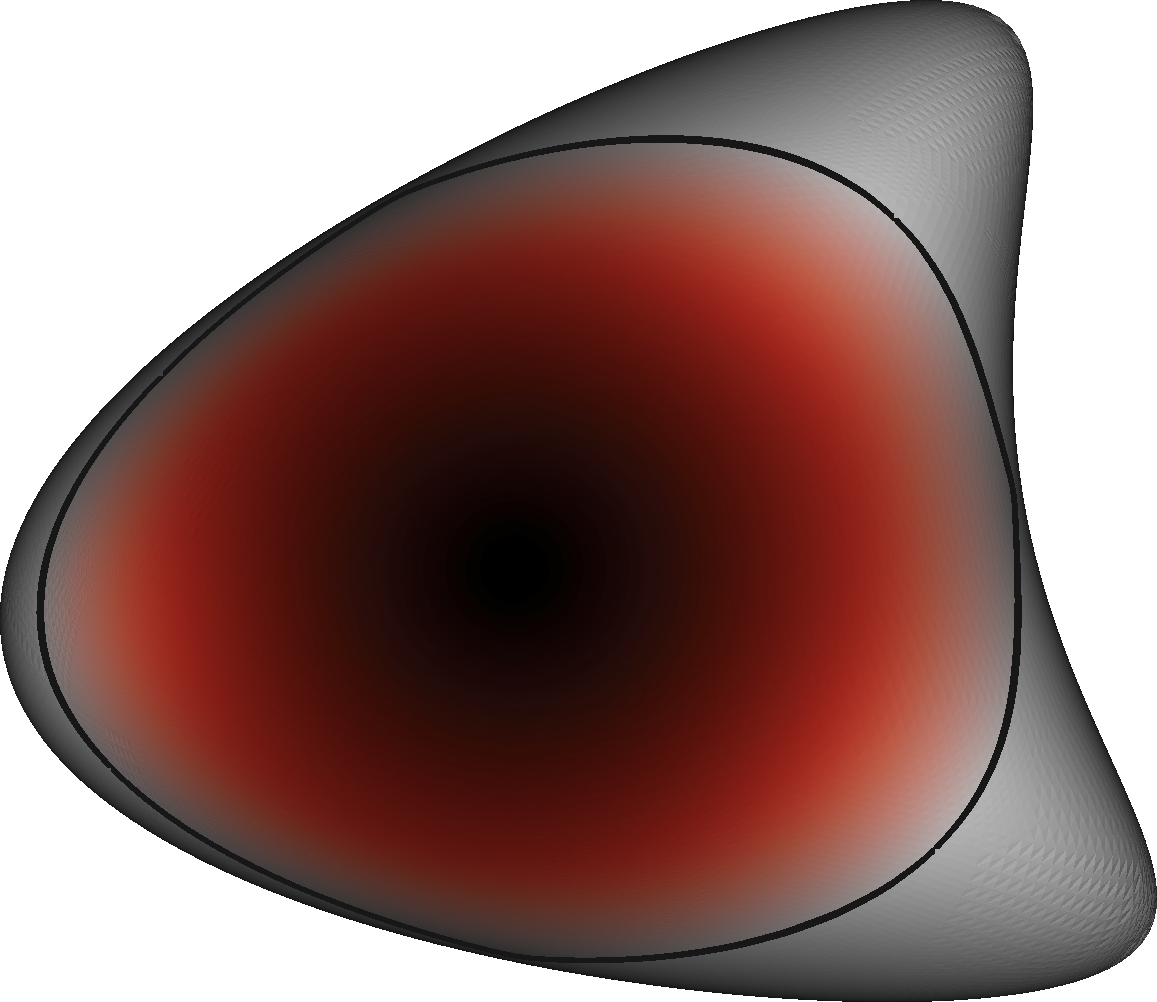} \\ 3.494
    \end{tabular}
    \begin{tabular}{c}
      triangle \\ \includegraphics[width=1in,height=0.9in,keepaspectratio]{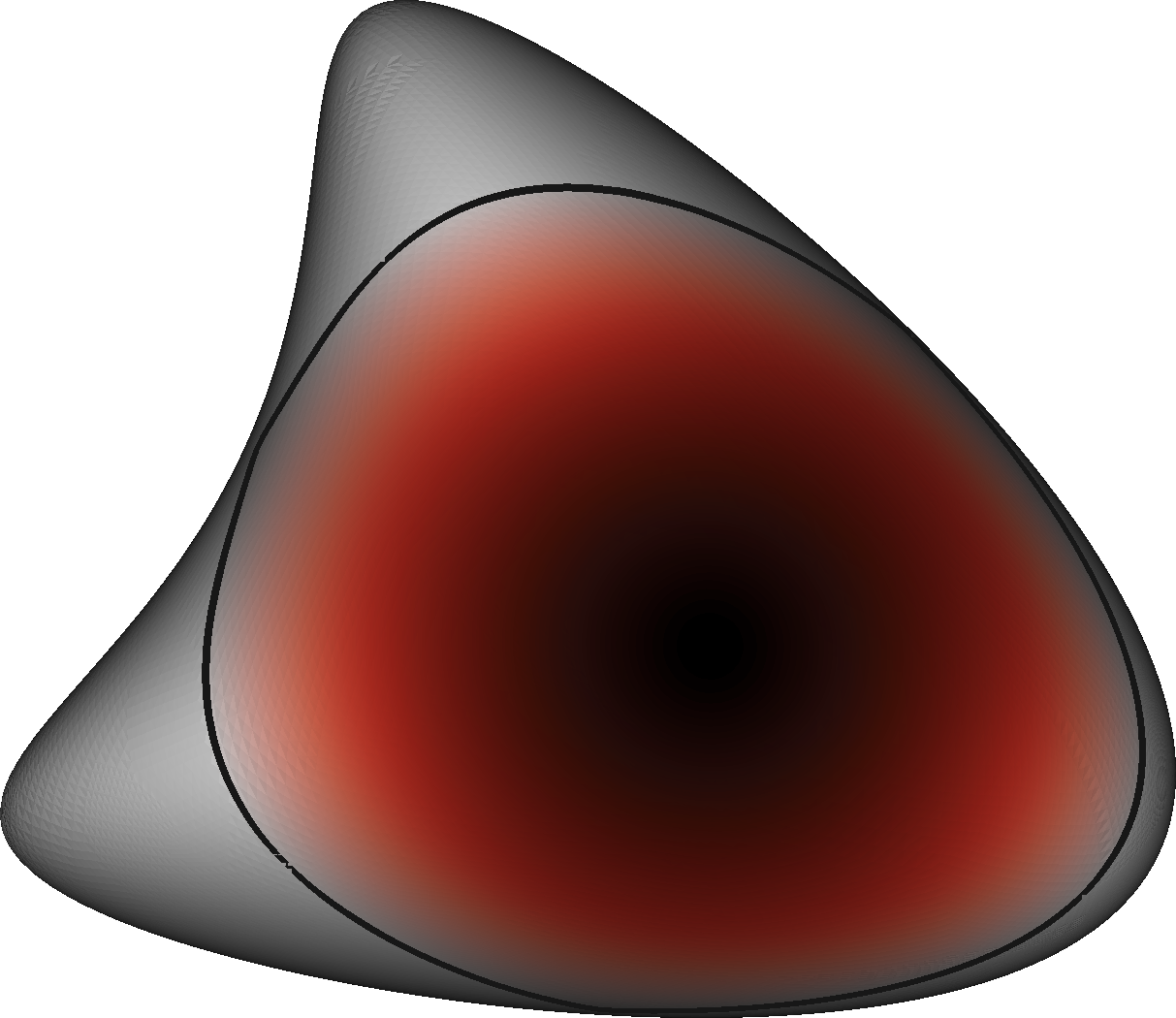} \\ 4.008
    \end{tabular}
    \begin{tabular}{c}
      triangle \\ \includegraphics[width=1in,height=0.9in,keepaspectratio]{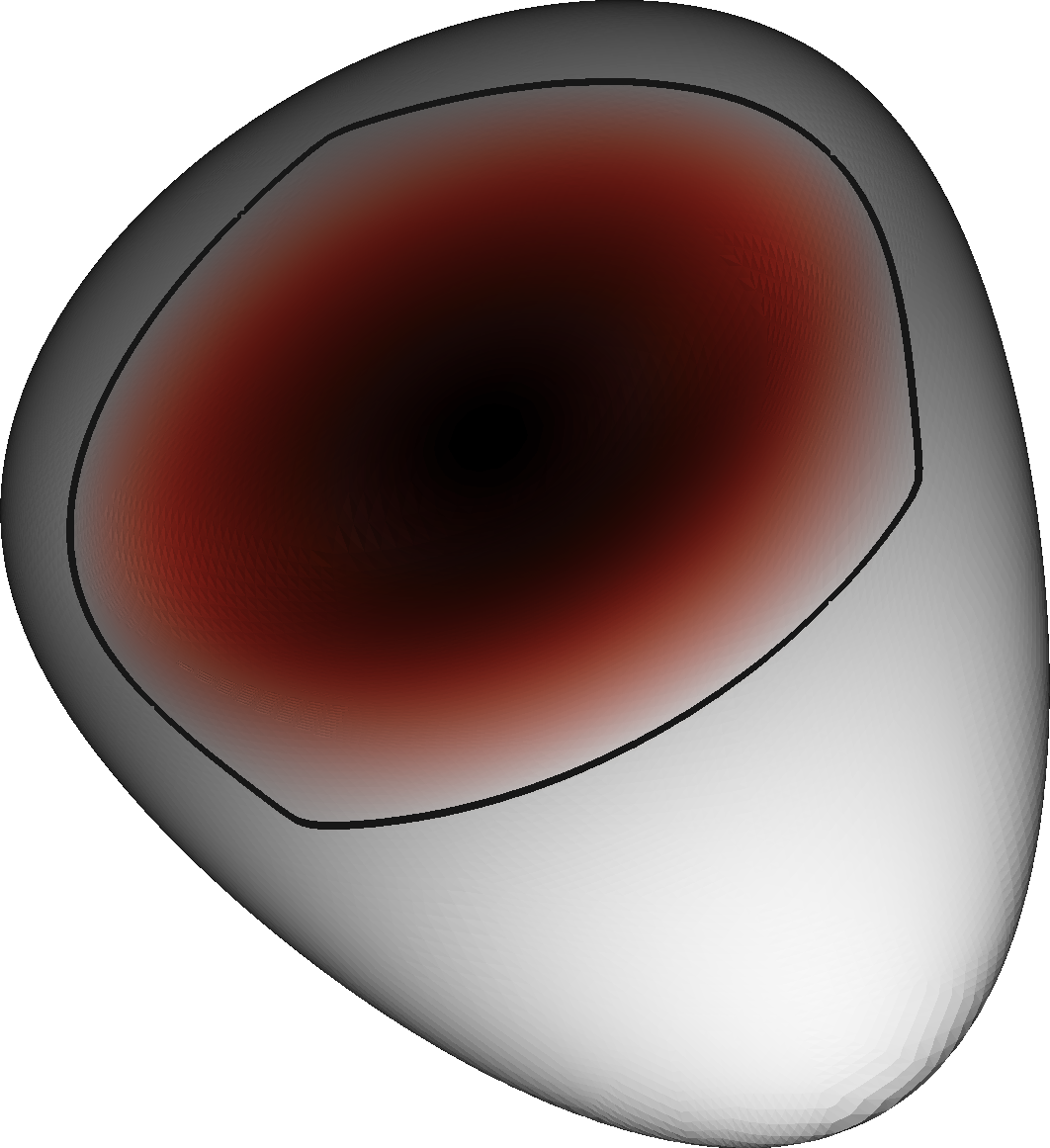} \\ 3.952
    \end{tabular}

    $S_\eps$: $0.103$

    Total energy: $15.051$
    \\
    \hline
    5 & %
    \begin{tabular}{c}
      triangle \\ \includegraphics[width=1in,height=0.9in,keepaspectratio]{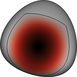} \\ 5.843
    \end{tabular}
    \begin{tabular}{c}
      triangle \\ \includegraphics[width=1in,height=0.9in,keepaspectratio]{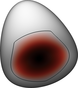} \\ 5.125
    \end{tabular}
    \begin{tabular}{c}
      quadrilateral \\ \includegraphics[width=1in,height=0.9in,keepaspectratio]{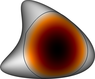} \\ 6.004
    \end{tabular}
    \begin{tabular}{c}
      quadrilateral \\ \includegraphics[width=1in,height=0.9in,keepaspectratio]{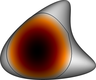} \\ 5.944
    \end{tabular}
    \begin{tabular}{c}
      quadrilateral \\ \includegraphics[width=1in,height=0.9in,keepaspectratio]{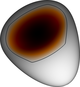} \\ 3.942
    \end{tabular}

    $S_\eps$: $0.312852$

    Total energy: $27.072$
    \\
    \hline
    6 & %
    \begin{tabular}{c}
      quadrilateral \\ \includegraphics[width=1in,height=0.9in,keepaspectratio]{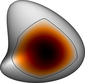} \\ 7.808
    \end{tabular}
    \begin{tabular}{c}
      quadrilateral \\ \includegraphics[width=1in,height=0.9in,keepaspectratio]{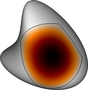} \\ 7.241
    \end{tabular}
    \begin{tabular}{c}
      quadrilateral \\ \includegraphics[width=1in,height=0.9in,keepaspectratio]{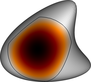} \\ 7.093
    \end{tabular}
    \begin{tabular}{c}
      quadrilateral \\ \includegraphics[width=1in,height=0.9in,keepaspectratio]{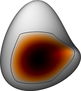} \\ 6.753
    \end{tabular}
    \begin{tabular}{c}
      quadrilateral \\ \includegraphics[width=1in,height=0.9in,keepaspectratio]{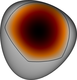} \\ 6.730
    \end{tabular}
    \begin{tabular}{c}
      quadrilateral \\ \includegraphics[width=1in,height=0.9in,keepaspectratio]{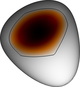} \\ 5.443
    \end{tabular}

    $S_\eps$: $0.753$

    Total energy: $41.821$ \\
    \hline
  \end{tabular}
\end{table}

\begin{table}[p]
  \centering
  \newcolumntype{C}[1]{>{\centering\let\newline\\\arraybackslash\hspace{0pt}}m{#1}}
  \begin{tabular}{|C{0.025\textwidth}|C{0.9725\textwidth}|}
    \hline
    7 & %
    \begin{tabular}{c}
      quadrilateral \\ \includegraphics[width=1in,height=0.9in,keepaspectratio]{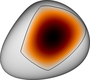} \\ 9.569
    \end{tabular}
    \begin{tabular}{c}
      quadrilateral \\ \includegraphics[width=1in,height=0.9in,keepaspectratio]{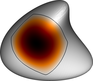} \\ 9.556
    \end{tabular}
    \begin{tabular}{c}
      quadrilateral \\ \includegraphics[width=1in,height=0.9in,keepaspectratio]{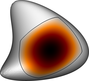} \\ 9.275
    \end{tabular}
    \begin{tabular}{c}
      quadrilateral \\ \includegraphics[width=1in,height=0.9in,keepaspectratio]{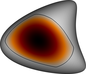} \\ 8.748
    \end{tabular}
    \begin{tabular}{c}
      quadrilateral \\ \includegraphics[width=1in,height=0.9in,keepaspectratio]{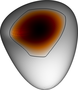} \\ 7.780
    \end{tabular}

    \begin{tabular}{c}
      pentagon \\ \includegraphics[width=1in,height=0.9in,keepaspectratio]{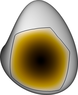} \\ 8.009
    \end{tabular}
    \begin{tabular}{c}
      pentagon \\ \includegraphics[width=1in,height=0.9in,keepaspectratio]{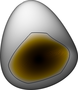} \\ 6.058
    \end{tabular}

    $S_\eps$: $1.102$

    Total energy: $60.096$
    \\
    \hline
    8 & %
    \begin{tabular}{c}
      quadrilateral \\ \includegraphics[width=1in,height=0.9in,keepaspectratio]{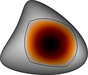} \\ 10.1602
    \end{tabular}
    \begin{tabular}{c}
      quadrilateral \\ \includegraphics[width=1in,height=0.9in,keepaspectratio]{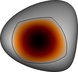} \\ 9.83384
    \end{tabular}
    \begin{tabular}{c}
      quadrilateral \\ \includegraphics[width=1in,height=0.9in,keepaspectratio]{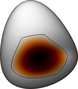} \\ 8.09237
    \end{tabular}
    \begin{tabular}{c}
      quadrilateral \\ \includegraphics[width=1in,height=0.9in,keepaspectratio]{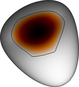} \\ 7.978
    \end{tabular}

    \begin{tabular}{c}
      pentagon \\ \includegraphics[width=1in,height=0.9in,keepaspectratio]{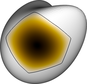} \\ 10.128
    \end{tabular}
    \begin{tabular}{c}
      pentagon \\ \includegraphics[width=1in,height=0.9in,keepaspectratio]{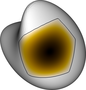} \\ 10.034
    \end{tabular}
    \begin{tabular}{c}
      pentagon \\ \includegraphics[width=1in,height=0.9in,keepaspectratio]{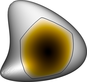} \\ 9.965
    \end{tabular}
    \begin{tabular}{c}
      pentagon \\ \includegraphics[width=1in,height=0.9in,keepaspectratio]{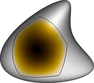} \\ 9.539
    \end{tabular}

    $S_\eps$: $1.63602$

    Total energy: $77.367$
    \\
    \hline
  \end{tabular}

  \caption{More details of optimal partitions on the surface (\emph{D}). In the small plots, we plot the corresponding $u_i^{\eps,h}$ with a black contour at $v_i^{\eps,h} = 0$.}
  \label{tab:many-m-dziuk}
\end{table}

\begin{table}[p]
  \centering
 \newcolumntype{C}[1]{>{\centering\let\newline\\\arraybackslash\hspace{0pt}}m{#1}}
  \begin{tabular}{|C{0.025\textwidth}|C{0.9725\textwidth}|}
    \hline
    $m$ & Partition \\
    \hline
    3 & %
    \begin{tabular}{c}
      cylinder \\ \includegraphics[width=1in,height=0.9in,keepaspectratio]{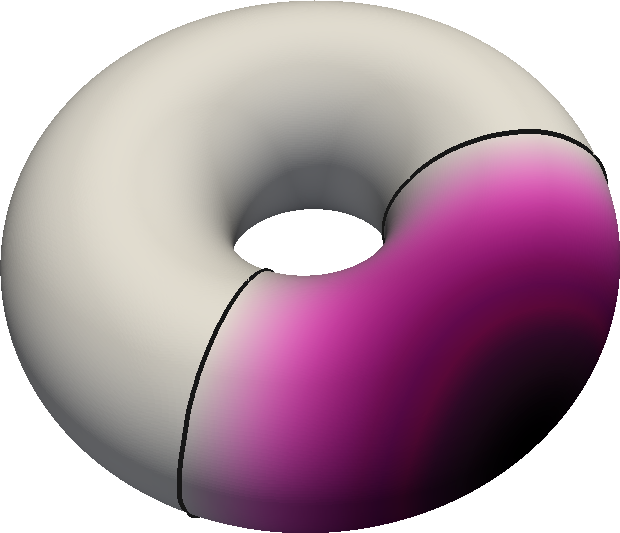} \\ 1.725
    \end{tabular}
    \begin{tabular}{c}
      cylinder \\ \includegraphics[width=1in,height=0.9in,keepaspectratio]{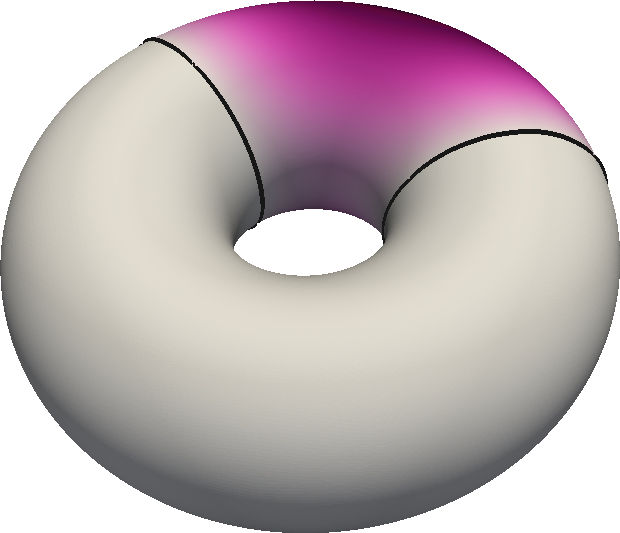} \\ 1.703
    \end{tabular}
    \begin{tabular}{c}
      cylinder \\ \includegraphics[width=1in,height=0.9in,keepaspectratio]{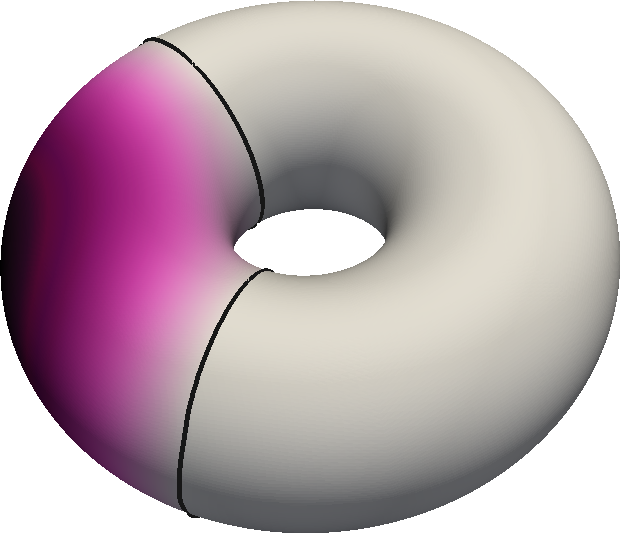} \\ 1.717
    \end{tabular}
    
    $S_\eps$: $1.207$

    Total energy: $6.353$
    \\
    \hline
    4 & %
    \begin{tabular}{c}
      cylinder \\ \includegraphics[width=1in,height=0.9in,keepaspectratio]{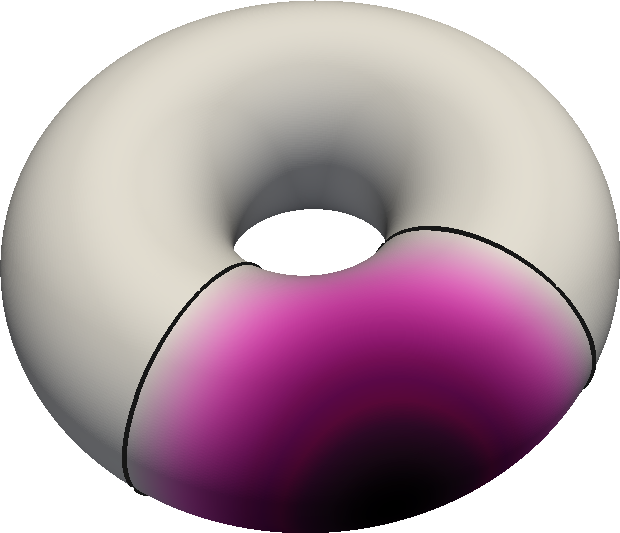} \\ 2.758
    \end{tabular}
    \begin{tabular}{c}
      cylinder \\ \includegraphics[width=1in,height=0.9in,keepaspectratio]{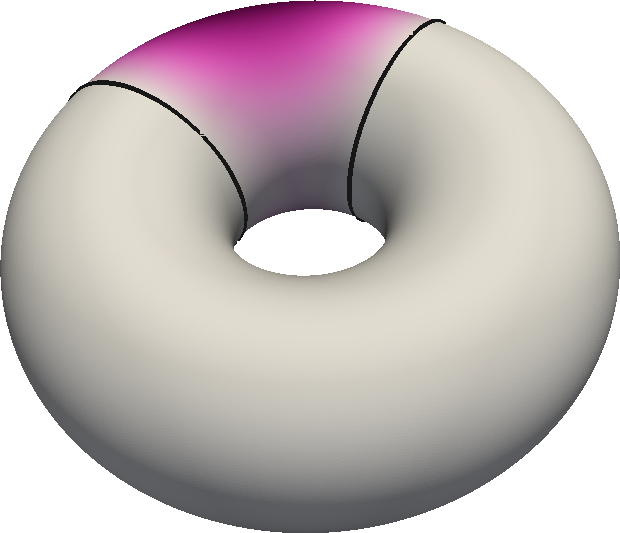} \\ 2.637
    \end{tabular}
    \begin{tabular}{c}
      cylinder \\ \includegraphics[width=1in,height=0.9in,keepaspectratio]{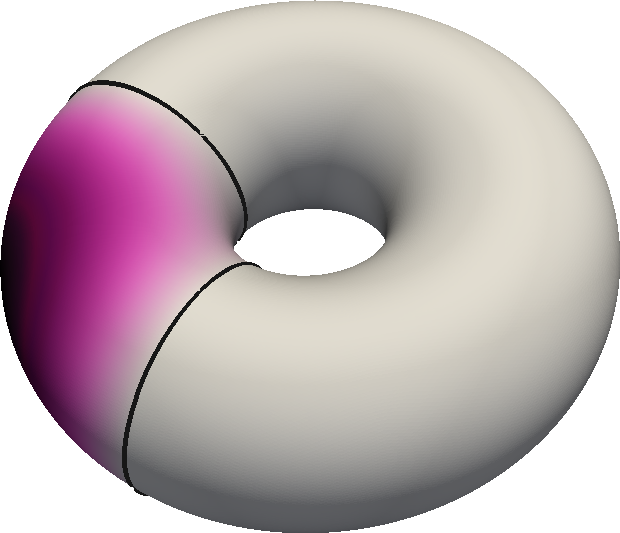} \\ 2.595
    \end{tabular}
    \begin{tabular}{c}
      cylinder \\ \includegraphics[width=1in,height=0.9in,keepaspectratio]{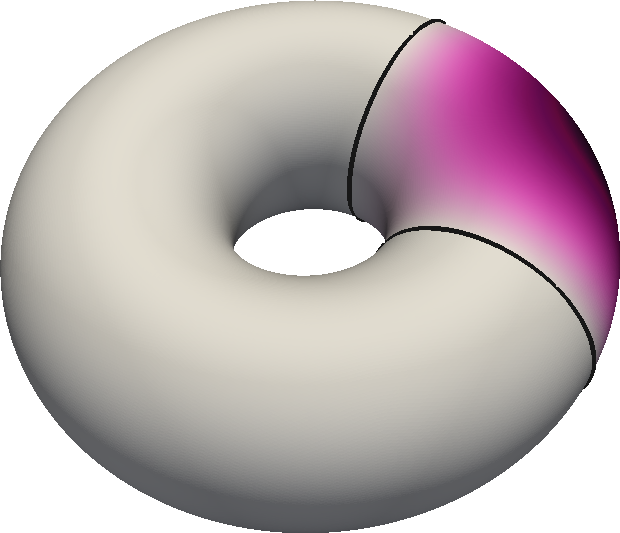} \\ 2.595
    \end{tabular}
    
    $S_\eps$: $0.106$

    Total energy: $10.890$
    \\
    \hline
    5 & %
    \begin{tabular}{c}
      two sided shape \\ \includegraphics[width=1in,height=0.9in,keepaspectratio]{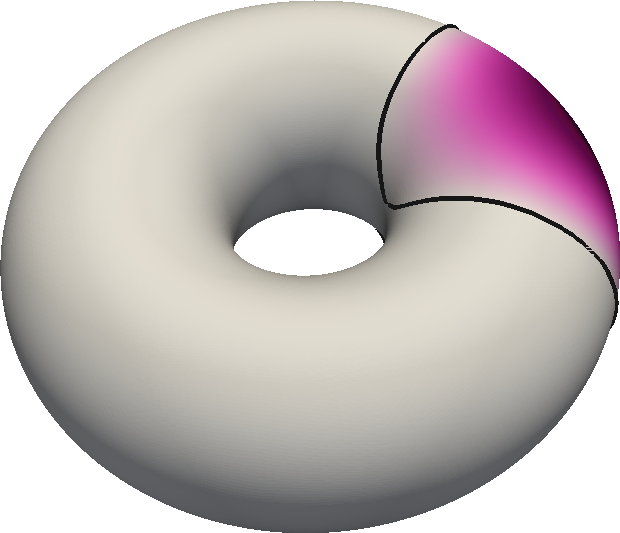} \\ 3.772
    \end{tabular}
    \begin{tabular}{c}
      two sided shape \\ \includegraphics[width=1in,height=0.9in,keepaspectratio]{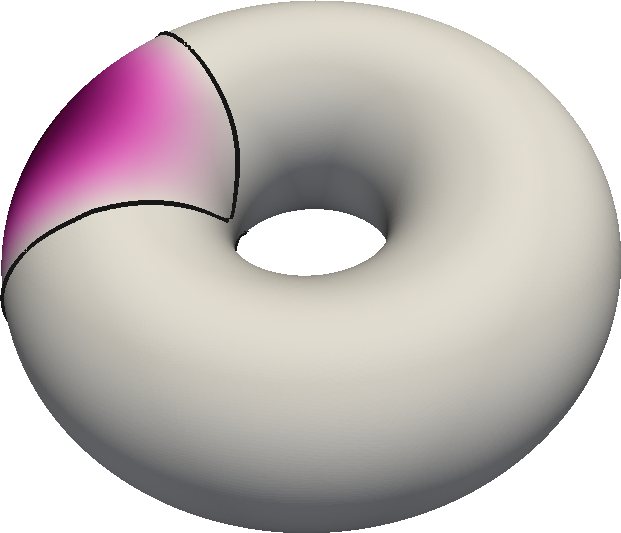} \\ 3.940
    \end{tabular}
    \begin{tabular}{c}
      two sided shape \\ \includegraphics[width=1in,height=0.9in,keepaspectratio]{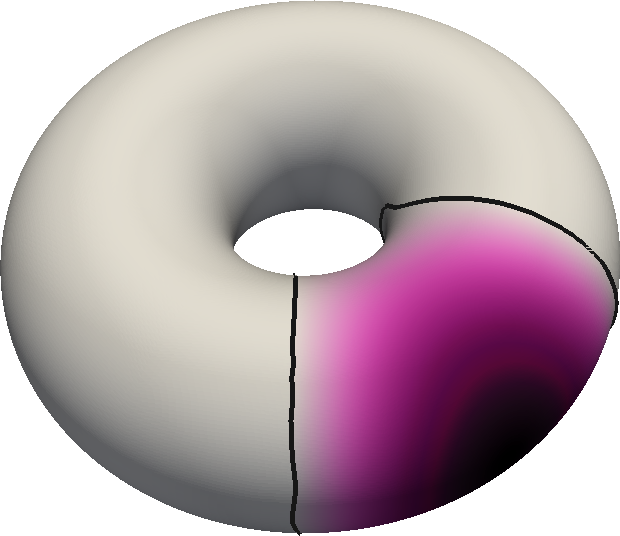} \\ 3.683
    \end{tabular}
    \begin{tabular}{c}
      two sided shape \\ \includegraphics[width=1in,height=0.9in,keepaspectratio]{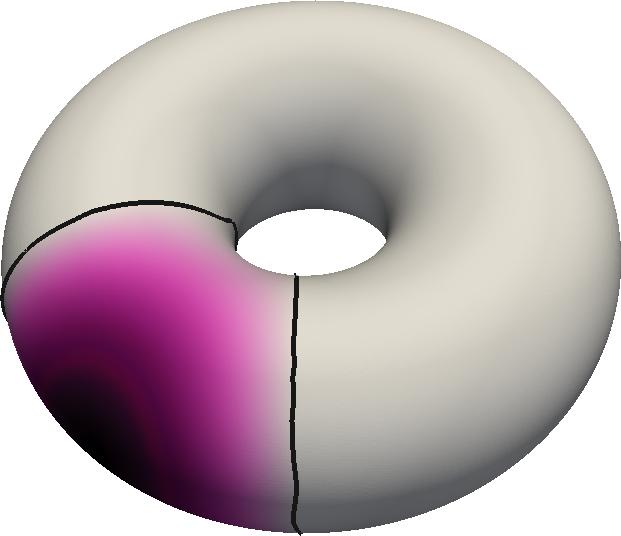} \\ 3.914
    \end{tabular}
    \begin{tabular}{c}
      quadrilateral \\ \includegraphics[width=1in,height=0.9in,keepaspectratio]{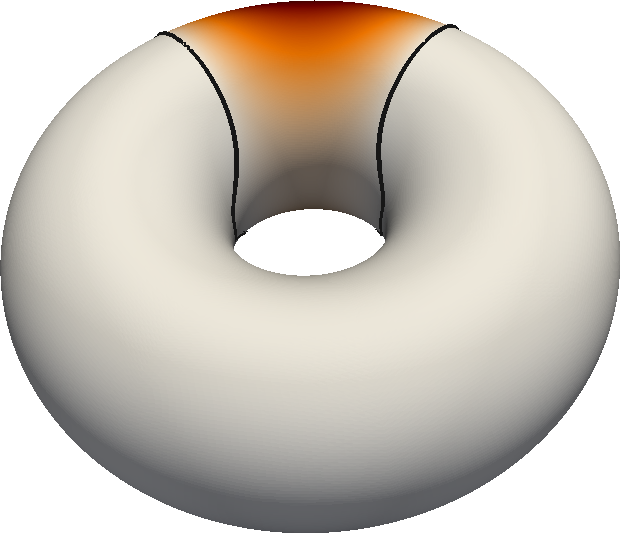} \\ 3.812
    \end{tabular}
    
    $S_\eps$: $0.595$

    Total energy: $19.717$
    \\
    \hline
    6 & %
    \begin{tabular}{c}
      hexagon \\ \includegraphics[width=1in,height=0.9in,keepaspectratio]{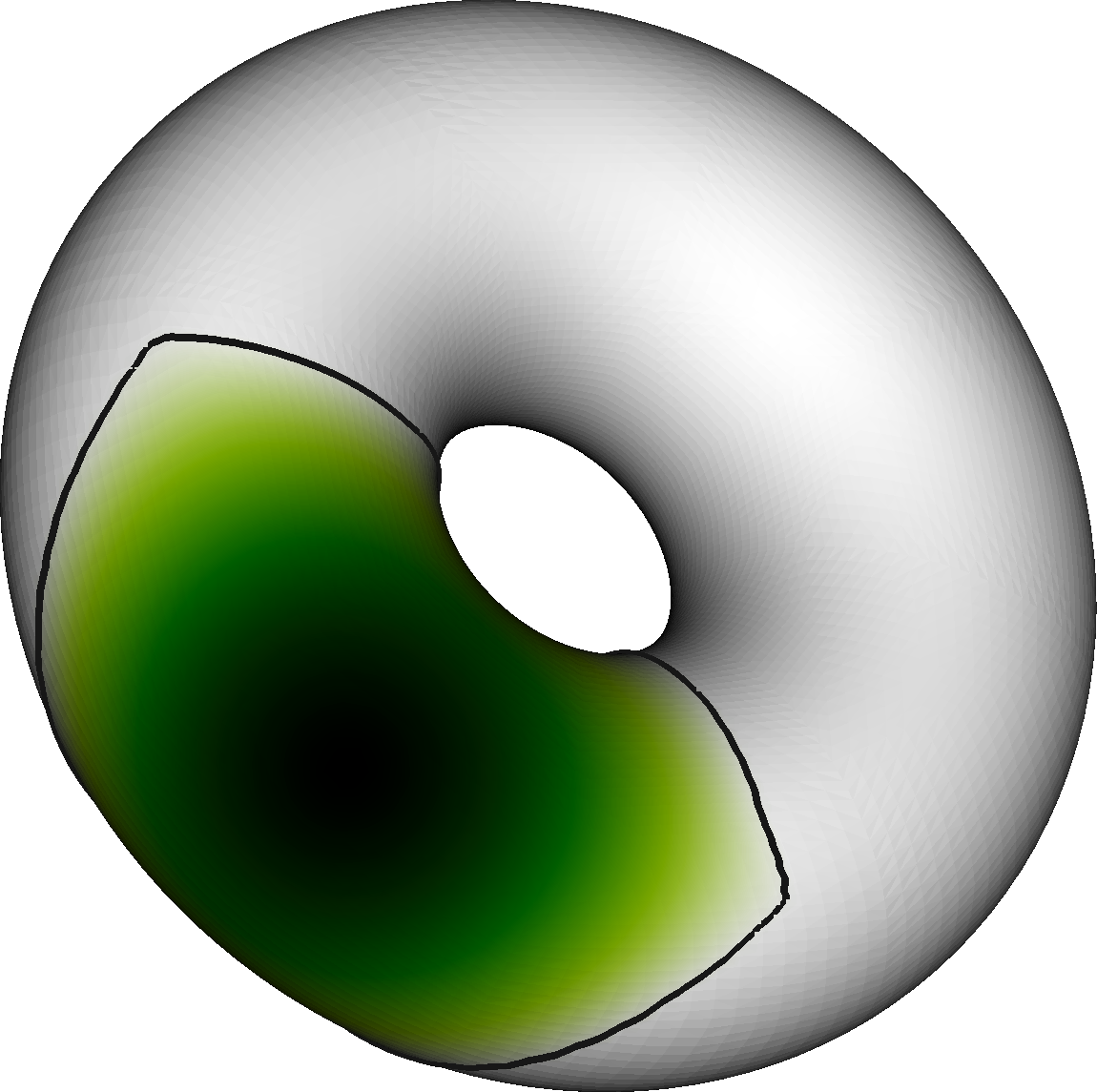} \\ 4.215
    \end{tabular}
    \begin{tabular}{c}
      hexagon \\ \includegraphics[width=1in,height=0.9in,keepaspectratio]{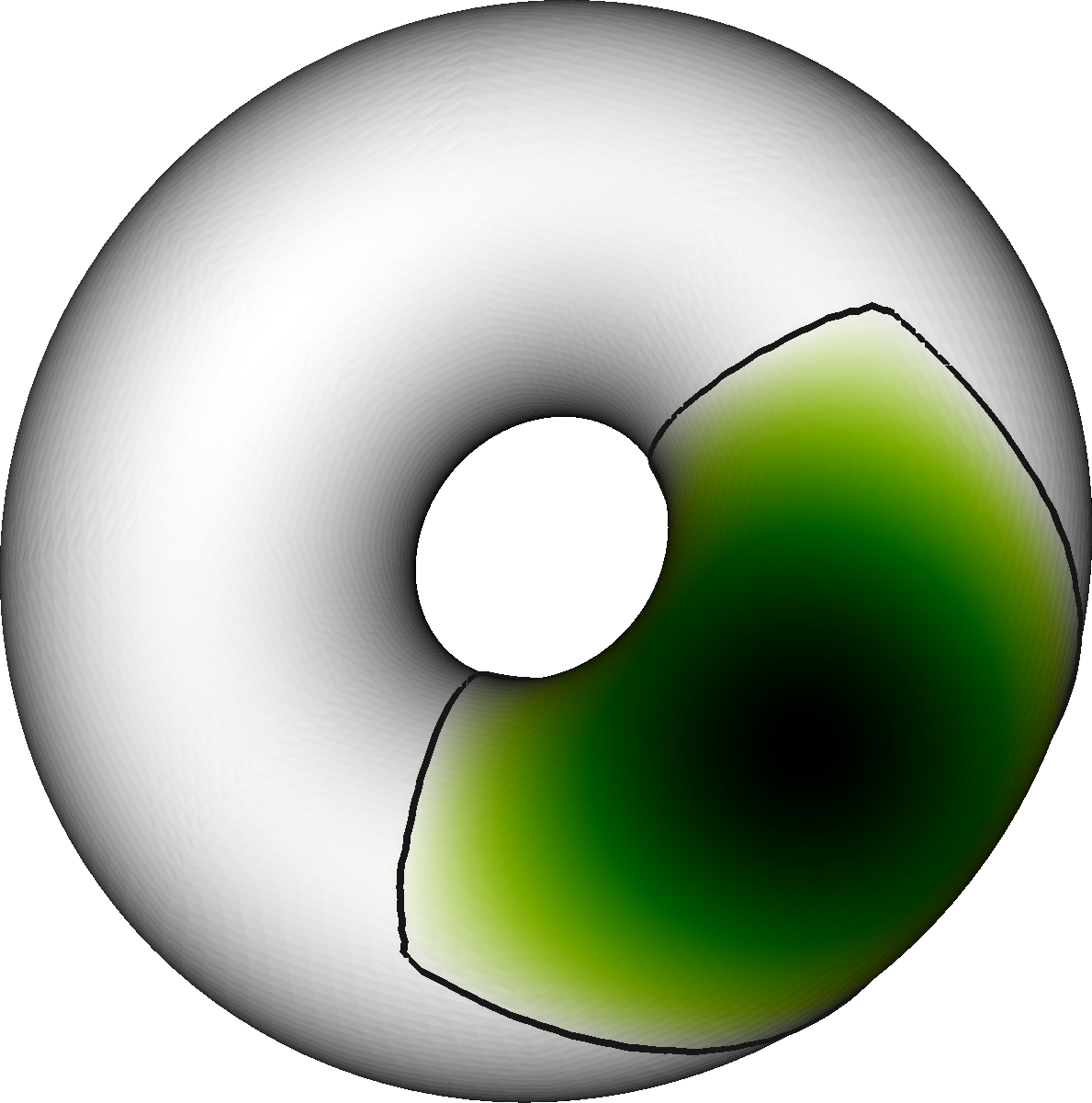} \\ 4.481
    \end{tabular}
    \begin{tabular}{c}
      hexagon \\ \includegraphics[width=1in,height=0.9in,keepaspectratio]{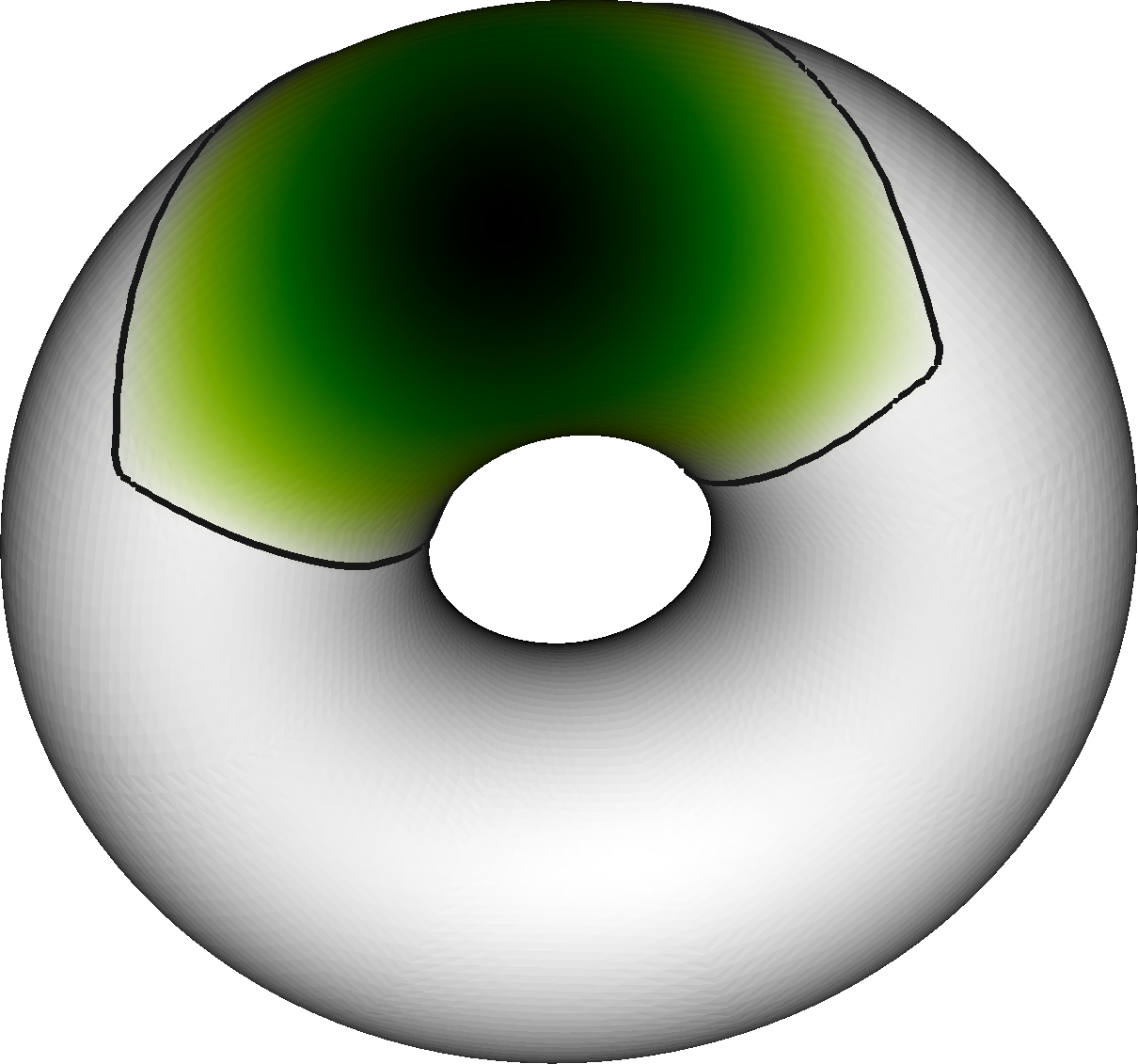} \\ 4.319
    \end{tabular}

    \begin{tabular}{c}
      hexagon \\ \includegraphics[width=1in,height=0.9in,keepaspectratio]{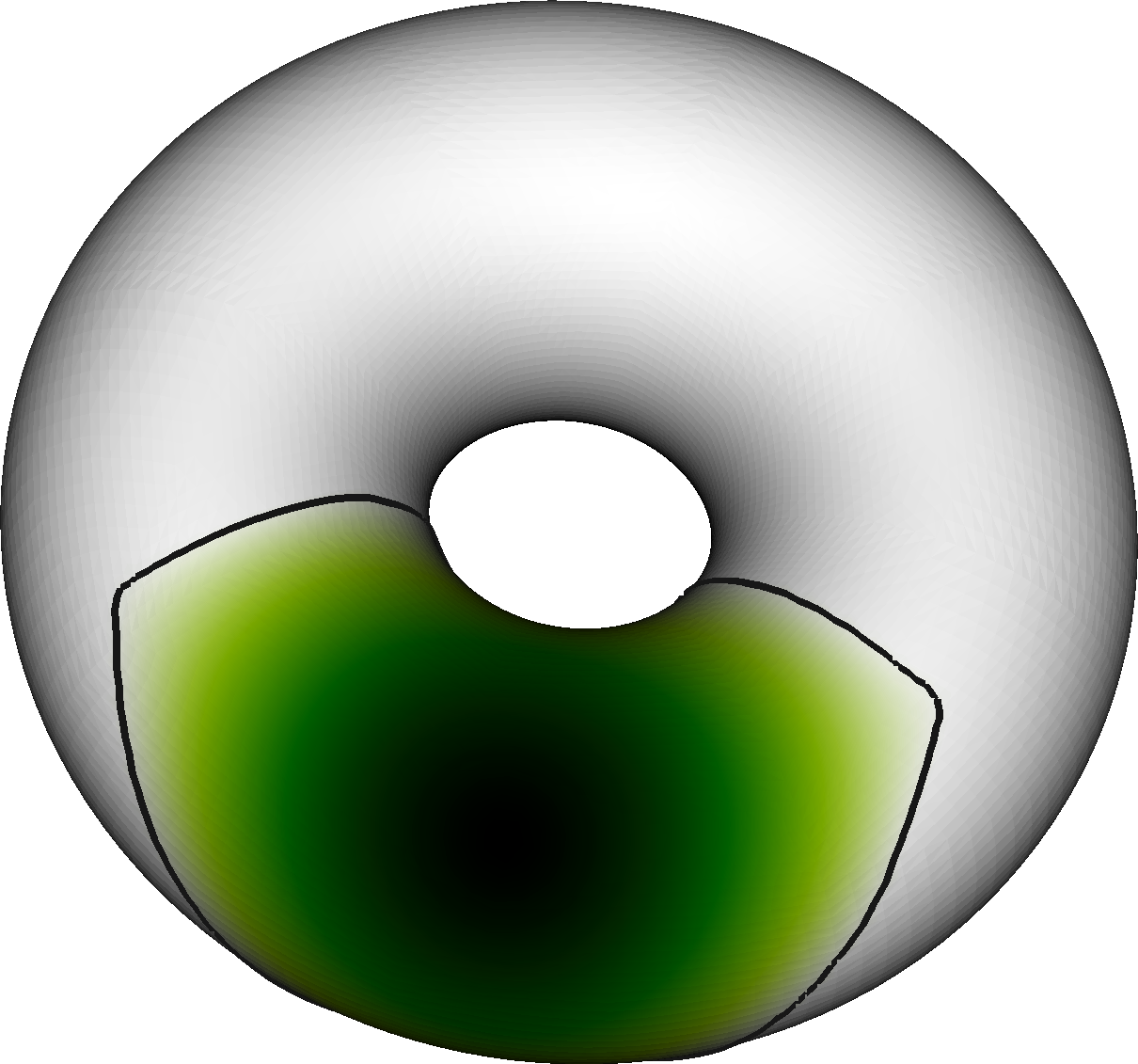} \\ 4.319
    \end{tabular}
    \begin{tabular}{c}
      hexagon \\ \includegraphics[width=1in,height=0.9in,keepaspectratio]{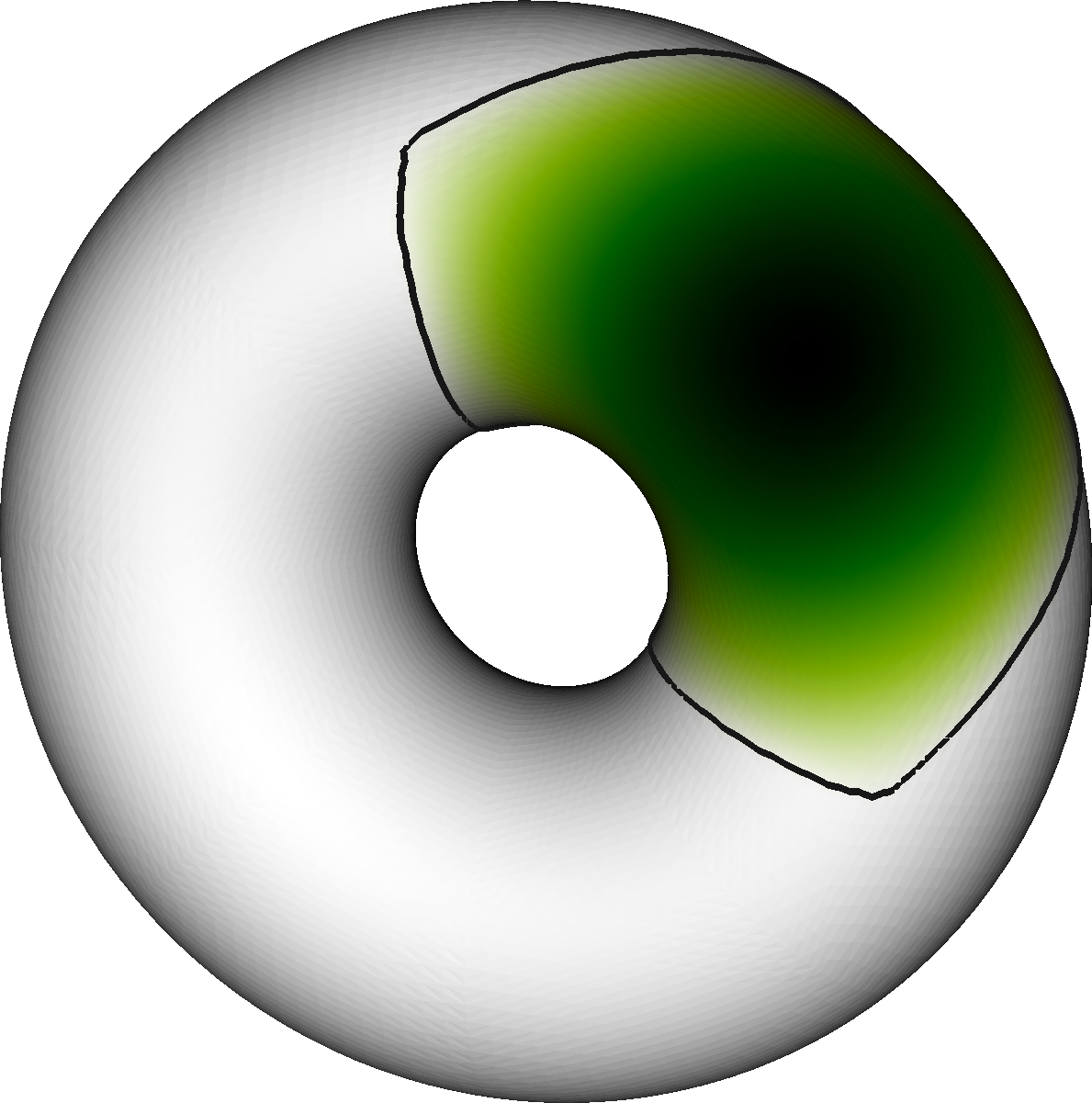} \\ 4.480
    \end{tabular}
    \begin{tabular}{c}
      hexagon \\ \includegraphics[width=1in,height=0.9in,keepaspectratio]{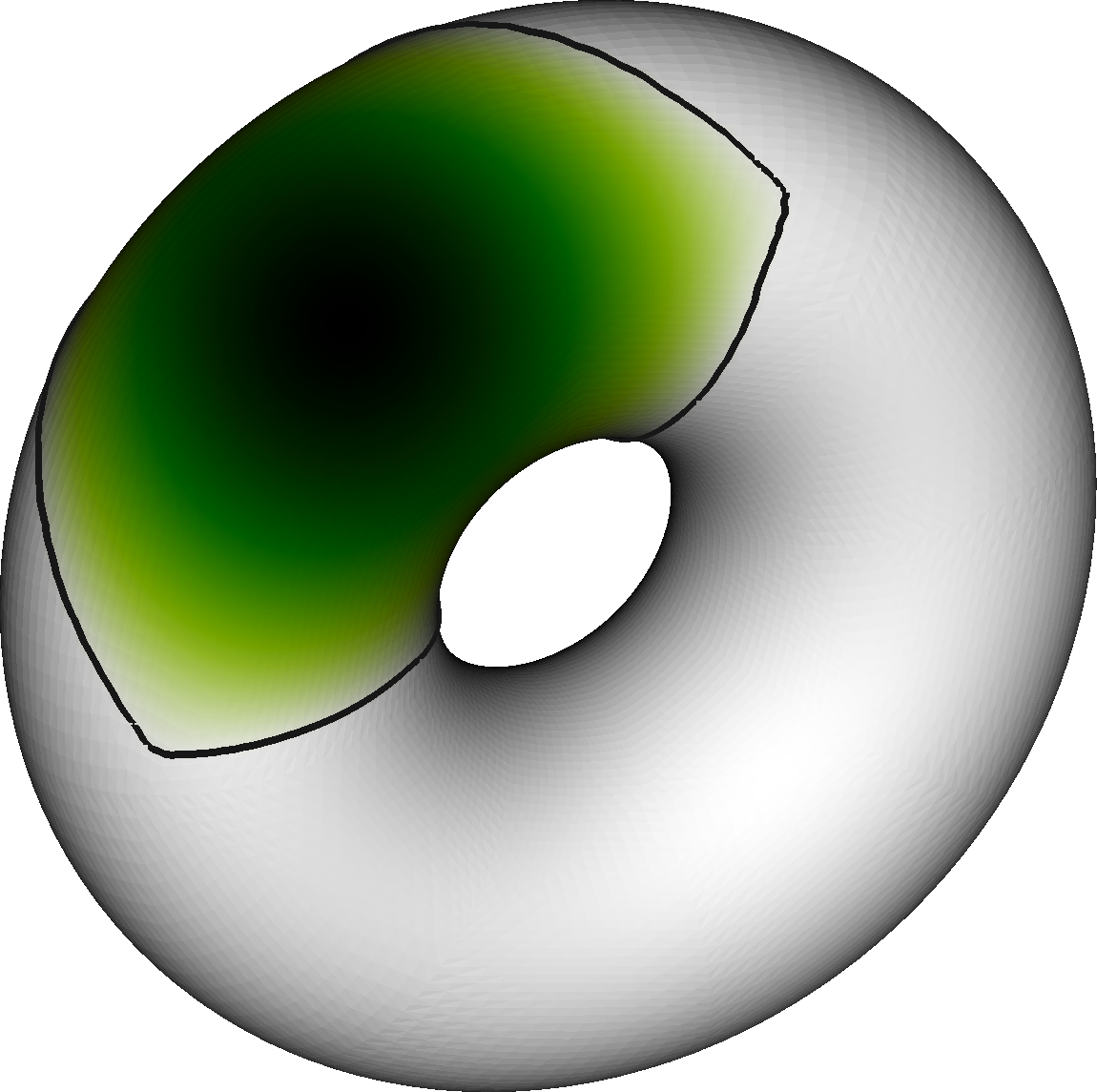} \\ 4.215
    \end{tabular}

    $S_\eps$: $1.005$

    Total energy: $27.035$
    \\
    \hline
  \end{tabular}
\end{table}

\begin{table}[p]
  \centering
  \newcolumntype{C}[1]{>{\centering\let\newline\\\arraybackslash\hspace{0pt}}m{#1}}
  \begin{tabular}{|C{0.025\textwidth}|C{0.9725\textwidth}|}
    \hline
    7 & %
    \begin{tabular}{c}
      quadrilateral \\ \includegraphics[width=1in,height=0.9in,keepaspectratio]{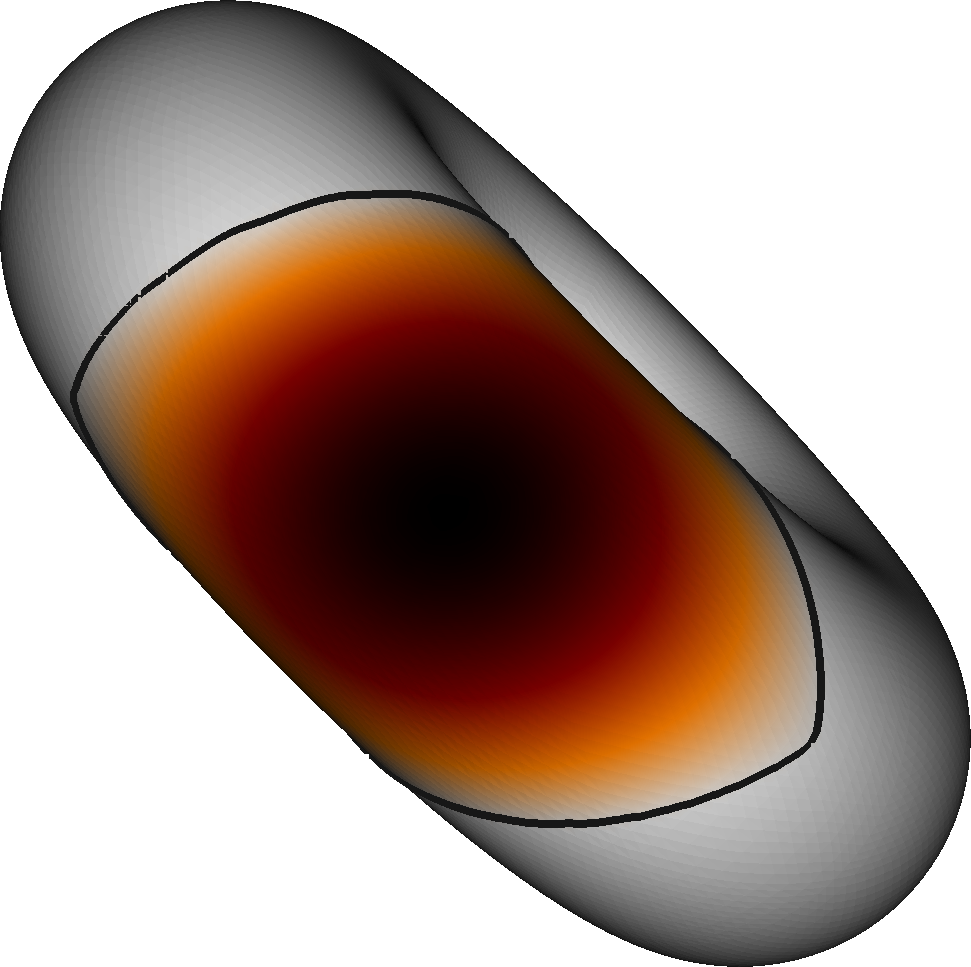} \\ 4.803
    \end{tabular}
    \begin{tabular}{c}
      quadrilateral \\ \includegraphics[width=1in,height=0.9in,keepaspectratio]{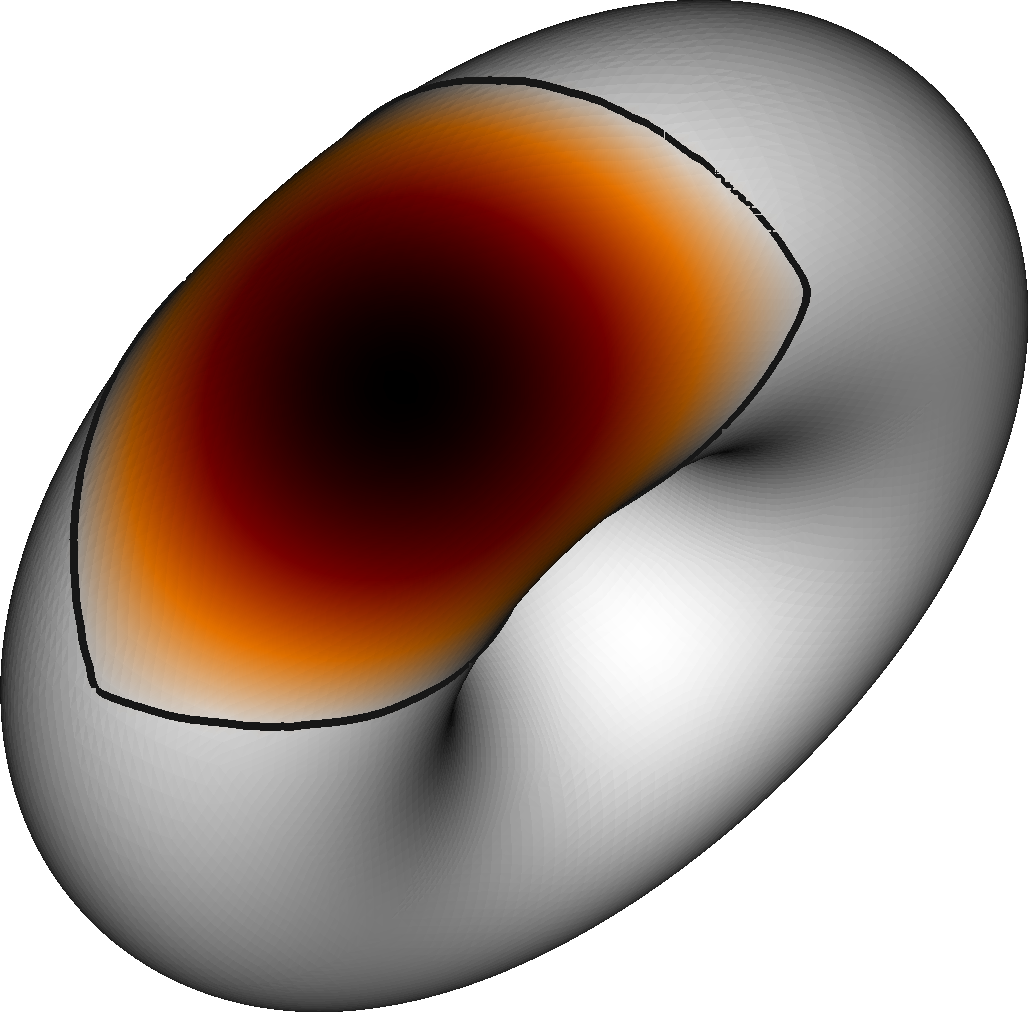} \\ 5.064
    \end{tabular}
    \begin{tabular}{c}
      pentagon \\ \includegraphics[width=1in,height=0.9in,keepaspectratio]{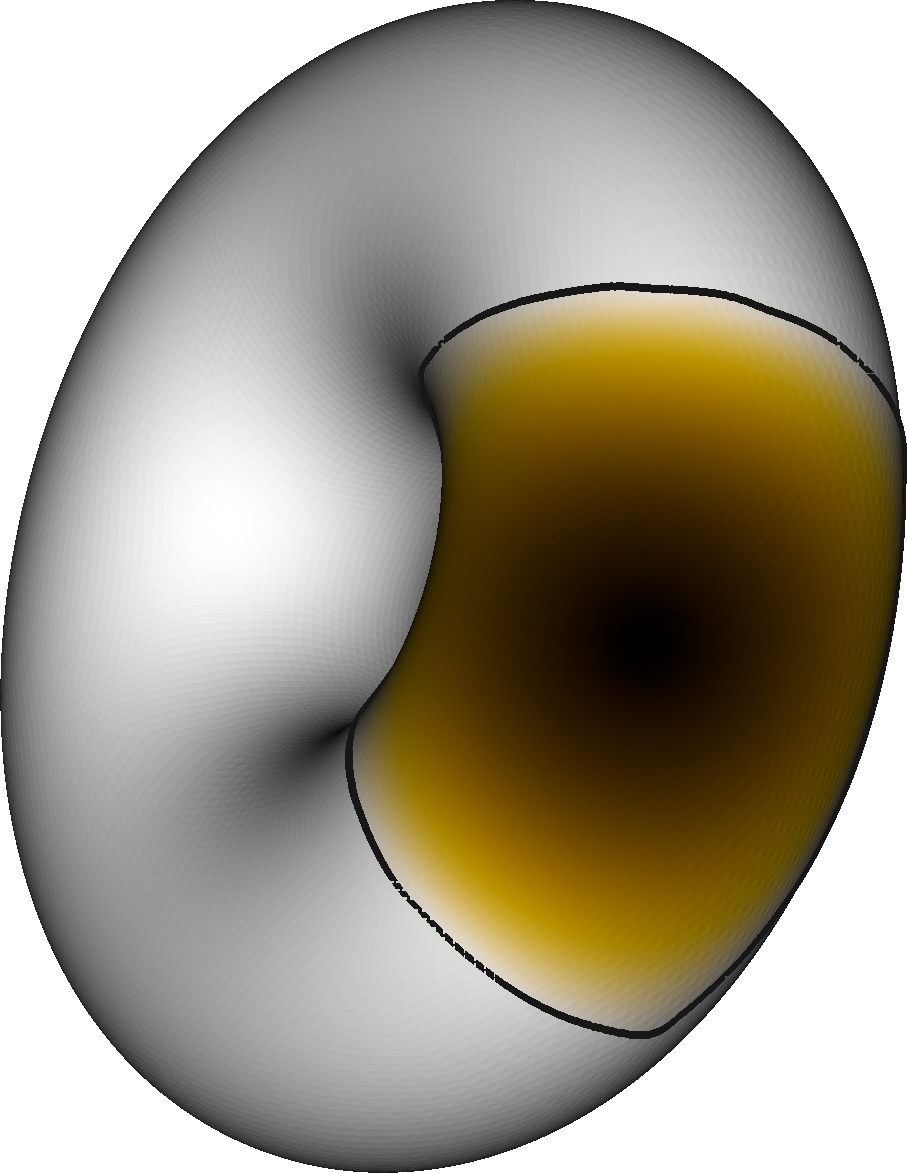} \\ 5.168
    \end{tabular}
    \begin{tabular}{c}
      pentagon \\ \includegraphics[width=1in,height=0.9in,keepaspectratio]{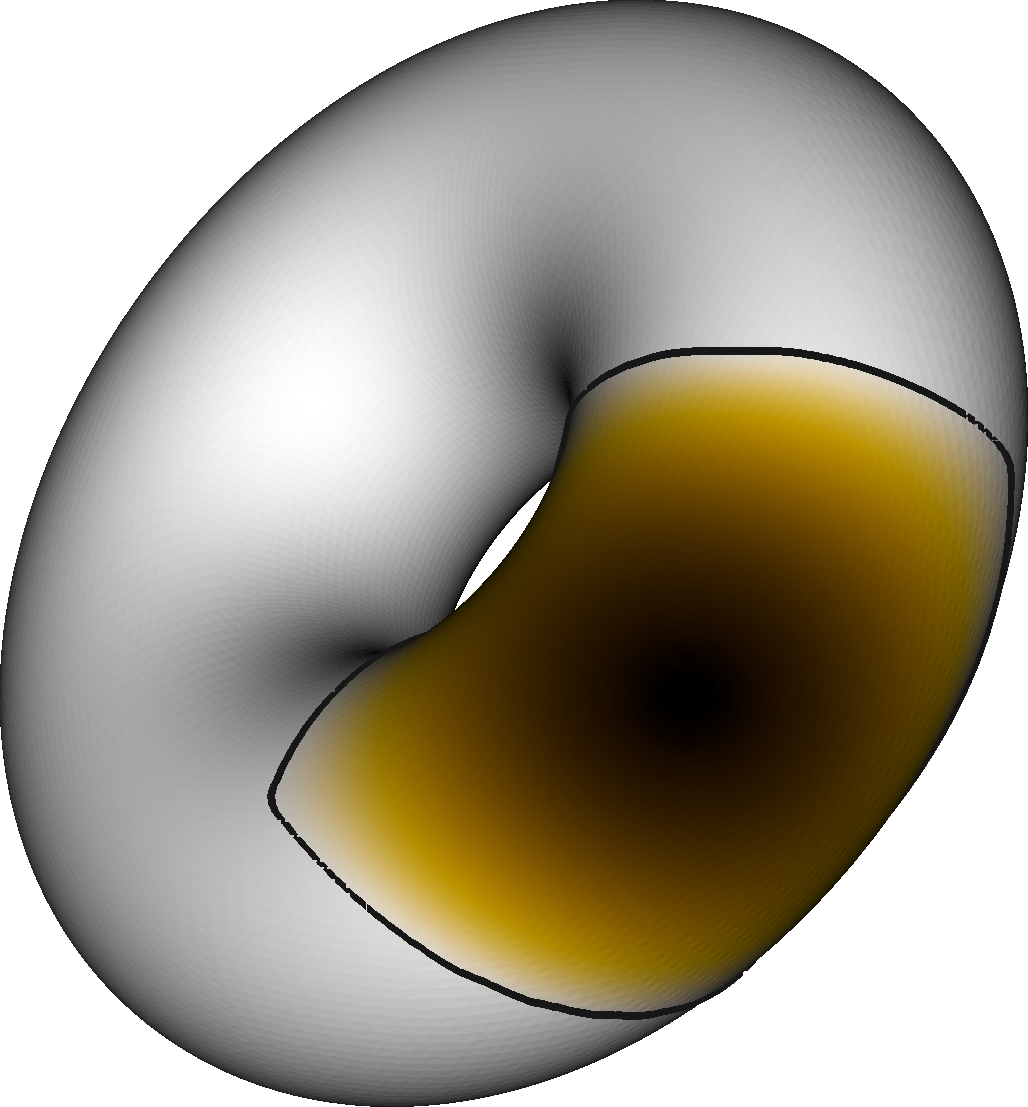} \\ 4.94272
    \end{tabular}

    \begin{tabular}{c}
      hexagon \\ \includegraphics[width=1in,height=0.9in,keepaspectratio]{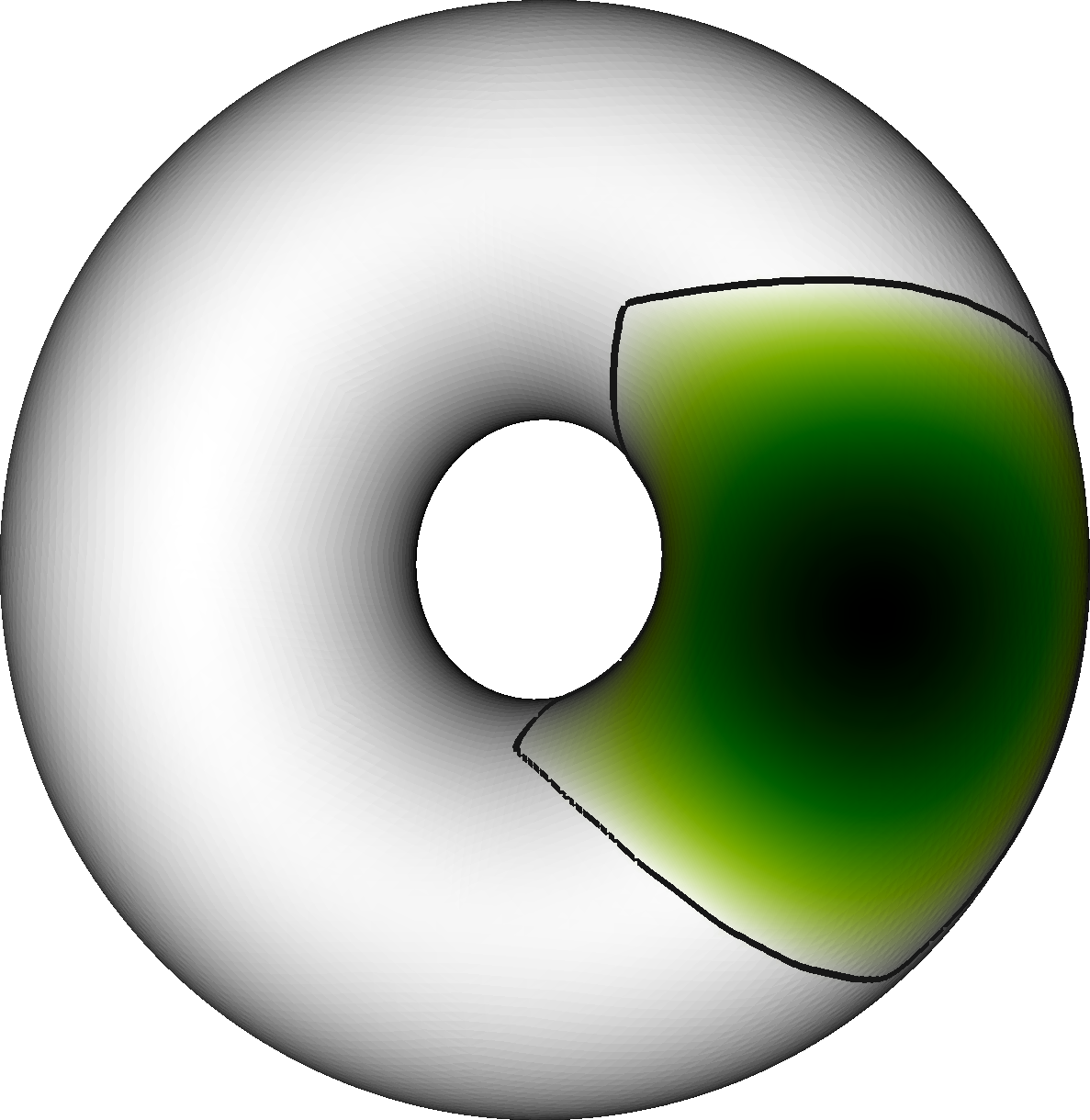} \\ 5.465
    \end{tabular}
    \begin{tabular}{c}
      octagon \\ \includegraphics[width=1in,height=0.9in,keepaspectratio]{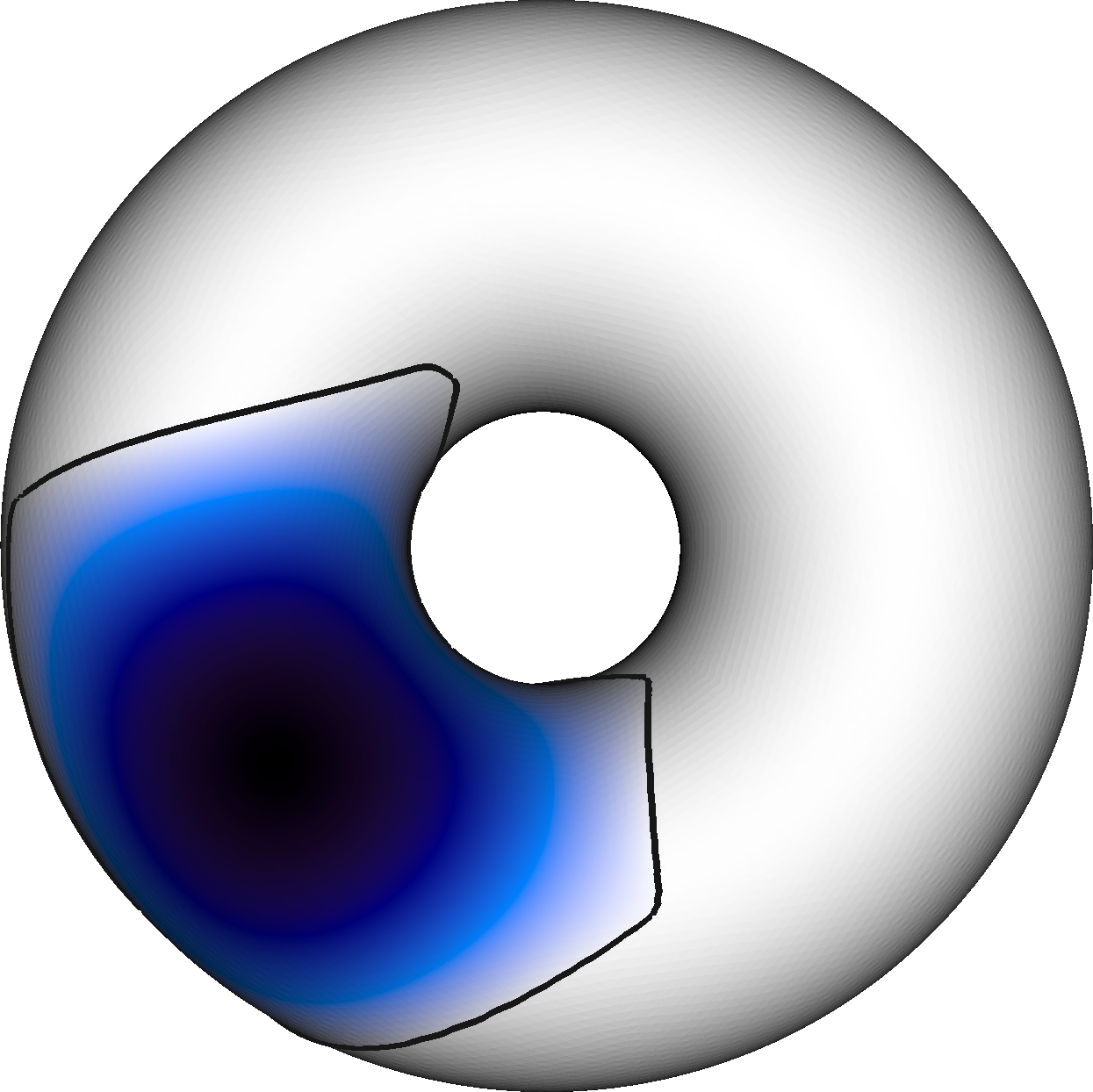} \\ 5.459
    \end{tabular}
    \begin{tabular}{c}
      decagon \\ \includegraphics[width=1in,height=0.9in,keepaspectratio]{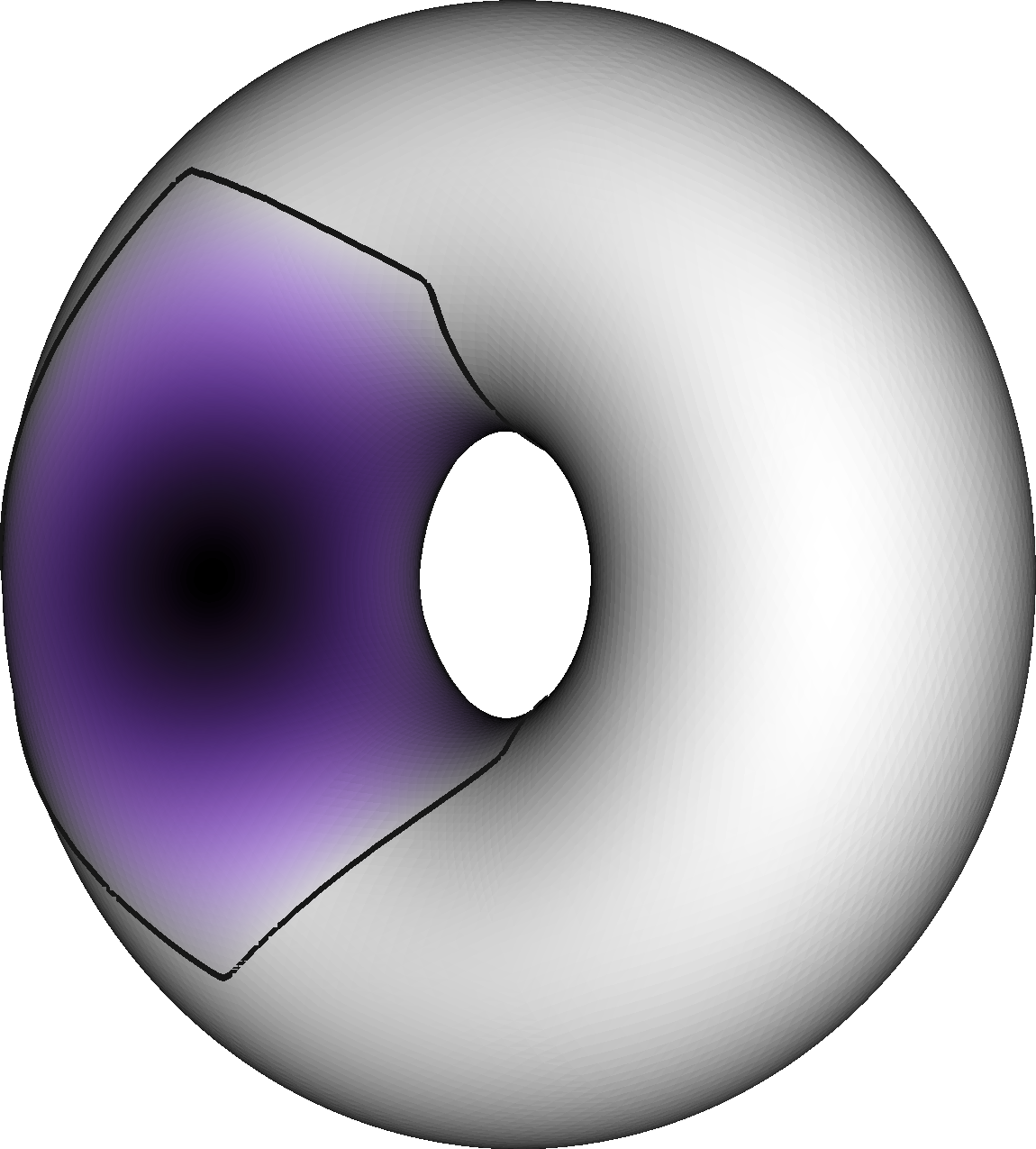} \\ 5.908
    \end{tabular}

    $S_\eps$: $0.81257$

    Total energy: $37.623$
    \\
    \hline
    8 & %
    \begin{tabular}{c}
      pentagon \\ \includegraphics[width=1in,height=0.9in,keepaspectratio]{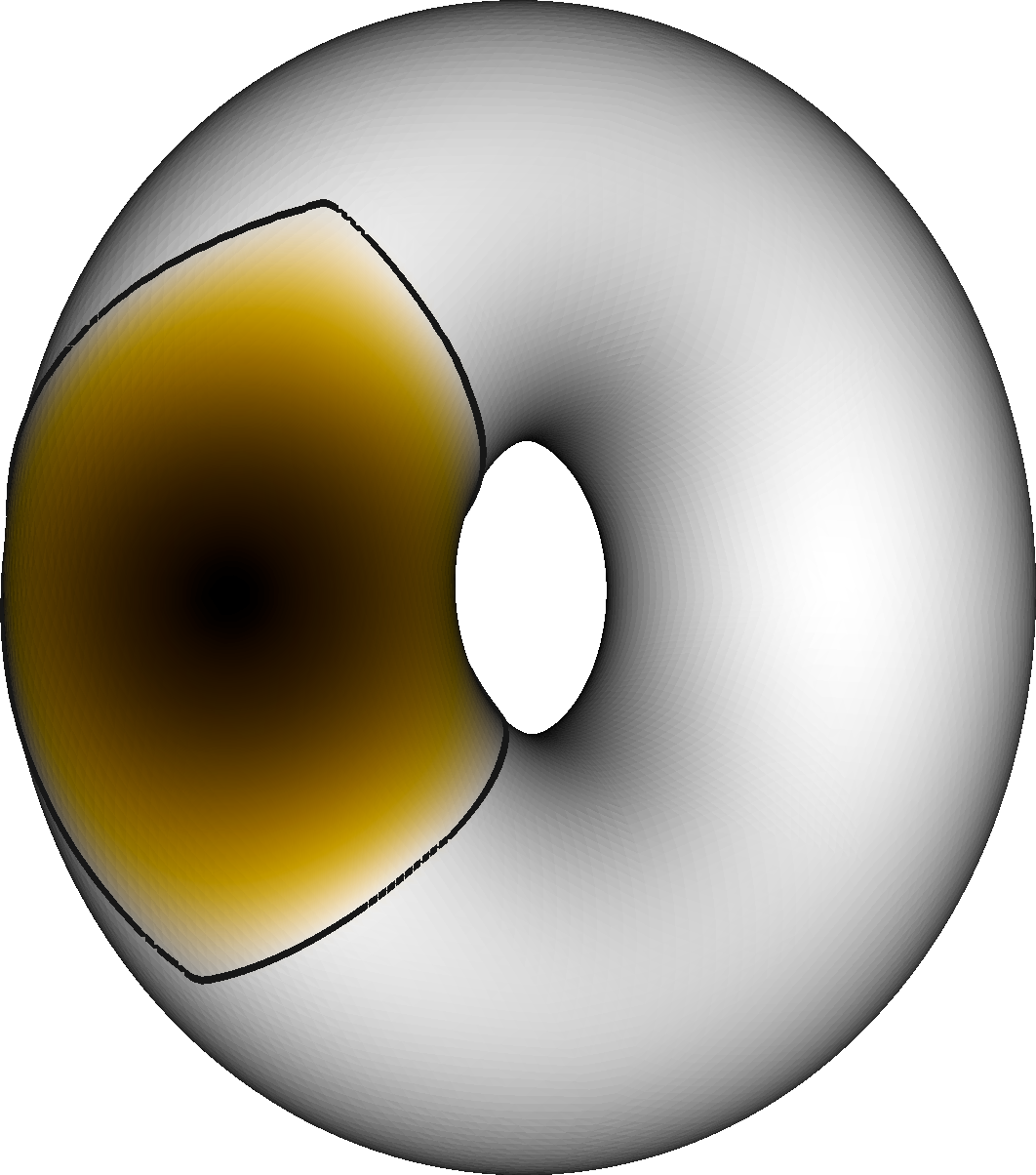} \\ 5.951
    \end{tabular}
    \begin{tabular}{c}
      pentagon \\ \includegraphics[width=1in,height=0.9in,keepaspectratio]{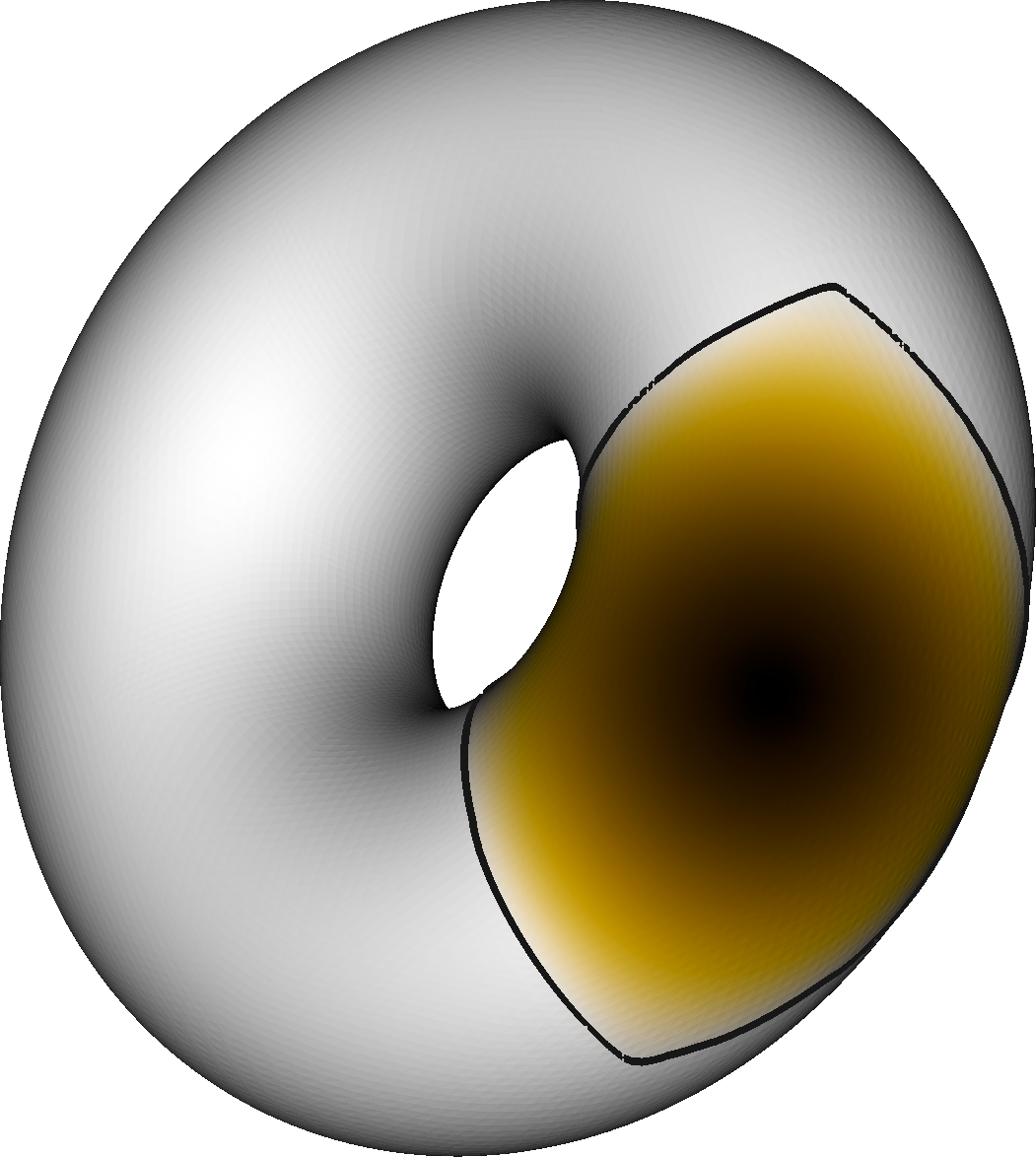} \\ 5.841
    \end{tabular}
    \begin{tabular}{c}
      pentagon \\ \includegraphics[width=1in,height=0.9in,keepaspectratio]{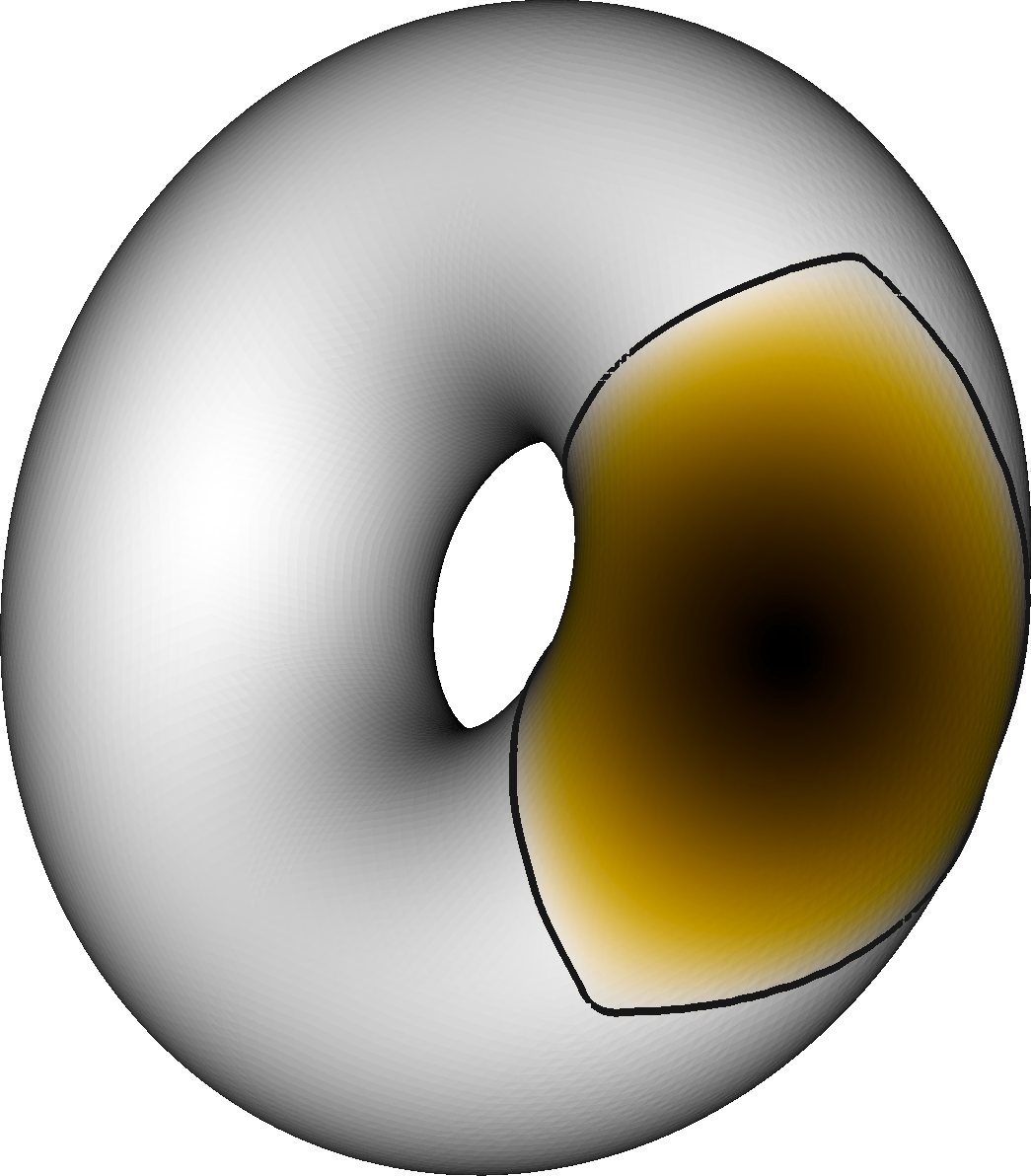} \\ 6.070
    \end{tabular}
    \begin{tabular}{c}
      pentagon \\ \includegraphics[width=1in,height=0.9in,keepaspectratio]{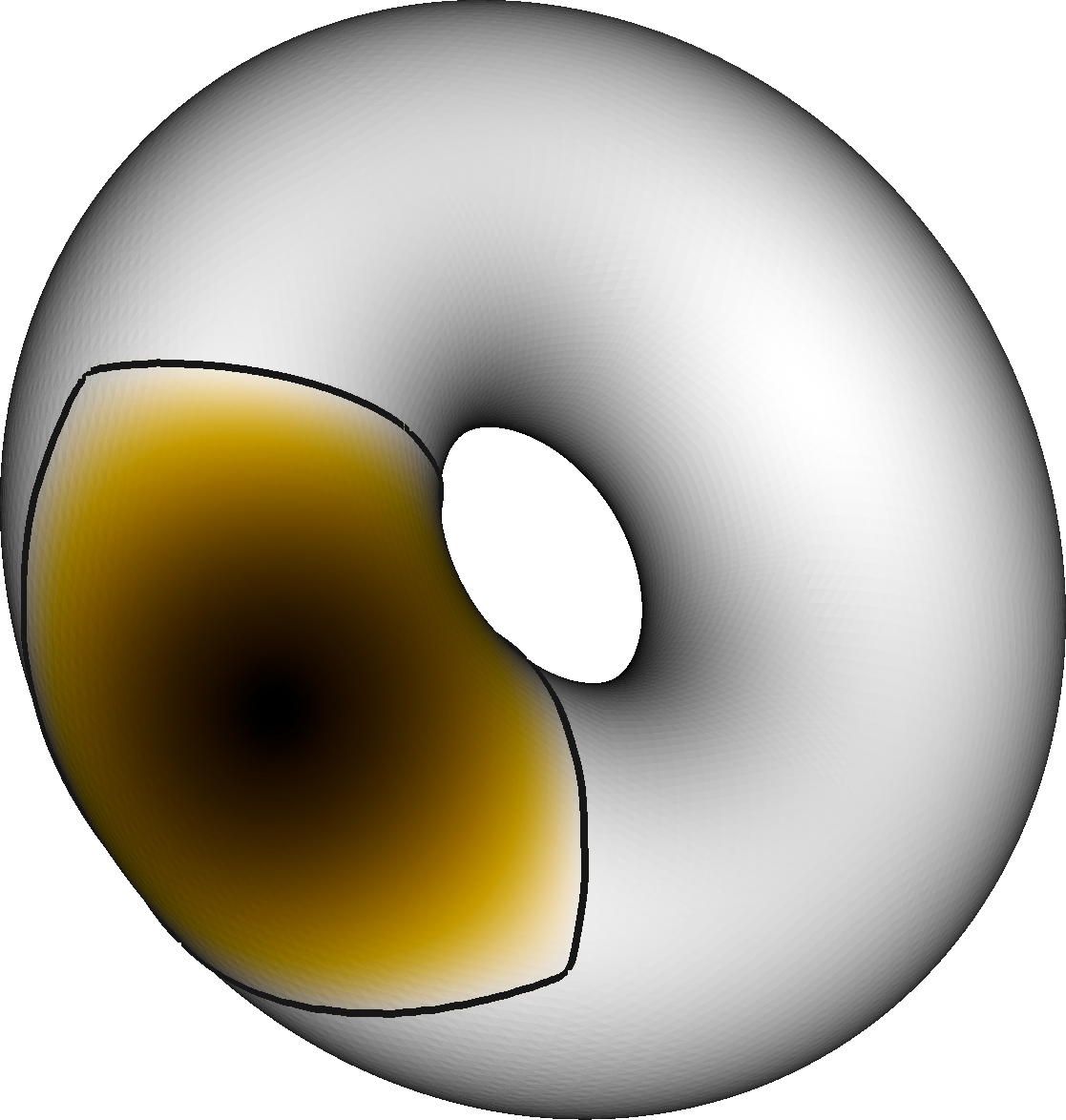} \\ 6.105
    \end{tabular}

    \begin{tabular}{c}
      hexagon \\ \includegraphics[width=1in,height=0.9in,keepaspectratio]{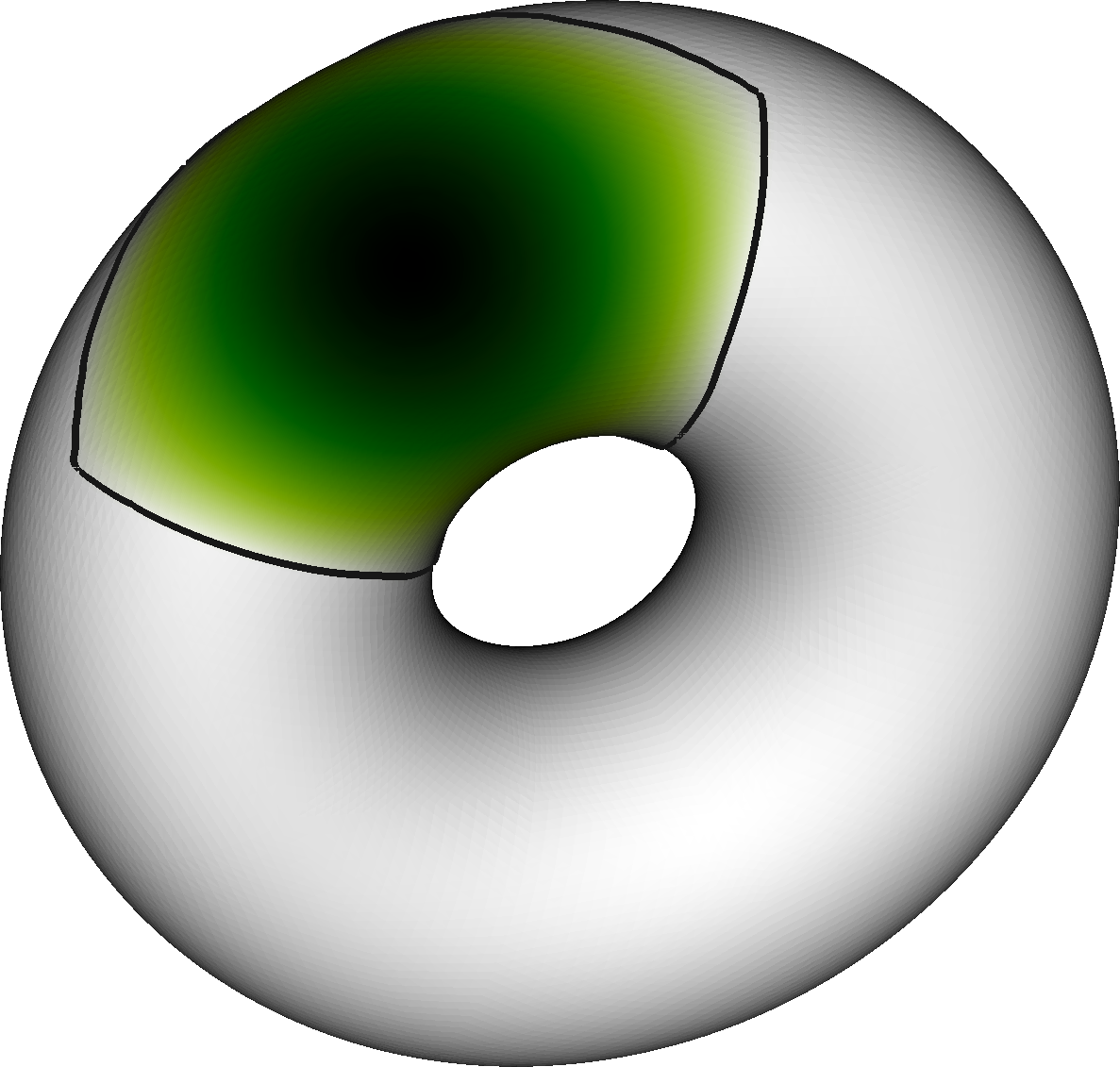} \\ 5.692
    \end{tabular}
    \begin{tabular}{c}
      heptagon \\ \includegraphics[width=1in,height=0.9in,keepaspectratio]{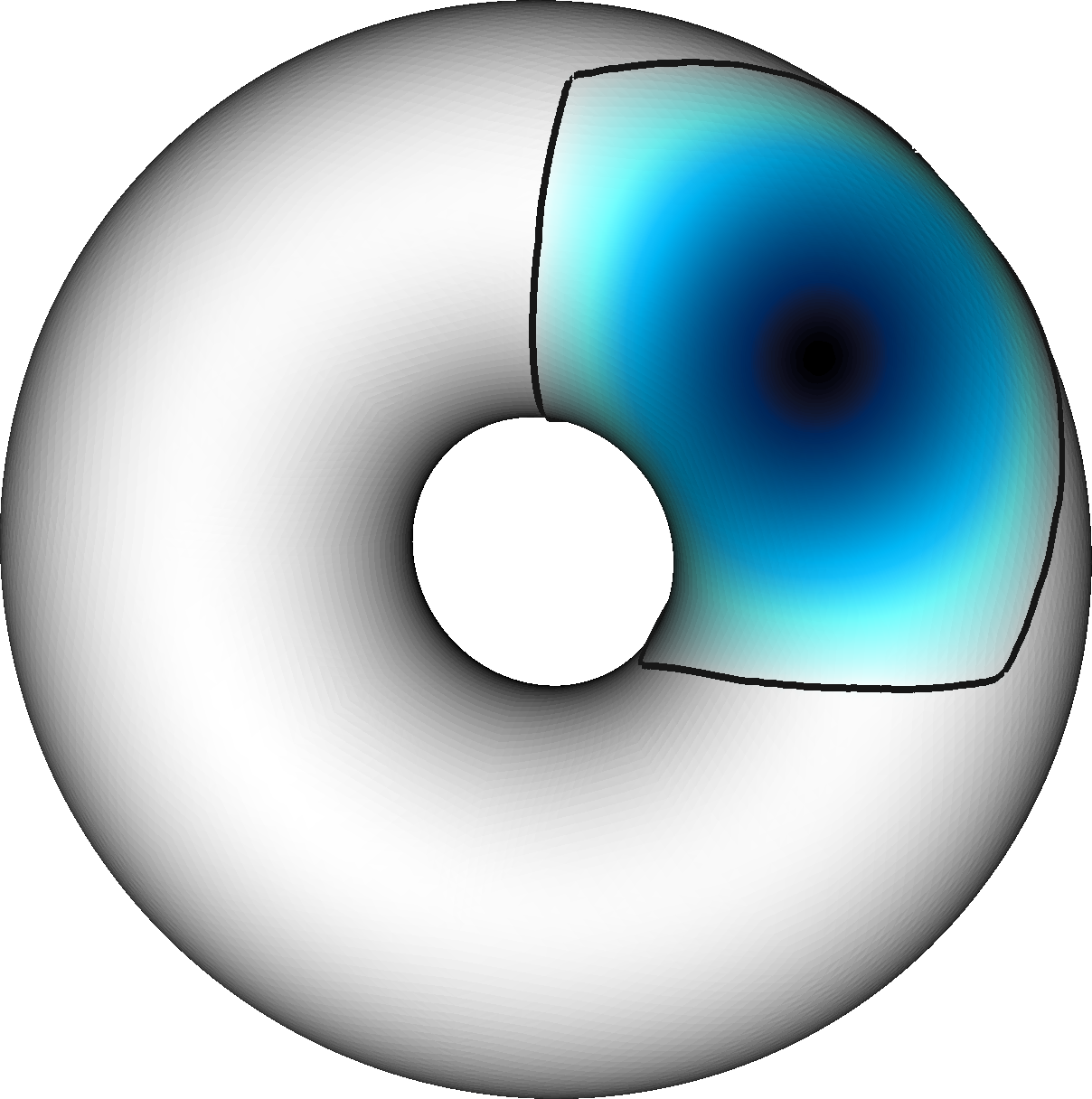} \\ 6.186
    \end{tabular}
    \begin{tabular}{c}
      heptagon \\ \includegraphics[width=1in,height=0.9in,keepaspectratio]{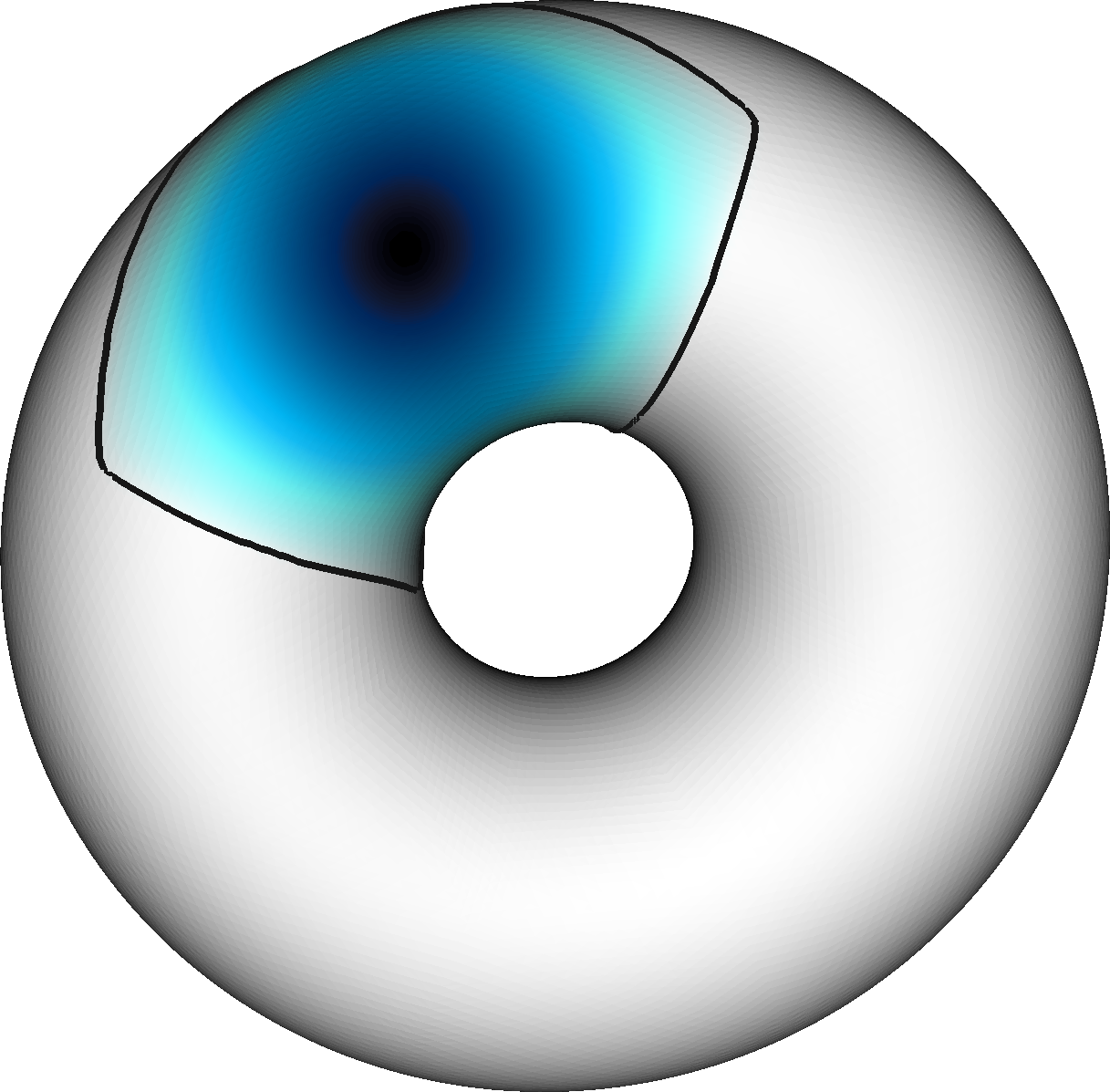} \\ 6.184
    \end{tabular}
    \begin{tabular}{c}
      octagon \\ \includegraphics[width=1in,height=0.9in,keepaspectratio]{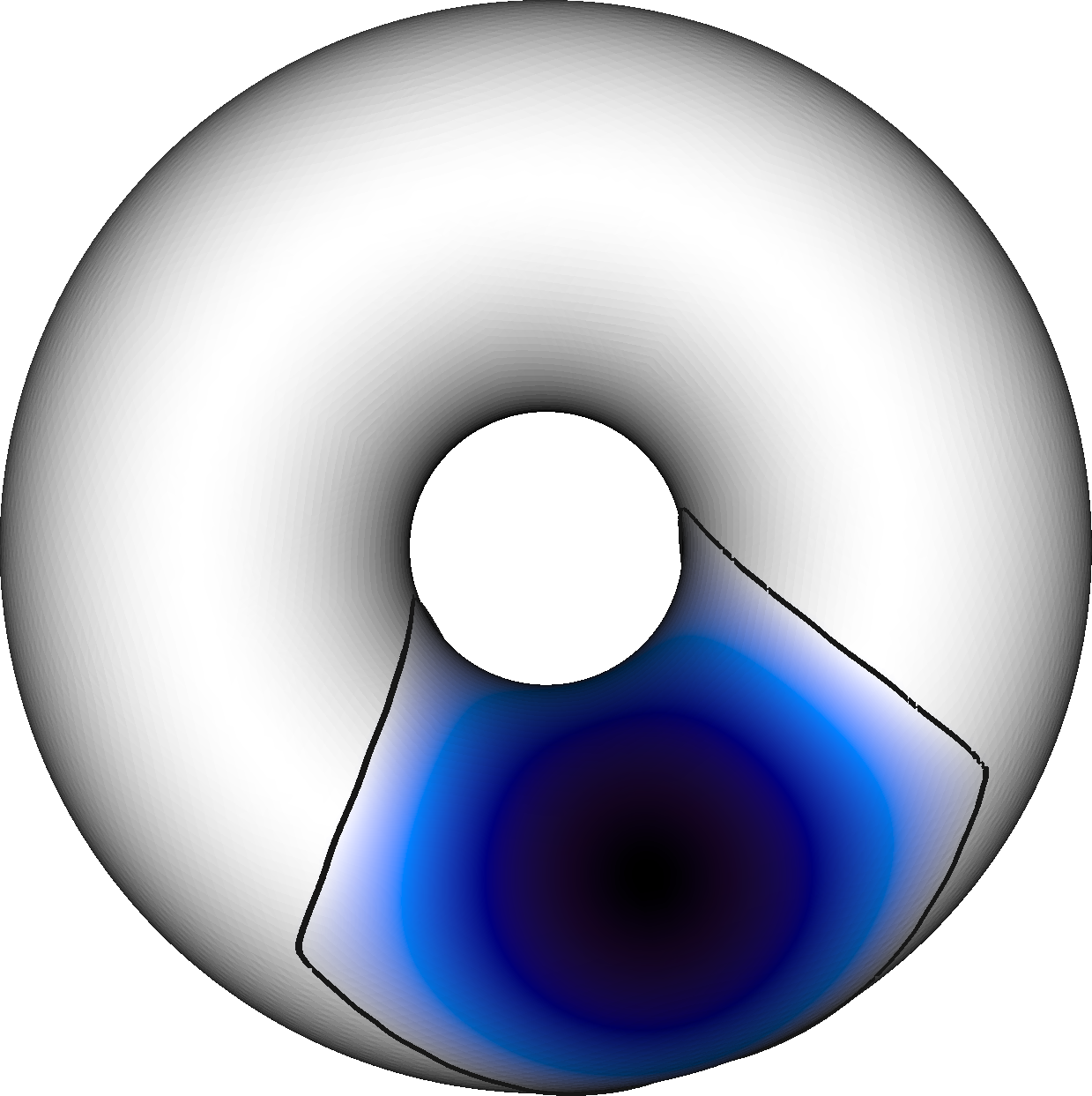} \\ 6.254
    \end{tabular}
  
    $S_\eps$: $1.257$

    Total energy: $49.540$
    \\
    \hline
  \end{tabular}

  \caption{More details of optimal partitions on the torus. In the small plots, we plot the corresponding $u_i^{\eps,h}$ with a black contour at $v_i^{\eps,h} = 0$.}
  \label{tab:many-m-torus}
\end{table}

\begin{figure}
  \centering
  \includegraphics[width=0.8\textwidth]{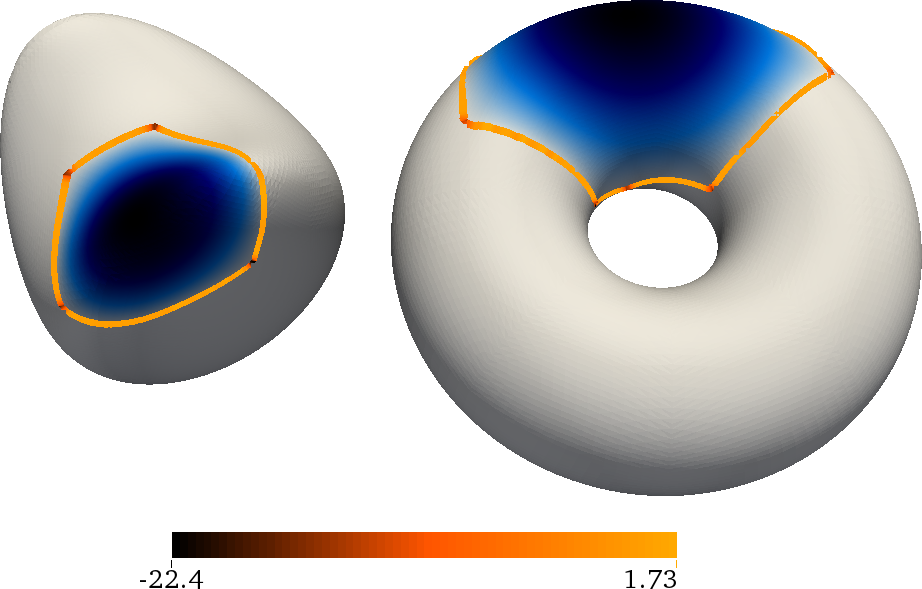}

  \caption{Plots of one partition and $\kappa_g$ for $m=8$ on the surface (\emph{D}) (left) and $m=6$ on the torus (right). The value of $u_i^{\eps,h}$ is shown on a black to white scale and $\kappa_g$ is plotted on the curve $\{ v_i^{\eps,h} = 0 \}$ on a black to orange scale.}
  \label{fig:kappag-other}
\end{figure}

By using $\Gamma_i^{\eps,h}$ and $v_i^{\eps,h}$ from \eqref{eq:vi}, we can define the boundary of partition on these surfaces also. This allows us to compute the geodesic curvature \eqref{eq:kappag} of the boundary of $\Gamma_i^{\eps,h}$; see Figure~\ref{fig:kappag-other} for computations. We again see that away from junctions the geodesic curvature is small. We also see that boundaries all meet at triple junction with the equal angle condition satisfied. We conjecture that on all surfaces optimal partitions have boundaries with zero geodesic curvature which meet at triple junctions with equal angles between each boundary.

On surface ({\it D}), the partition has exactly the same structure as for the sphere for $m \le 8$ but the eigenvalues do not group in the same way because of the variations in curvature. For large values of $m$ the structure changes. Now in regions with higher curvature we see partitions with few sides. In fact, for $m=16$, three partitions have four sides, which does not occur in the case of the sphere. The familiar pattern of pentagons and hexagons reoccurs for $m=32$ except now the pentagons are clustered in regions of high curvature. The number of sides of each partition is still limited to six. Because of \eqref{eq:gauss-bonnet}, for larger values of $m$ we expect to see $12$ pentagons and $m - 12$ hexagons with the pentagons clustered in the higher curvature regions.
We see that none of the partitions are equi-spectral.

On the torus ({\it T}),  the situation is very different. For $m \le 6$, we have very structured partitions which reflect the symmetry of the surface. For the case of $m=5$, we see all triple junctions occur in the center of the torus. For $m > 6$, we have partitions with more that $6$ sides. The formula \eqref{eq:gauss-bonnet} tells us that the numbers of partitions with more than six sides must balance the number of partitions with less than six sides. For the cases we see, the partitions with more than six sides cluster in the center and those with less than six sides cluster on the exterior.
As we increase $m$ we see an increase in the number of hexagons, however it is not clear whether the number of non-hexagonal partitions will decrease.
For smaller area partitions, for larger $m$, the curvature of the surface is less important and the problem becomes more like the flat problem, so we expect that for large values of $m$, we will  see a preponderance of hexagons.
We see that the partitions for $m = 3,4$ are almost equi-spectral and so conjecture that these partitions are also optimal for the case $p=\infty$.

%%% Local Variables: 
%%% mode: latex
%%% TeX-master: "../partition"
%%% End: 

\section{Discussion}

We have explored an eigenvalue partition problem on three different surfaces and for many different numbers of partitions. We have observed good convergence both with respect to discretisation parameters and also with respect to our choice of regularisation. From our results we make the following conjectures:

\begin{enumerate}
\item The optimal partition consists of curvilinear polygons whose edges have zero geodesic curvature.
\item Partitions either meet along edges or at triple junctions where edges meet at equal angles.
\item For genus zero surfaces, for large values of $m$ the optimal partition consists of $12$ pentagons and $m-12$ hexagons. If the curvature of the surface varies, the pentagons will be located where the curvature is highest.
\item For genus one surfaces, for large values of $m$ the optimal partition has a preponderance of  hexagons.
\end{enumerate}

%%% Local Variables:
%%% mode: latex
%%% TeX-master: "../partition"
%%% End:

\subsection*{Acknowledgments}

The research of TR was funded by the EPSRC (grant number EP/L504993/1).
This work was undertaken on ARC2, part of the High Performance Computing facilities at the University of Leeds. 
The authors were participants of the Isaac Newton Institute programme Free Boundary Problems and Related Topics (January--July 2014) when this article was written.

%%% Local Variables: 
%%% mode: latex
%%% TeX-master: "../partition"
%%% End: 

%% Use the widest label as parameter.

%% Reference items may be numbered or have labels of your choice.
%% The author's surname PRECEDES the initial of the first name.
%% The surnames are set in small caps.
%% Add \& before the last author's surname.
%% Book titles and journal names are italicized. 
%% Only journal volume numbers are boldfaced.
%% The issue number is only given when the issues are paginated separately.

%%%%%%%%%%% To ease editing, use normal size:

\normalsize
\baselineskip=17pt

%%%%%%%%%%%%%

\end{document}